\documentclass{elsarticle}

\usepackage{bookmark}
\usepackage{graphicx}
\usepackage{url}
\usepackage{amsmath,amsfonts}  
\usepackage{amssymb}
\usepackage{mathrsfs}
\usepackage{geometry}          
\usepackage{setspace}
\usepackage{hhline}
\usepackage{times}
\usepackage{color}
\usepackage{url}
\usepackage{fancyhdr}          
\usepackage{chapterbib}
\usepackage{supertabular}
\usepackage{array}
\usepackage{tabularx,ragged2e,booktabs,caption}
\usepackage{tikz}
\usetikzlibrary{shapes.geometric,arrows}
\usetikzlibrary{patterns}
\usepackage{subfig}
\usepackage{float}
\usepackage{multirow}
\usepackage{xfrac}
\usepackage{xcolor}

\newcommand{\bx}{\mathbf{x}}

\newcommand{\bu}{\mathbf{u}}
\newcommand{\BU}{\mathbf{U}}
\newcommand{\bw}{\mathbf{w}}
\newcommand{\BW}{\mathbf{W}}

\newcommand{\bv}{\mathbf{v}}

\newcommand{\nl}{\mathscr{L}}

\newcolumntype{C}[1]{>{\centering\let\newline\\\arraybackslash\hspace{0pt}}m{#1}}

\begin{document}

\begin{frontmatter}


\title{PetIGA-MF: a multi-field high-performance toolbox for structure-preserving B-splines spaces}

\author[numpor,amcs]{A.F.~Sarmiento\corref{cor1}}
\ead{adel.sarmientorodriguez@kaust.edu.sa}
\cortext[cor1]{Correponcence to: A.F.~Sarmiento, 4700 King Abdullah University of Science and Technology, al-Khawarizmi Bldg (Bldg 1), Office 4319WS12, Thuwal 23955-6900, Kingdom of Saudi Arabia}
\author[numpor]{A.M.A.~C\^{o}rtes} 
\ead{adrimacortes@gmail.com}
\author[bcam]{D.A.~Garcia}
\ead{dgarcia@bcamath.org}
\author[numpor,conicet]{L.~Dalcin}
\ead{dalcinl@gmail.com}
\author[ornl]{N.~Collier}
\ead{nathaniel.collier@gmail.com}
\author[numpor,amcs]{V.M.~Calo}
\ead{vmcalo@gmail.com}

\address[numpor]{ %
Numerical Porous Media Center (NumPor)\\
King Abdullah University of Science and Technology (KAUST)\\
Thuwal, Saudi Arabia}
\address[amcs]{ %
Applied Mathematics \& Computational Science (AMCS)\\
King Abdullah University of Science and Technology (KAUST)\\
Thuwal, Saudi Arabia}
\address[conicet]{ %
Centro de Investigaci\'on de M\'etodos Computacionales (CIMEC)\\
Consejo Nacional de Investigaciones Cient\'{\i}ficas y T\'ecnicas (CONICET)\\
Universidad Nacional del Litoral (UNL)\\
Santa Fe, Argentina}
\address[bcam]{ %
Computational Mathematics,\\
Basque Center for Applied Mathematics (BCAM)\\
Bilbao, Spain}
\address[ornl]{%
Computer Science and Mathematics Division\\
Oak Ridge National Laboratory\\
Oak Ridge, TN, USA}

\begin{abstract}
We describe the development of a high-performance solution framework for isogeometric discrete differential forms based on B-splines: PetIGA-MF. Built on top of PetIGA, PetIGA-MF is a general multi-field discretization tool. To test the capabilities of our implementation, we solve different viscous flow problems such as Darcy, Stokes, Brinkman, and Navier-Stokes equations. Several convergence benchmarks based on manufactured solutions are presented assuring optimal convergence rates of the approximations, showing the accuracy and robustness of our solver.
\end{abstract}

\begin{keyword}
isogeometric analysis, 
discrete differential forms, 
structure-preserving discrete spaces, 
multi-field discretizations, 
PetIGA, 
high-performance computing
\end{keyword}

\end{frontmatter}

\section{Introduction}
The theory of finite element exterior calculus and the underlying concept of discrete differential forms surveyed in \cite{ExtCalc} formalize the design of \textit{compatible} discrete schemes. By compatibility is meant that the discretization preserves the mathematical structure of the partial differential equation and the functional spaces underlying them, from the continuous to the discrete setting. One example where such compatibility property is a necessary requirement for the stability of the discrete scheme is the Maxwell equations system. In many cases, such compatibility condition is encoded on the commutativity of the de Rham diagram \cite{monk2003finite,Dem2000}.

Isogeometric Analysis (IGA)\cite{HBCBOOK} allows the definition of a family of discrete differential forms, based on splines functions, called isogeometric discrete differential forms. The isogeometric discrete differential forms theory, described in \cite{Buffa2011}, provides structure-preserving discrete spaces, namely, the gradient-, curl-, divergence- and integral-conforming spaces, which satisfy a discrete de Rham diagram. Curl-conforming spaces were first applied to approximate and solve Maxwell's equations \cite{Buffa2010}. Later the divergence- and integral-conforming spaces were used to approximate and solve the Stokes system in \cite{Buffa2010stokes}. In a series of papers \cite{EvansHughesStokes,EvansHughesNavierStokes,EvansNSUnsteady} Evans and Hughes further developed the theory and the application of these spaces to approximate different incompressible viscous flow problems such as Darcy, Stokes, Brinkman and Navier-Stokes equations. In this case, the compatibility of the divergence- and integral-conforming B-spline spaces, when used as a discrete velocity-pressure pair, engenders to important properties of the scheme, namely, the $\inf-\sup$ stability and a point-wise divergence-free discrete velocity field.

Using these ideas, we build a high-performance solver called PetIGA-MF, that is an extension of PetIGA \cite{petiga}, a high-performance isogeometric discretization framework that simplifies modelling and simulation of problems using IGA \cite{Rudraraju,Vignal,Wozniak,Yokota}. PetIGA-MF focuses on multiphysics and multi-field analysis using gradient-conforming spaces as well as curl-, divergence- and integral-conforming discretizations \cite{Cortes2015,nsch15a,nsch15b}.

The paper is organized as follows. In Section 2, we present the strong and weak forms of the generalized Navier-Stokes problem. Section 3 introduces B-spline basis functions, B-spline compatible spaces, and boundary condition imposition. In Section 4, we describe the implementation of our framework. In Section 5, we show the numerical results for all the test cases. We draw conclusions in Section 6.

\section{Generalized Navier-Stokes problem}
We start introducing the generalized Navier-Stokes problem to simplify the description of the incompressible flow problems we address in this paper, which are Darcy, Brinkman, Stokes and Navier-Stokes problems, the difference between them being which physical feature we plan to take into account by the model. The Darcy equation models viscous flows through porous media, whereas Brinkman equation models flow through porous media with an effective viscosity representing high permeability contrasts, for example, when large cavities are present in the medium. The Stokes system model highly viscous flows, while the Navier-Stokes system model flows where the advection is not negligible compared to the diffusivity. These generalizations result in a coupled nonlinear system of partial differential equations for the conservation of linear momentum and mass.

Assuming a steady state system in a bounded open domain $\Omega\in\mathbb{R}^d$ ($d\!=\!2,3$), the problem in its strong form is to find $\BU\!=\!\{\bu,p\}$, with $\bu:\Omega\rightarrow\mathbb{R}^d$, and $p:\Omega\rightarrow\mathbb{R}$ such that:
\begin{align*}
\alpha\nabla\cdot (\bu\otimes \bu) + \beta\bu- \nabla\cdot\boldsymbol{\sigma}(\bu,p)
 &= \mathbf{f}            &\text{in}\     \ &\Omega\       \ \\
\nabla\cdot \bu &= 0      &\text{in}\     \ &\Omega\       \ \\
\bu&=\mathbf{g}	          &\text{on}\  \ &\partial\Omega,
\end{align*}
where $\bu$ is the fluid velocity field, $p$ is the fluid pressure field, $\boldsymbol{\sigma}(\bu,p)\,{=}\,-p\mathbb{I}+2\nu\nabla^{s}\bu$ is the Cauchy stress tensor for an incompressible fluid, with $\mathbb{I}$ being the identity matrix, and $\nabla^{s}\bu$ the symmetric part of the velocity gradient (strain rate), $\nu$ is the kinematic viscosity, $\beta$ is the reaction rate, $\mathbf{f}$ is the body force, and $\mathbf{g}$ is the Dirichlet boundary condition for the velocity. The remaining coefficient, $\alpha$, is used to incorporate or not advective effect on the models, namely, for $\alpha=0$ we have non-advective flows, like Stokes, Darcy, and Brinkman, while $\alpha=1$ incorporates it on the flows, like Navier-Stokes. The different equations models are recovered by varying the coefficients $\alpha, \beta, \nu$. Having $\alpha\!=\!0$ and $\beta\!\gg\!\nu$ represents the Darcy equations, $\alpha\!=\!0$ and $\beta\!\simeq\!\nu$ the Brinkman equations, $\alpha\!=\!0$ and $\beta\!\ll\!\nu$ the Stokes equations, and $\alpha\!=\!1$ and $\beta\!\ll\!\nu$ the Navier-Stokes equations.

Let $(\cdot,\cdot)_\Omega$ denote the $L_2$ inner product in $\Omega$. The trial and weighting spaces for velocity are defined by $\mathcal{V}_{g}$ and $\mathcal{V}_0\!=\!\{\bv\in\mathbf{H}^1(\Omega) :\bv=0\>\text{ on }\>\partial\Omega\}$ respectively, where $\bu \in \mathcal{V}_{g}$ is a lift of a function in $\mathcal{V}$, that is, $\bu = \bv + \mathbf{g}$ for $\bv \in \mathcal{V}_0$. The trial and weighting spaces for pressure is $\mathcal{Q} = L^2(\Omega)$. With these notations the weak form of the problem is to find $\BU\!=\!\{\bu,p\}$, where $\bu \in \mathcal{V}_{g}$ and $p \in \mathcal{Q}$, such that $\forall\BW\!=\!\{\bw,q\}$, where $\bw \in \mathcal{V}_{0}$ and $q \in \mathcal{Q}$:
\begin{equation}\label{eq:wf}
  (\BW,\nl\BU)=B_1(\BW,\BU)+B_2(\BW,\BU,\BU)=L(\BW)\\[0.15cm]
\end{equation}
where
\begin{align*}
B_1(\BW,\BU)=&(\nabla^{s}\bw,2\nu\nabla^{s}\bu)_\Omega+(\bw,\beta\bu)_\Omega-(\nabla\cdot\bw,p)_\Omega+(q,\nabla\cdot\bu)_\Omega\\[0.15cm]
B_2(\BW,\BU,\BU)=&-(\nabla\bw,\alpha(\bu\otimes\bu))_\Omega\\[0.15cm]
L(\BW)=&(\bw,\mathbf{f})_\Omega
\end{align*}
here the bilinear operator $B_1(\cdot,\cdot)$ represents the diffusive and reactive terms of the problem, the trilinear operator $B_2(\cdot,\cdot,\cdot)$ represents the advective term, and the linear operator $L(\cdot)$ represents the forcing term.

\section{Discretization} \label{sec:discretization}
We discretize the weak form of the problem \eqref{eq:wf} using compatible B-spline spaces, namely, using a divergence-conforming space for the velocity and an integral-conforming for the pressure. To simplify the description of such discrete approximation spaces, we give a brief introduction to B-spline functions, and then describe the compatible B-spline spaces as presented in ~\cite{Buffa2011}.

\subsection{B-splines basis functions}
\begin{figure}
\centering
\includegraphics[width=0.7\textwidth]{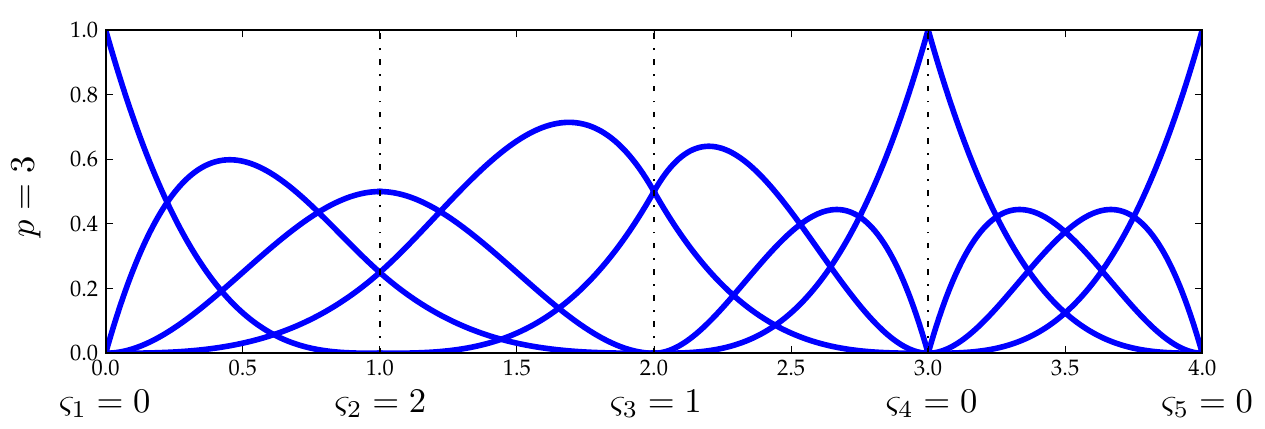}
\captionsetup{belowskip=-0.5cm}
\caption{Example of a cubic ($p\!=\!3$) B-splines basis functions with varying smoothness, quantified by $\varsigma_i$. Dashed lines mark the elements.}
\label{fig:bsplinex}
\end{figure}

B-spline basis functions are piecewise polynomials of degree $p$, defined by specifying the number $n$ of basis functions wanted, the polynomial degree $p$ of the basis, and a knot vector $\Xi=\{ 0=\xi_1,\ldots,\xi_{n+p+1}\!=\!1\}$, which is a finite nondecreasing sequence of real numbers. Additionally, knot multiplicity can be used to control the basis smoothness, see Figure \ref{fig:bsplinex}. The set of B-splines $\left\lbrace B^p_1,\ldots,B^p_n \right\rbrace$ defines a basis with all the properties wanted for analysis purposes \cite{HBCBOOK}. The space spanned by these B-splines is denoted by,
\begin{equation*}
\mathcal{S}^p_{\boldsymbol{\varsigma}} := \text{span} \left\lbrace B^p_i \right\rbrace_{i=1}^n.
\end{equation*}
where $\boldsymbol{\varsigma} := \{\varsigma_1,\ldots,\varsigma_m\}$ is the vector that collects the basis continuity \cite{LesPiegl} at each element boundary.

We describe the trivariate case, that is, when the parametric space is in $\mathbb{R}^3$. The bivariate case follows in a straightforward manner. Given the polynomial orders $p_1, p_2, p_3$, and the numbers of basis $n_1, n_2, n_3$, the trivariate B-spline basis functions are defined by the tensor product of univariate ones as
\begin{equation*}
B^{p_1,p_2,p_3}_{i_1,i_2,i_3} := B^{p_1}_{i_1,1} \otimes B^{p_2}_{i_2,2} \otimes B^{p_3}_{i_3,3}, \quad i_1 = 1,\ldots,n_1;~ i_2 = 1,\ldots,n_2;~ i_3 = 1,\ldots,n_3.
\end{equation*}
Defining the regularity vectors $\boldsymbol{\varsigma}_1, \boldsymbol{\varsigma}_2, \boldsymbol{\varsigma}_3$ in each direction, the trivariate B-spline space is defined by
\begin{equation*}
\mathcal{S}^{p_1,p_2,p_3}_{\boldsymbol{\varsigma}_1,\boldsymbol{\varsigma}_2,\boldsymbol{\varsigma}_3} := \text{span} \left\lbrace B^{p_1,p_2,p_3}_{i_1,i_2,i_3} \right\rbrace_{i_1,i_2, i_3=1}^{n_1,n_2,n_3}.
\end{equation*}
We assume that the regularity vectors $\boldsymbol{\varsigma}_i$ are constant, with components equal to $\varsigma$ (except $\varsigma_1\!=\!\varsigma_m\!=\!0$), unless stated otherwise.

\subsection{Isogeometric (B-spline) differential forms}
The discrete differential forms concept in the context of the finite element method, also known as finite element exterior calculus, is surveyed in~\cite{ExtCalc}. The key aspect of the theory is the use of algebraic topology tools, realized by the existence of de Rham diagrams (exact sequences) relating functional spaces and the image and the kernel of a differential operator between them. These relations are known to hold on the continuous setting, but to inherent such relations on the a discrete setting is a challenging accomplishment since it requires the definition of interpolation and projection operators that renders the commutativity of the de Rham diagrams from the continuous to the discrete setting.

Based on the isogeometric analysis discretization framework Buffa et al. first introduced the isogeometric differential forms in the context of Maxwell equations~\cite{Buffa2010} and Stokes equations~\cite{Buffa2010stokes}, and later developed the general theory in~\cite{Buffa2011}. At the same time, Evans and Hughes~\cite{EvansHughesStokes,EvansHughesNavierStokes,EvansNSUnsteady} applied it to the Generalized Stokes and Navier-Stokes equations. Thus, the isogeometric differential forms, based on B-splines, generate an exact sequence of discrete gradient-, curl-, divergence-, and integral-conforming spaces, that together with the proper interpolation and projection operators, defined in~\cite{Buffa2011}, renders a commutative de Rham diagram. For the construction of de Rham commuting diagram in the context of \textit{hp} finite elements see~\cite{Dem2000}. The novelty of using the isogeometric framework is the possibility of an exact description of the geometry~\cite{HBCBOOK}.

We use divergence- and integral-conforming spaces for the velocity and pressure, respectively, to solve the generalized Navier-Stokes problem. These spaces are defined in the parametric domain as follows:
\begin{center}
\begin{tabular}{lcc}
& Divergence-conforming & Integral-conforming \\[0.20cm]
2D:
	&$\mathcal{S}_{\varsigma_1+1,\varsigma_2}^{p_1+1,p_2}     \times \mathcal{S}_{\varsigma_1,\varsigma_2+1           }^{p_1,p_2+1}$
	&$\mathcal{S}_{\varsigma_1,\varsigma_2}^{p_1,p_2}$\\[0.25cm]
3D:
    &$\mathcal{S}_{\varsigma_1+1,\varsigma_2,\varsigma_3}^{p_1+1,p_2,p_3} \times \mathcal{S}_{\varsigma_1,\varsigma_2+1,\varsigma_3}^{p_1,p_2+1,p_3} \times \mathcal{S}_{\varsigma_1,\varsigma_2,\varsigma_3+1}^{p_1,p_2,p_3+1}$
    &$\mathcal{S}_{\varsigma_1,\varsigma_2,\varsigma_3}^{p_1,p_2,p_3}$\\[0.25cm]
\end{tabular}
\end{center}

\begin{figure}
\centering
\includegraphics[width=0.7\textwidth]{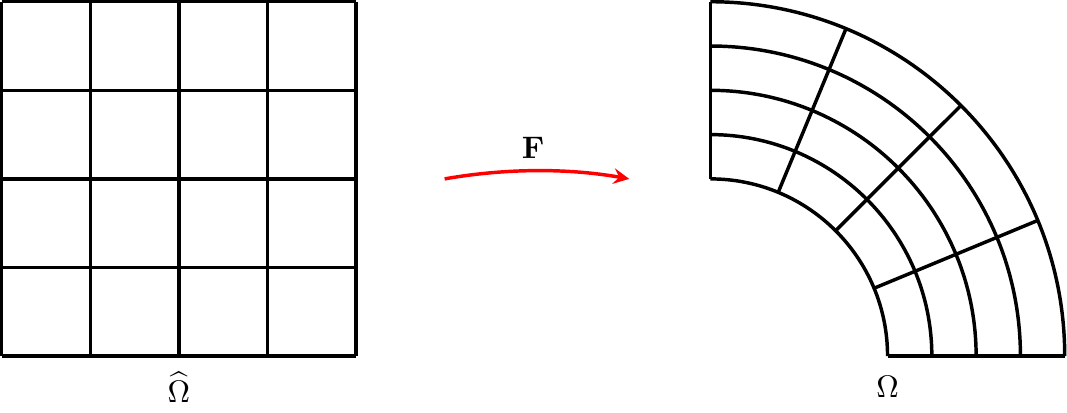}
\caption{Geometric mapping $\mathbf{F}$.}
\label{fig:geomap}
\end{figure}

The main consequence of those definitions is that one can prove that the divergence operator is surjective for each pair of spaces above. In order to move the definitions from the parametric to the physical domain, a preserving push-forward mapping is used for every space in the sequence, guaranteeing that the de Rham diagram for the spaces defined on the physical domain (see~\cite{Buffa2011}) also commutes in the discrete setting, that is, the spaces mapped to the physical domain also define an exact sequence. In our case, the relevant preserving mappings for our spaces of interest are the pullbacks: 
\begin{align*}
\iota_\bu(\bv) &= \mathrm{det}\left(D\mathbf{F}\right)\left(D\mathbf{F}\right)^{-1} \left(\bv\circ\mathbf{F}\right) & \bv \in \mathbf{H}(\mathrm{div};\Omega), \\
\iota_p(q)  &=  \mathrm{det}\left(D\mathbf{F}\right)\left(q\circ\mathbf{F}\right) & q \in \mathrm{L^2}(\Omega),
\end{align*}
where $\mathbf{F}$ is the geometric mapping from the parametric domain $\widehat{\Omega}$ onto the physical domain $\Omega$ (Figure \ref{fig:geomap}), and $D\mathbf{F}$ is the gradient of the geometric mapping. The divergence-preserving map, $\iota_\bu$, is the Piola transformation \cite{Gonz}, and $\iota_p$ is the integral preserving transformation. As one can infer the polynomial order and continuity of our basis functions in the physical domain, depend not only on the basis functions adopted in the parametric space but also on the geometric mapping used. 

Additionally, the commutativity of the de Rham diagram, with respect to the projections, for the case of the discrete velocity and pressure pair above guarantees the stability of the scheme, that is, the discrete $\inf-\sup$ condition of the pair, and that the satisfaction of the weak incompressibility condition implies it holds strongly, that is, $\nabla \cdot \bu = 0$ pointwise (for a proof see~\cite{EvansHughesStokes}).

\subsection{Boundary Condition Imposition}
We impose the normal boundary conditions on the velocity strongly but doing the same for the tangential boundary conditions on the velocity with divergence-conforming basis functions may lead to unstable discretizations in domains with corners. Thus, we use Nitsche's method for weak boundary imposition to avoid this problem, alleviating the necessity for highly refined meshes to reproduce the layer effect on the no-slip boundary conditions \cite{WeakImpos1,WeakImpos2}. The weak imposition of tangential boundary conditions modifies the operators $B_1(\BW,\BU)$ and $L(\BW)$, introducing the adjoint consistency and the penalization terms, we then get
\begin{align*}
\widehat{B}_1(\BW,\BU)& =B_1(\BW,\BU)          \\
            & -\left(\bw,2\nu\nabla^{s} \bu \cdot \mathbf{n} \right)_\Gamma   &&\text{Consistency}\\
            & -\textcolor{blue}{\left(\bu,2\nu\nabla^{s} \bw \cdot \mathbf{n} \right)_\Gamma}   &&\text{Adjoint consistency}\\
            & +\textcolor{red}{\left(\bw,2\nu\nabla^{s} \bu \alpha_p  \right)_\Gamma}          &&\text{Penalization,}\\[0.15cm]
\end{align*}
\vspace*{-1cm}
\begin{align*}
\widehat{L}(\BW)& =L(\BW)                                            \\
            & -\textcolor{blue}{\left(\mathbf{g},2\nu\nabla^{s} \bw \cdot \mathbf{n} \right)_\Gamma}   &&\text{Adjoint consistency}\\
            & +\textcolor{red}{\left(\bw,\nu\alpha_p \mathbf{g}   \right)_\Gamma}                     &&\text{Penalization.}
\end{align*}
where $\alpha_p\,{=}\,C_{pen}/h_f$, here $C_{pen}\!=\!5(p+1)$ is the penalty term parameter depending on the polynomial order $p$ of the discretization, and $h_f$ is the wall normal mesh size \cite{WeakImpos1}.

\section{Implementation}
In this section, we describe an extension of PetIGA, which adds a flexible and scalable parallel implementation of multi-field isogeometric discretizations, where the discrete fields can belong to any conforming space of the B-spline de Rham sequence, that is, gradient-, curl-, divergence- and integral-conforming spaces. We first introduce the basic structures of PetIGA, and then present PetIGA-MF and the new structures that implement multi-fields discretizations in a user-friendly manner.

\subsection{PetIGA}
PetIGA \cite{petiga} is a framework based on PETSc \cite{petsc}, which uses its parallel tools to solve a discrete variational formulation (Galerkin or collocation method) of partial differential equations. The discretization is built using B-spline functions and a patch-wise isoparametric mapping. Different structures built in PetIGA contain all the information the user needs to code the discrete variational formulation at a quadrature/collocation point. Regarding the data-structure for a structured mesh and its partitioning, PetIGA implements its data-structure, similar to a DM in the PETSc jargon. In this case, it is tailored to the specifics of isogeometric analysis, particularly the use of high continuous basis functions, with possibly arbitrary continuity orders across elements boundaries, and the respective connectivity array to promote the assembly of the global matrices from their local contributions. With respect to the synergy between geometry description and finite element analysis, the data-structure called IGA provides the abstraction of a spline patch, together with the elemental and quadrature information needed to integrate a variational form when we use Galerkin's method or collocation schemes~\cite{Collocation}.

\begin{figure}[H]
\captionsetup{aboveskip=0.1cm}
\captionsetup{belowskip=-0.3cm}
\centering
\includegraphics[width=0.70\textwidth]{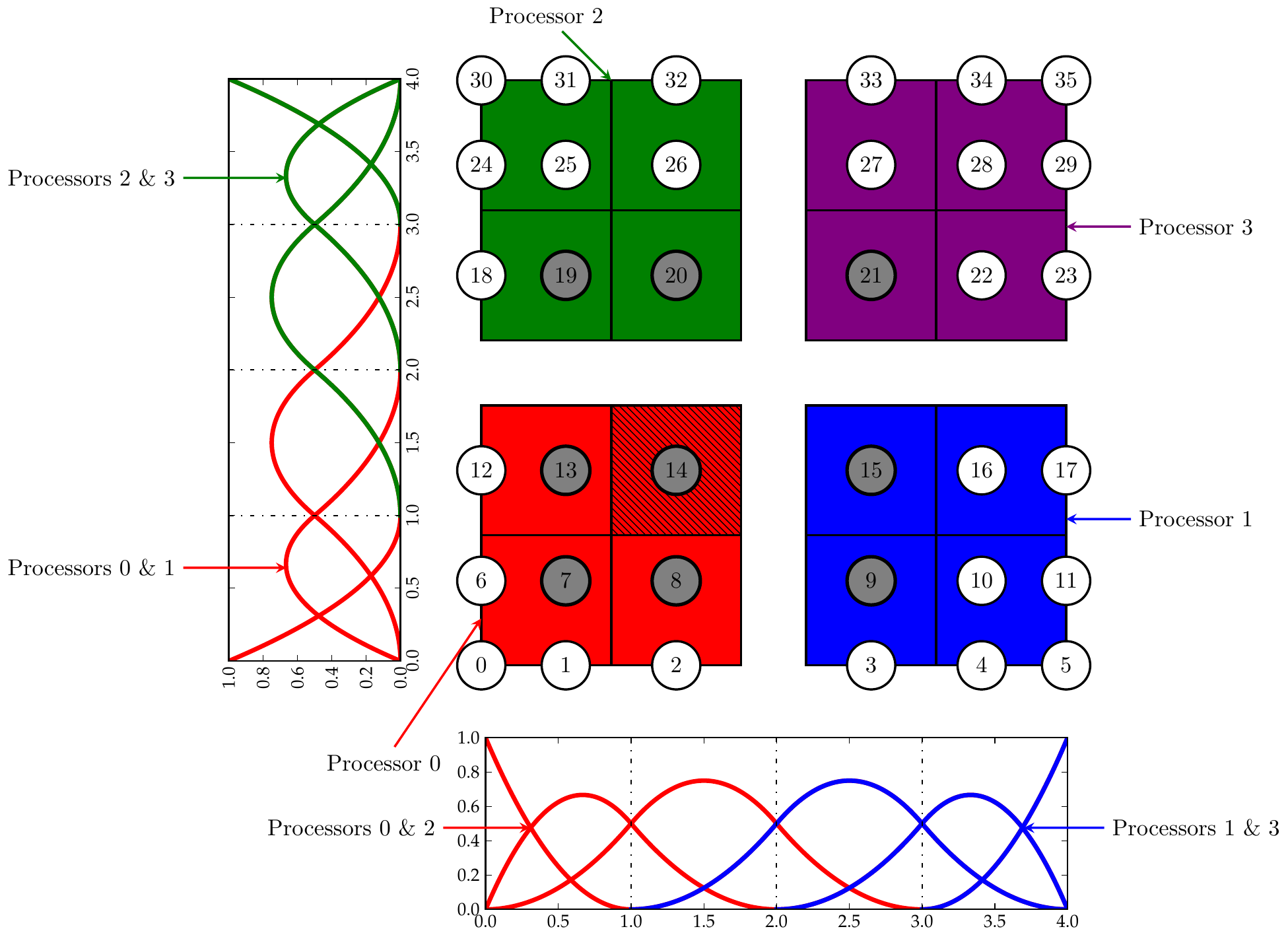}
\caption{Distribution of elements for a mesh of $4\!\times\!4$ elements on a grid of $2\!\times\!2$ processors, for a spline space with $p_1\!=p_2\!=\!2$ and $\varsigma_1\!=\varsigma_2\!=\!1$ regularity. Grey-filled nodes represent the basis functions with support on the dashed element.}
\label{fig:meshpart}
\end{figure}

The mesh is split, to balance the workload between the processors, according to a calculation of a box stencil, distributing the elements through the grid of processors, and then assigning the degrees of freedom that lie on the interfaces to one of the neighboring processors. When having an uneven distribution of elements, PetIGA is programmed to assign the higher workload to the next processor in the grid, to the left or the bottom, depending on the interface. Figure \ref{fig:meshpart} shows an example of a spline space $\mathcal{S}_{\varsigma_1,\varsigma_2}^{p_1,p_2}$ defined over a mesh of $4\!\times\!4$ elements $p_1\!=p_2\!=\!2$ and $\varsigma_1\!=\varsigma_2\!=\!1$ basis, and its splitting through a grid of $2\!\times\!2$ processors.

\definecolor{dgreen}{HTML}{008000}
\definecolor{dmagenta}{HTML}{800080}
\definecolor{lblue}{HTML}{2A98FF}
\definecolor{lgreen}{HTML}{00FF08}
\definecolor{lmagenta}{HTML}{FF0464}

\begin{figure}[H]
\captionsetup{aboveskip=0.1cm}
\captionsetup{belowskip=-0.3cm}
\centering
\subfloat[Natural numbering.]{\includegraphics[width=0.225\textwidth]{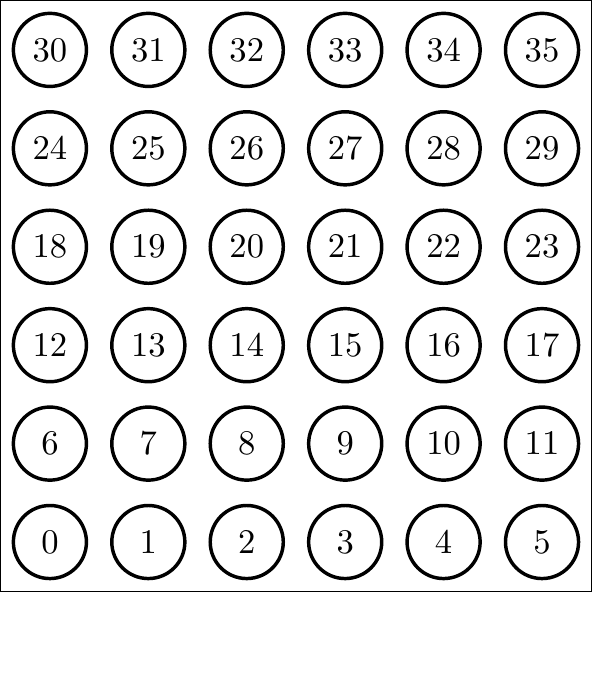} \label{fig:natural-numbering}}
\quad
\subfloat[Global numbering.]{\includegraphics[width=0.225\textwidth]{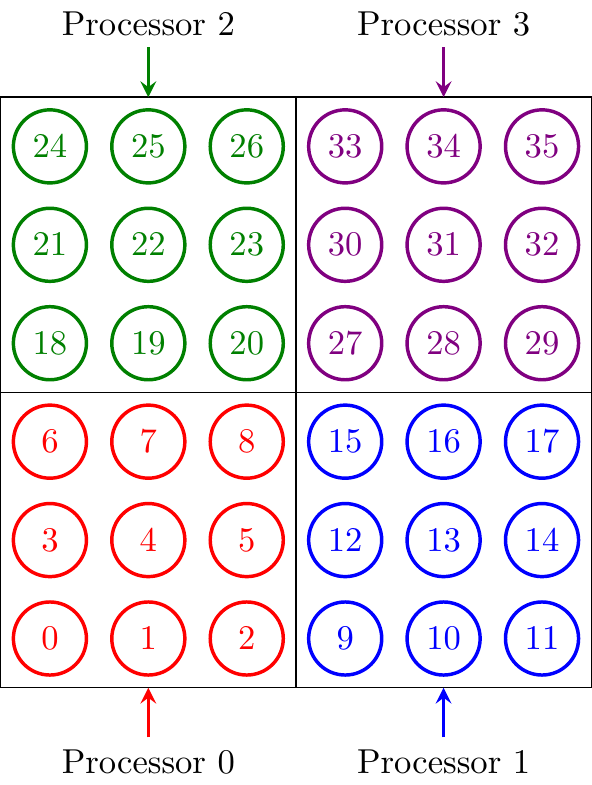} \label{fig:global-numbering}}
\quad
\subfloat[Local numbering.]{\includegraphics[width=0.226\textwidth]{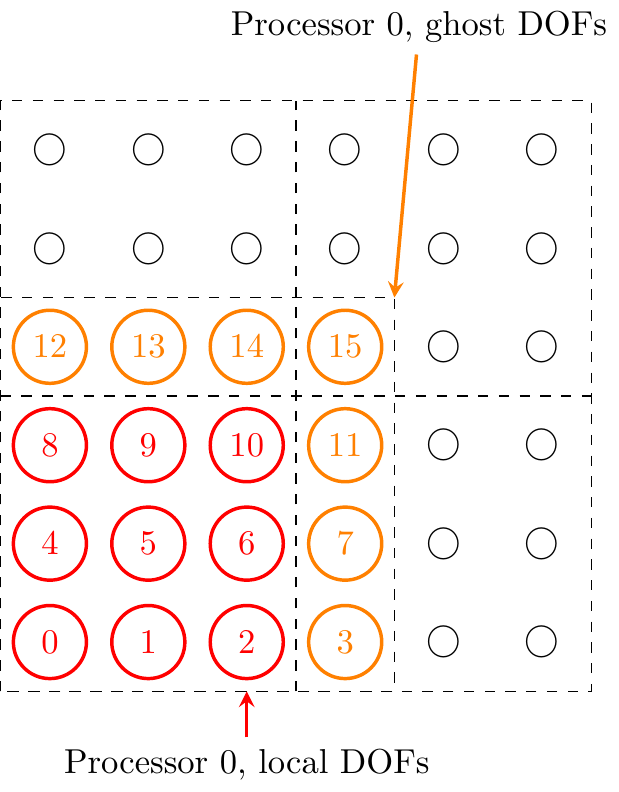} \label{fig:local-numbering}}
\caption{Natural, global and local numbering.}
\label{fig:numbering}
\end{figure}

\begin{figure}
\captionsetup{aboveskip=0.1cm}
\captionsetup{belowskip=-0.3cm}
\centering
\subfloat[Global vector.]{\includegraphics[height=0.45\textwidth]{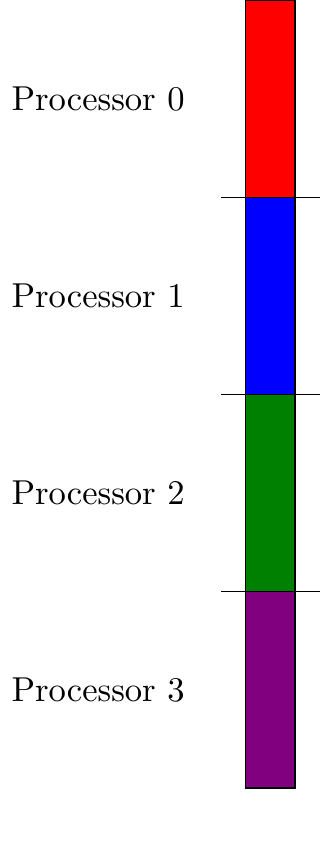} \label{fig:global-vecs}}
\hspace{2cm}
\subfloat[Local vectors.]{\includegraphics[height=0.45\textwidth]{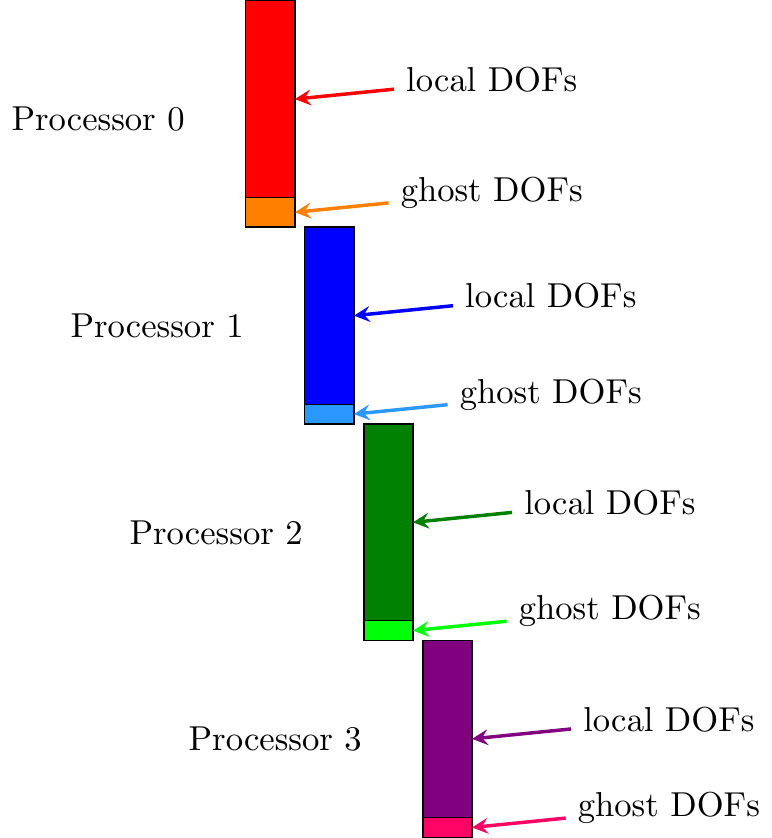} \label{fig:local-vecs}}
\caption{Global and local vectors.}
\label{fig:vecs}
\end{figure}

As shown in Figure \ref{fig:meshpart}, and also reproduced in Figure \ref{fig:natural-numbering}, the basis functions are naturally ordered in a lexicographic way, called natural numbering in PETSc jargon. Once in parallel such numbering is not convenient anymore, and the mesh splitting among processors induces a new numbering where the degrees of freedom that belong to the same processor are numbered first (see Figure \ref{fig:global-numbering}). It is referred to as global numbering. Global vectors (see Figure \ref{fig:global-vecs}) are associated with this numbering. For processors to be able to solve in parallel, the information of the ``ghost degrees of freedom'' must communicate from neighboring processors. For such task a local numbering is more convenient as Figure \ref{fig:local-numbering} shows. Local vectors are associated with this numbering as shown in Figure \ref{fig:local-vecs} where lighter colors represent the ghost degrees of freedom. The amount of communication between processors depends on the continuity of the basis. All processor communications are hidden from the user and managed internally by the PetIGA data-structures.


\subsection{PetIGA-MF}
PetIGA-MF is a multi-field extension of PetIGA, where different discretization spaces can be used for each field, making it suitable for solving multi-physics problems. All scalar and vector structure-preserving B-spline discrete spaces mentioned in section \ref{sec:discretization} are available. To simplify the access to the information of the different fields, we create new structures on top of the ones already existent in PetIGA, combining the single field data-structures to work in a multi-field framework.

\begin{figure}
\captionsetup{aboveskip=0.1cm}
\captionsetup{belowskip=-0.3cm}
\centering
\includegraphics[width=0.7\textwidth]{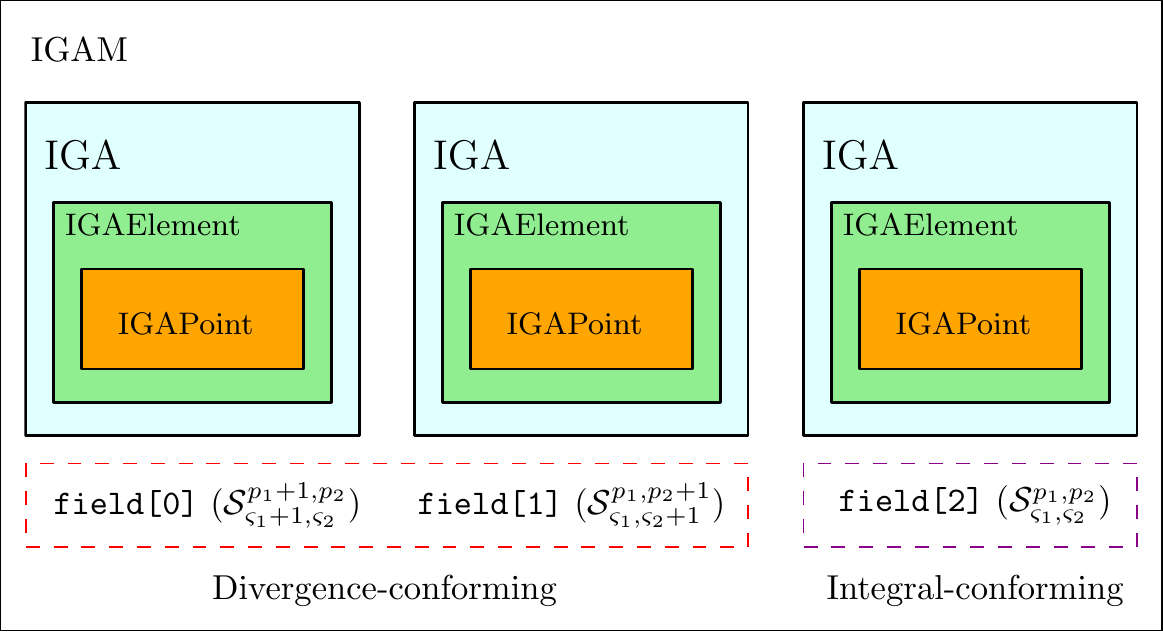}
\caption{Discrete velocity and pressure spaces abstraction used in PetIGA-MF.}
\label{fig:petiga2}
\end{figure}

\begin{figure}
\captionsetup{aboveskip=0.1cm}
\captionsetup{belowskip=-0.3cm}
\centering
\subfloat[$4\!\times\!4$ mesh used to define all the fields]{\includegraphics[width=0.27\textwidth]{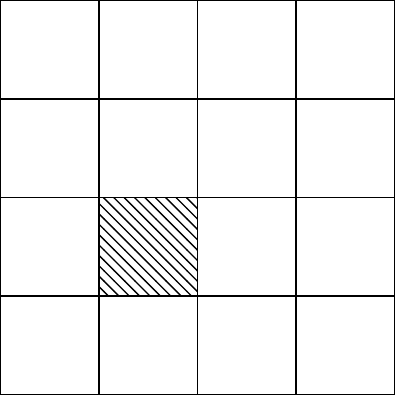} \label{fig:mesh-all-fields}}
\\
\subfloat[Natural numbering for the space $\mathcal{S}_{1,0}^{2,1}$ (\texttt{field[0]}).]{\includegraphics[width=0.27\textwidth]{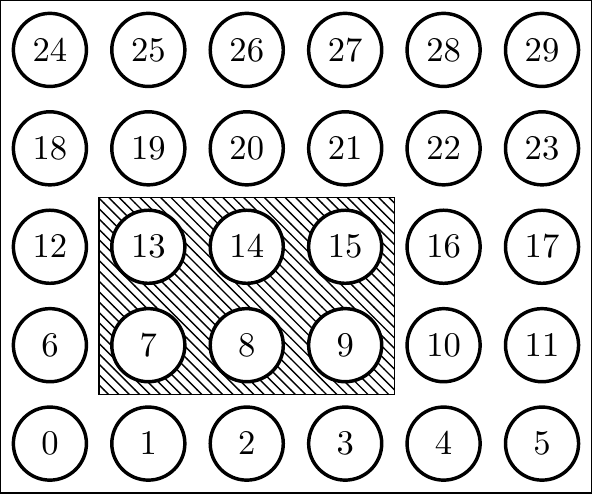} \label{fig:dof-field-0}}
\quad
\subfloat[Natural numbering for the space $\mathcal{S}_{0,1}^{1,2}$ (\texttt{field[1]}).]{\includegraphics[width=0.27\textwidth]{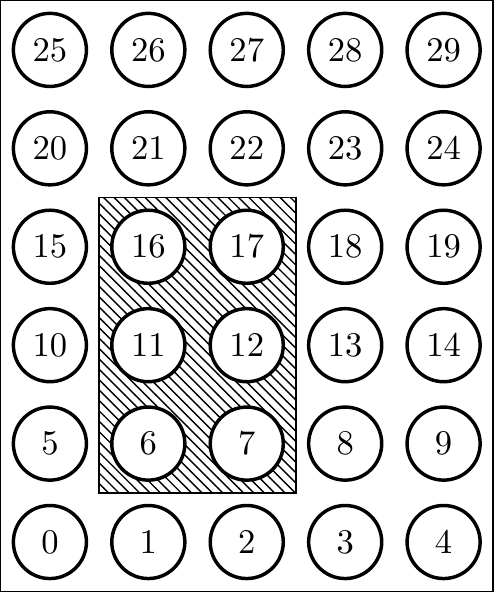} \label{fig:dof-field-1}}
\quad
\subfloat[Natural numbering for the space $\mathcal{S}_{0,0}^{1,1}$ (\texttt{field[2]}).]{\includegraphics[width=0.27\textwidth]{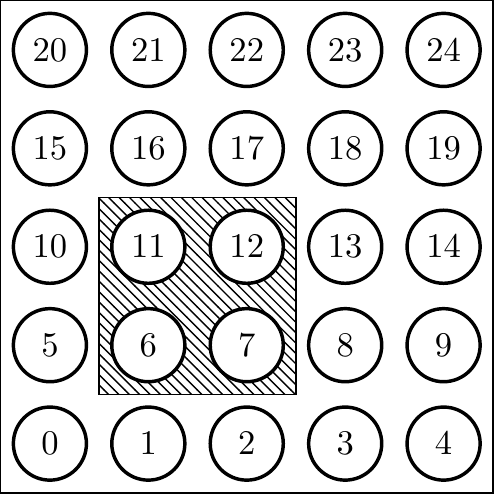} \label{fig:dof-field-2}}
\caption{Natural numbering for the degrees of freedom of all fields and basis functions numbers with support on the dashed element.}
\label{fig:multidof}
\end{figure}

PETSc provides a data management subclass, called DMComposite, which allows one to pack several fields in a monolithic blocked way for multi-field and multi-physics discretizations. A new IGAM class packs the IGAs for each field together with an instance of a DMComposite. Once we create a IGAM object, the discrete spaces are set by assigning a type of structure-preserving space (gradient-conforming is the default type), and the corresponding fields of it. Figure \ref{fig:petiga2} gives a schematic representation for the case of a two-dimensional divergence- and integral-conforming velocity-pressure pair of B-spline spaces. Two constraints built into PetIGA-MF are that all the fields, that are B-spline spaces, need to be defined on a mesh with the same number of elements and to use the same number of quadrature points per element. Figure \ref{fig:multidof} illustrates the natural numbering of the three fields for the divergence- and integral-conforming pair of spaces shown in Figure \ref{fig:petiga2}, defined on a mesh of $4\!\times\!4$ elements (see Figure \ref{fig:mesh-all-fields}) with $p_1\!=p_2\!=\!1$ and $\varsigma_1\!=\varsigma_2\!=\!0$. Figures \ref{fig:dof-field-0}, \ref{fig:dof-field-1} and \ref{fig:dof-field-2} emphasize the basis functions with support on the dashed element in Figure \ref{fig:mesh-all-fields} for every field, that have a direct impact on the parallel partitioning of the degrees of freedom of each field.

\begin{figure}
\captionsetup{aboveskip=0.1cm}
\captionsetup{belowskip=-0.3cm}
\centering
\subfloat[Global vector.]{\includegraphics[height=0.45\textwidth]{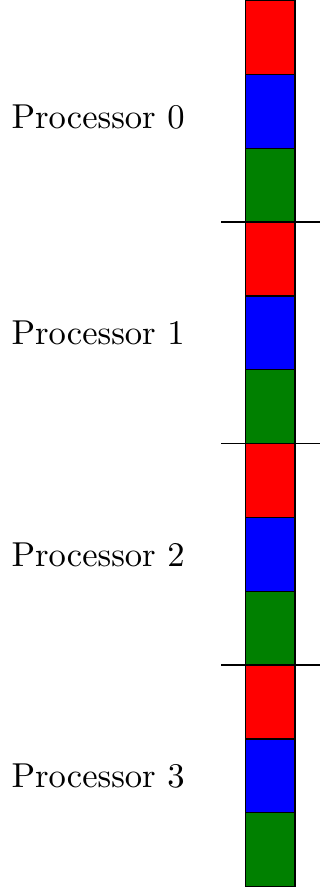} \label{fig:multivec-global}}
\hspace{1.0cm}
\subfloat[Split global vectors.]{\includegraphics[height=0.4\textwidth]{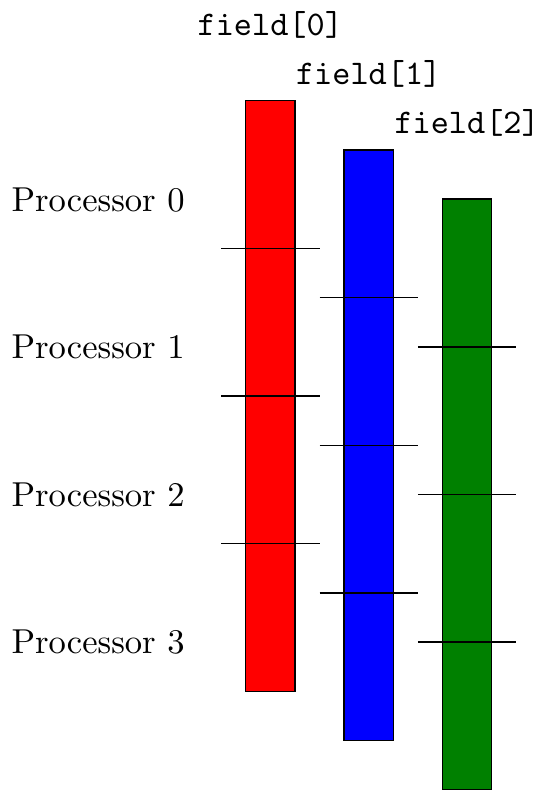} \label{fig:multivec-split-global}}
\hspace{1.0cm}
\subfloat[Split local vectors.]{\includegraphics[height=0.5\textwidth]{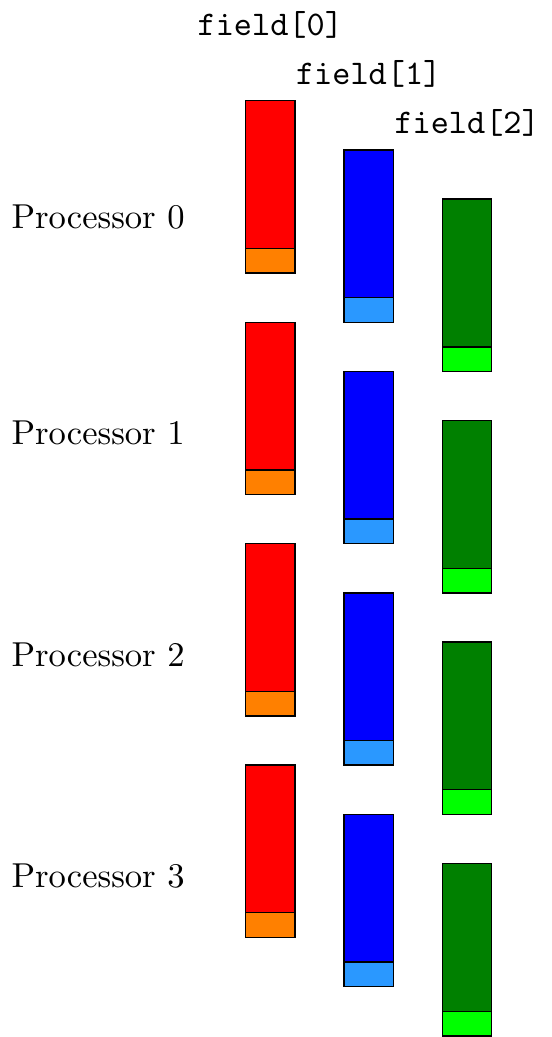} \label{fig:multivec-split-local}}
\caption{Global, split global and split local vectors for three fields (now represented by the different colors) and four processors.}
\label{fig:multivec}
\end{figure}

In PetIGA-MF every processor owns the part of every field that corresponds to its part of the mesh, and with respect to the global numbering the vector for a multi-field problem is schematically represented as in Figure \ref{fig:multivec-global}. To solve the fields in parallel, a first step is to create an independent vector for each field. We create these vectors by splitting the global vector into fields, obtaining the split global vectors as Figure \ref{fig:multivec-split-global} shows. The second step is to follow the same procedure as in PetIGA, for every split global vector, we obtain the split local vectors, that incorporates the \textit{ghost} degrees of freedom as shown in Figure \ref{fig:multivec-split-local}.

\subsubsection{Mapped basis functions.}
We stuck with PetIGA's philosophy, namely, that the framework delivers to the user the basis functions and their derivatives already mapped to the physical space, called in this case shape functions. In this way, the user can directly code the variational formulation. Since in the multi-field setting we can have a mixing of scalar and vector discrete spaces, we create a three indexed array of pointers, $\texttt{*shape[d][i][j]}$, to store the shape functions and their derivatives evaluated at the quadrature points of an element. The index $\texttt{d} = 0,1,2$ selects: the shape function, $\texttt{d} = 0$, its first derivative, $\texttt{d} = 1$, and its second derivative, $\texttt{d} = 2$. The indices \texttt{i} and \texttt{j} stand for the field components. Such an indexing is needed because of the use of mapped vector basis functions, for example, the divergence-conforming space on the physical domain.

Given the nature of the mappings used for the discrete vector spaces, for example the Piola transformation in the case of the divergence-conforming spaces, a component of a vector basis function in the parametric space, is mapped into a linear combination of all the parametric components, coupling them all in the physical space. Indeed, consider the example of the divergence-conforming space depicted in Figure \ref{fig:petiga2}, and let $\{ \widehat{N}_{1}^{u}, \widehat{N}_{2}^{u}, \ldots, \widehat{N}_{n_u}^{u} \}$ and $\{ \widehat{N}_{1}^{v}, \widehat{N}_{2}^{v}, \ldots, \widehat{N}_{n_v}^{v} \}$ represent the basis functions with support on an element for the spaces $\mathcal{S}_{\varsigma_1+1,\varsigma_2}^{p_1+1,p_2}$ (\texttt{field[0]}) and $\mathcal{S}_{\varsigma_1,\varsigma_2+1}^{p_1,p_2+1}$ (\texttt{field[1]}) respectively. The vector basis functions of the parametric space $\mathcal{S}_{\varsigma_1+1,\varsigma_2}^{p_1+1,p_2} \times \mathcal{S}_{\varsigma_1,\varsigma_2+1  }^{p_1,p_2+1}$ with support on the same element will be
\begin{align}
\left\lbrace 
\begin{pmatrix}
\widehat{N}_{1}^{u} \\ 0
\end{pmatrix},
\begin{pmatrix}
\widehat{N}_{2}^{u} \\ 0
\end{pmatrix},
\ldots,
\begin{pmatrix}
\widehat{N}_{n_u}^{u} \\ 0
\end{pmatrix},
\begin{pmatrix}
0 \\ \widehat{N}_{1}^{v}
\end{pmatrix},
\begin{pmatrix}
0 \\ \widehat{N}_{2}^{v}
\end{pmatrix},
\ldots,
\begin{pmatrix}
0 \\ \widehat{N}_{n_v}^{v}
\end{pmatrix}
\right\rbrace.
\label{eq:parametric-basis-func}
\end{align}

Applying the push-forward transformation $\iota_\mathbf{u}^{-1}(\bu) = \mathrm{det}\left(D\mathbf{F}\right)^{-1} \left(D\mathbf{F}\right) (\bu)$ to the set of parametric basis functions \ref{eq:parametric-basis-func} we obtain the mapped basis function
\begin{align}
&\left\lbrace
\iota_\mathbf{u}^{-1} 
\begin{pmatrix}
\widehat{N}_{1}^{u} \\ 0
\end{pmatrix},
\iota_\mathbf{u}^{-1}
\begin{pmatrix}
\widehat{N}_{2}^{u} \\ 0
\end{pmatrix},
\ldots,
\iota_\mathbf{u}^{-1}
\begin{pmatrix}
\widehat{N}_{n_u}^{u} \\ 0
\end{pmatrix},
\iota_\mathbf{u}^{-1}
\begin{pmatrix}
0 \\ \widehat{N}_{1}^{v}
\end{pmatrix},
\iota_\mathbf{u}^{-1}
\begin{pmatrix}
0 \\ \widehat{N}_{2}^{v}
\end{pmatrix},
\ldots,
\iota_\mathbf{u}^{-1}
\begin{pmatrix}
0 \\ \widehat{N}_{n_v}^{v}
\end{pmatrix}
\right\rbrace = \\
& = \left\lbrace 
\begin{pmatrix}
\mathbf{N}_{1}^{u} \\ \mathbf{N}_{1}^{v}
\end{pmatrix},
\begin{pmatrix}
\mathbf{N}_{2}^{u} \\ \mathbf{N}_{2}^{v}
\end{pmatrix},
\ldots,
\begin{pmatrix}
\mathbf{N}_{n_u}^{u} \\ \mathbf{N}_{n_u}^{v}
\end{pmatrix},
\begin{pmatrix}
\mathbf{N}_{n_u + 1}^{u} \\ \mathbf{N}_{n_u + 1}^{v}
\end{pmatrix},
\begin{pmatrix}
\mathbf{N}_{n_u + 2}^{u} \\ \mathbf{N}_{n_u + 2}^{v}
\end{pmatrix},
\ldots,
\begin{pmatrix}
\mathbf{N}_{n_u + n_v}^{u} \\ \mathbf{N}_{n_u + n_v}^{v}
\end{pmatrix}
\right\rbrace.
\label{eq:physical-basis-func}
\end{align}

In terms of implementation the pointers to the shape functions can be schematically represented like in Table \ref{tab:pointers-basis}. We emphasize that for the gradient-conforming vector basis functions, that is the standard $\mathbf{H}^1$ basis used in a stabilized formulation, such coupling of the components does not occur and the pointers to the shape functions are also diagonal on the physical space.   

\setlength\arraycolsep{10pt}
\begin{table}[h]
\centering
\begin{tabular}{@{}cc@{}}
\toprule
\multicolumn{1}{c}{Parametric space} & \multicolumn{1}{c}{Physical space} \\ \midrule \\
$\texttt{*shape[0][i][j]} =$ & $\texttt{*shape[0][i][j]} =$ \\
 & \\
$\begin{bmatrix}
\mathcal{S}_{\varsigma_1+1,\varsigma_2}^{p_1+1,p_2} & \times & \times \\[0.25cm]
\times & \mathcal{S}_{\varsigma_1,\varsigma_2+1}^{p_1,p_2+1} & \times\\[0.25cm]
\times & \times & \mathcal{S}_{\varsigma_1,\varsigma_2}^{p_1,p_2}
\end{bmatrix}$ 
&
$
\begin{bmatrix}
	J^{-1} \left(D\mathbf{F}\right)_{x,X}\mathcal{S}_{\varsigma_1+1,\varsigma_2}^{p_1+1,p_2}&
	J^{-1} \left(D\mathbf{F}\right)_{x,Y}\mathcal{S}_{\varsigma_1,\varsigma_2+1}^{p_1,p_2+1}&
	\times \\[0.25cm]
	J^{-1} \left(D\mathbf{F}\right)_{y,X}\mathcal{S}_{\varsigma_1+1,\varsigma_2}^{p_1+1,p_2}&
	J^{-1} \left(D\mathbf{F}\right)_{y,Y}\mathcal{S}_{\varsigma_1,\varsigma_2+1}^{p_1,p_2+1}&
	\times \\[0.25cm]
	\times & \times & J^{-1} \mathcal{S}_{\varsigma_1,\varsigma_2}^{p_1,p_2}
\end{bmatrix}
$                   
\end{tabular}
\caption{Schematic representation of the pointers to the basis functions before and after being mapped to the physical domain for velocity-pressure pair depicted in Figure \ref{fig:petiga2}. Here $J = \mathrm{det}\left(D\mathbf{F}\right)$ means the Jacobian, and $\times$ the \texttt{NULL} pointer.}
\label{tab:pointers-basis}
\end{table}

\section{Numerical Results}

In this section, we present convergence tests using manufactured solutions in 2D and 3D for our implementation, showing optimal convergence rates for both parametric and physical domains. Results for different Reynolds numbers and Damk\"{o}hler numbers were obtained, where the Reynolds and Damk\"{o}hler numbers are defined as
\begin{align*}
Re=\frac{UL}{\nu},\quad
Da=\frac{\beta L^2}{\nu}
\end{align*}
these variations cover all the different problems proposed, where having $Da\!=\!1000$ and $\alpha\!=\!0$ represents the Darcy flow model, $Da\!=\!1$ and $\alpha\!=\!0$ the Brinkman flow model, $Da\!=\!0$ and $\alpha\!=\!0$ the Stokes flow model, and $Da\!=\!0$ and $\alpha\!=\!1$ represents the Navier-Stokes flow model. All test cases consider equal polynomial order $p$ in every direction, and maximum continuity ($\varsigma\!=\!p\!-\!1$) for the pressure space. We ran all the test cases on a workstation (2 Hex-core Xeon X5650, 48 Gb RAM).

\subsection{Solution in a unitary square}

\begin{figure}
\centering
\subfloat[Velocity magnitude.]{\includegraphics[width=0.45\textwidth]{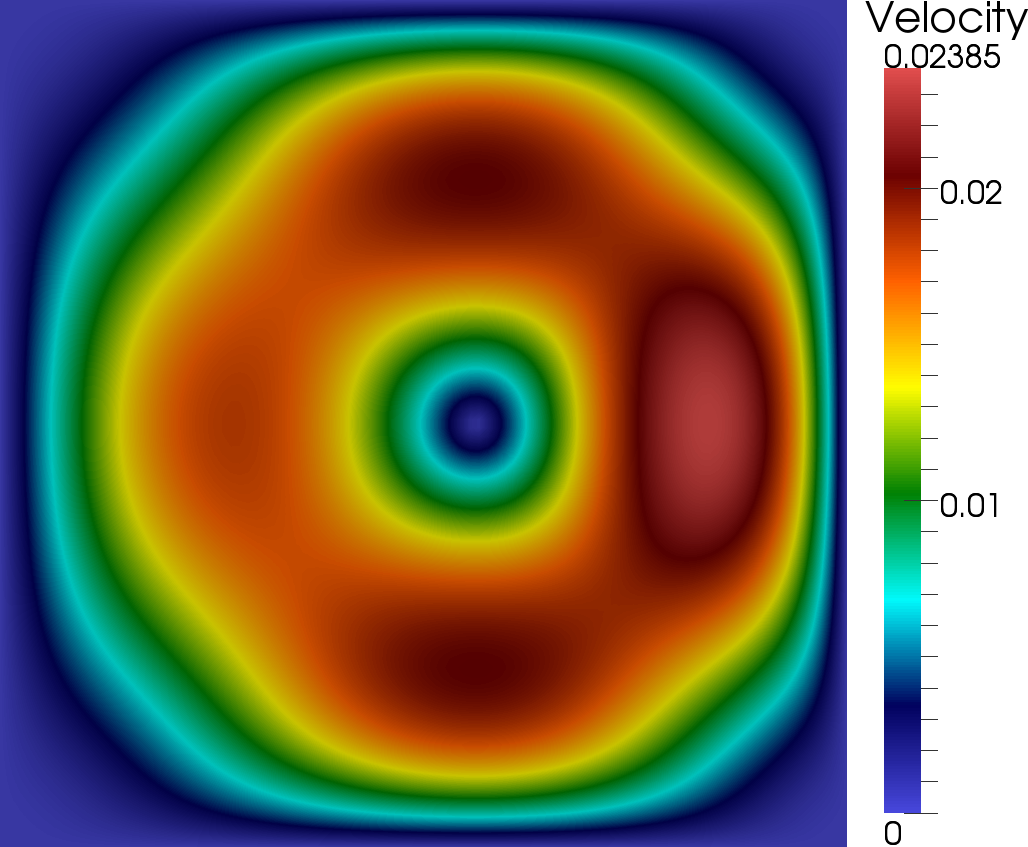}}
\hspace{0.2cm}
\subfloat[Pressure.]{\includegraphics[width=0.45\textwidth]{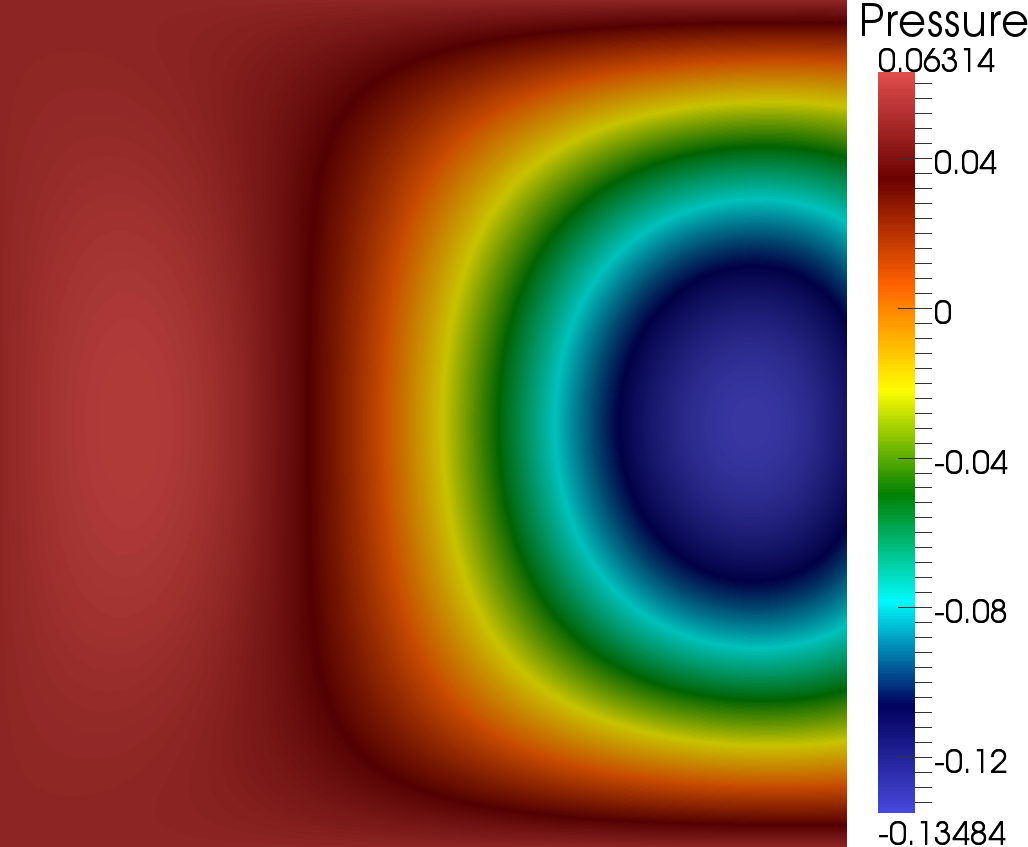}}
\caption{Analytical solution for the square problem.}
\label{fig:squareresults}
\end{figure}

Here we present the solution of the two-dimensional flow in a unitary square shown in \cite{Buffa2010stokes}. We compute $\bu$ and $p$, when a force $\mathbf{f}$ is imposed, and compare the numerical solution with the analytical solution $\overline{\bu}$ and  $\overline{p}$ (Figure \ref{fig:squareresults}).

\begin{align*}
\overline{\bu}=&\begin{bmatrix}
 2e^x(-1 + x)^2 x^2 (y^2 - y)(-1 + 2y)\\
 (-e^x (-1 + x)x(-2 + x(3 + x))(-1 + y)^2 y^2 )
\end{bmatrix}\\[0.25cm]
\overline{p}=&(-424 + 156e + (y^2 - y)(-456 + e^x (456 + x^2 (228 - 5(y^2 - y))+\\
&2x(-228 + (y^2 - y)) + 2x^3 (-36 + (y^2 - y)) + x^4 (12 + (y^2 - y)))))
\end{align*}

\begin{figure}[H]
\centering
\subfloat{\includegraphics[width=0.3\textwidth]{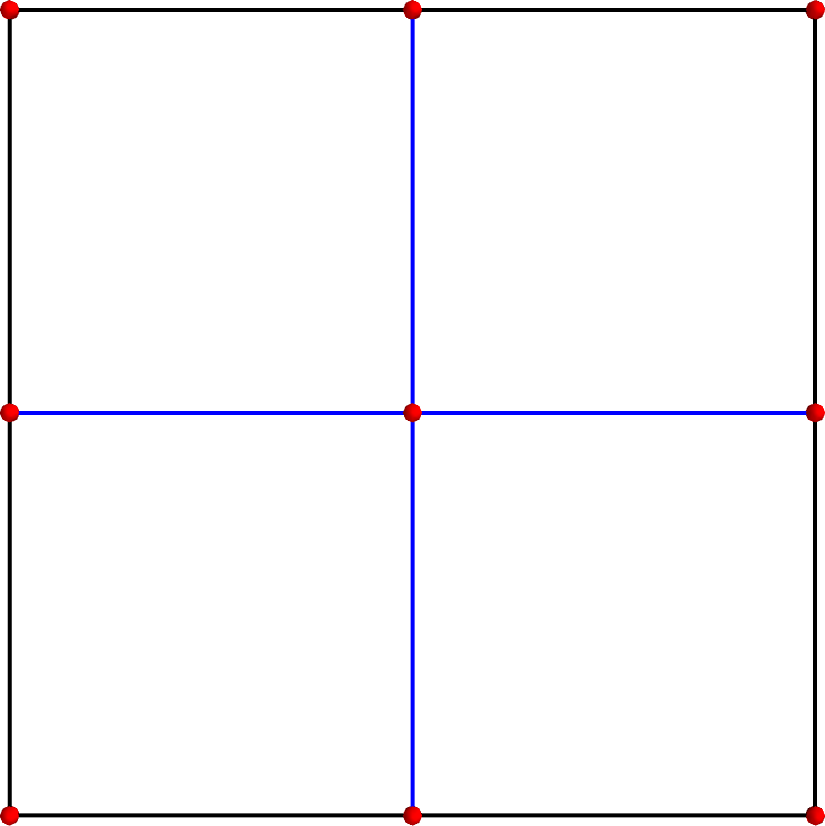}}
\hspace{1cm}
\subfloat{
\begin{tikzpicture}[scale=0.82]
\node[inner sep=0pt] (russell) at (0,0)
    {\includegraphics[width=0.3\textwidth]{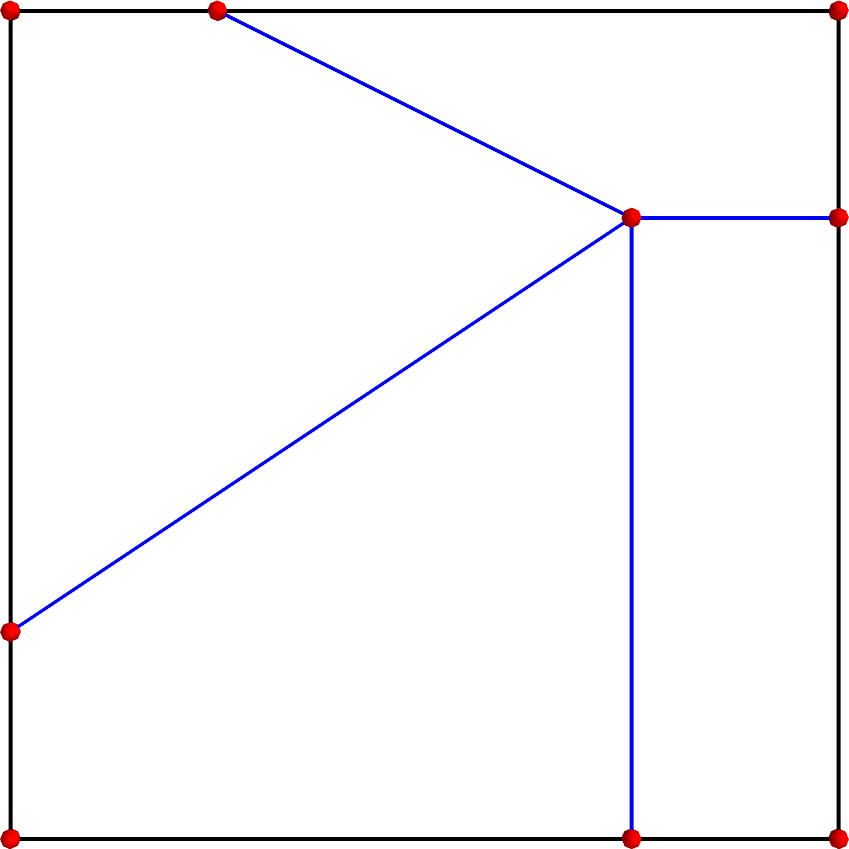}};
    
\draw[dashed,blue,line width=0.4pt](   0cm,-2.7cm)--(   0cm, 2.7cm);
\draw[dashed,blue,line width=0.4pt](-2.7cm,   0cm)--( 2.7cm,   0cm);

\draw[>=stealth,<->,red,line width=1.5pt](-2.68cm,    0cm)--(-2.68cm, -1.3cm);
\draw[>=stealth,<->,red,line width=1.5pt](  2.7cm,    0cm)--(  2.7cm, 1.28cm);
\draw[>=stealth,<->,red,line width=1.5pt](    0cm,    0cm)--(  1.3cm,  1.3cm);
\draw[>=stealth,<->,red,line width=1.5pt](    0cm, 2.69cm)--(-1.28cm, 2.69cm);
\draw[>=stealth,<->,red,line width=1.5pt](    0cm, -2.7cm)--(  1.3cm, -2.7cm);

\draw (-2.5,-0.7) node {d};
\draw ( 2.5, 0.6) node {d};
\draw (-0.6, 2.5) node {d};
\draw ( 0.7,-2.5) node {d};

\end{tikzpicture}
}
\\ \setcounter{subfigure}{0}
\subfloat[Undistorted mesh.]{\includegraphics[width=0.3\textwidth]{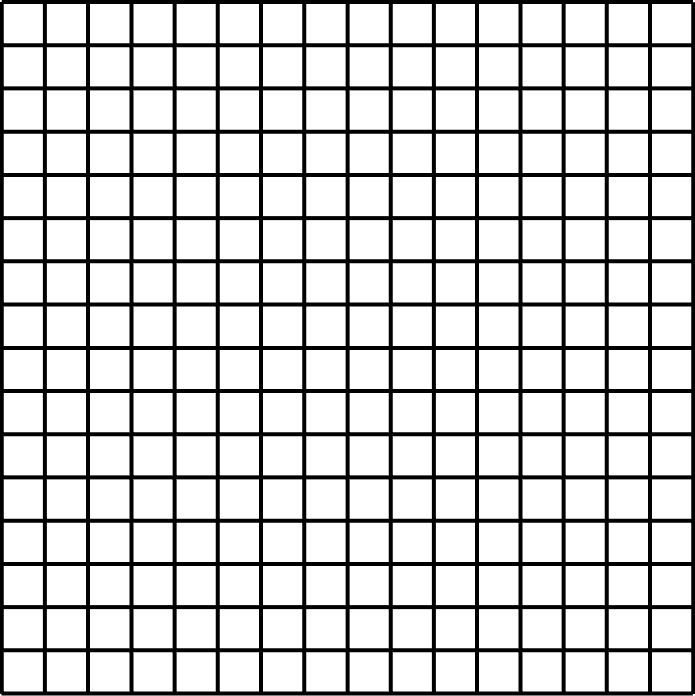}}
\hspace{1.1cm}
\subfloat[Distorted mesh.]{\includegraphics[width=0.3\textwidth]{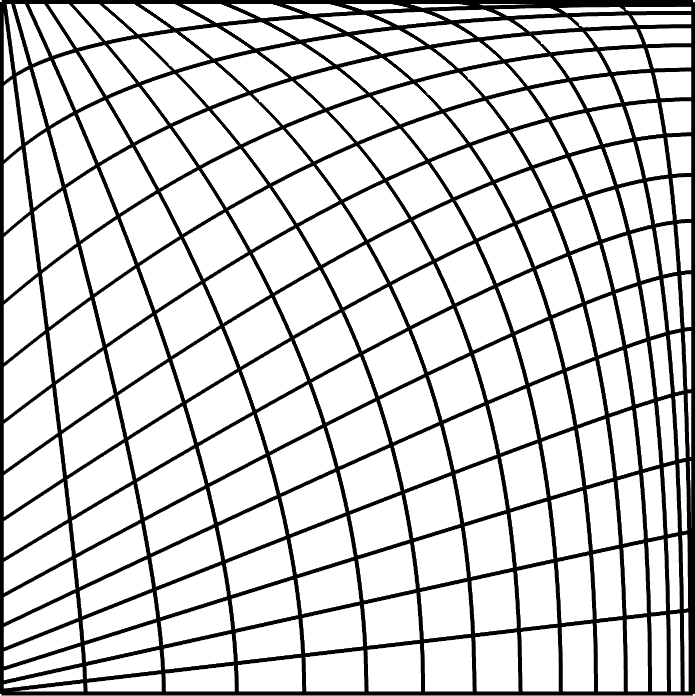}}
\caption{Meshes of 16$\times$16 elements used to discretize the square physical domain. The undistorted mesh is used to test the convergence in the parametric domain $\widehat{\Omega}$, and the distorted mesh to test the convergence in the physical domain $\Omega$.}
\label{fig:meshdist}
\end{figure}
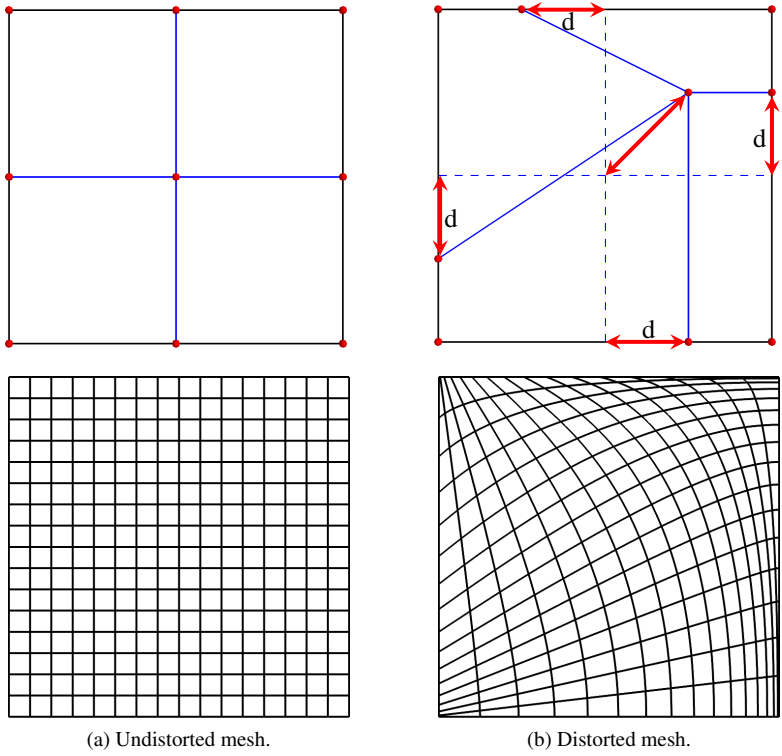

\begin{figure}[H]
\captionsetup{aboveskip=0.1cm}
\captionsetup{belowskip=-0.3cm}
\centering
\subfloat{
\includegraphics[trim={0.2cm 0 1.8cm 0.75cm},clip=true,width=7cm]{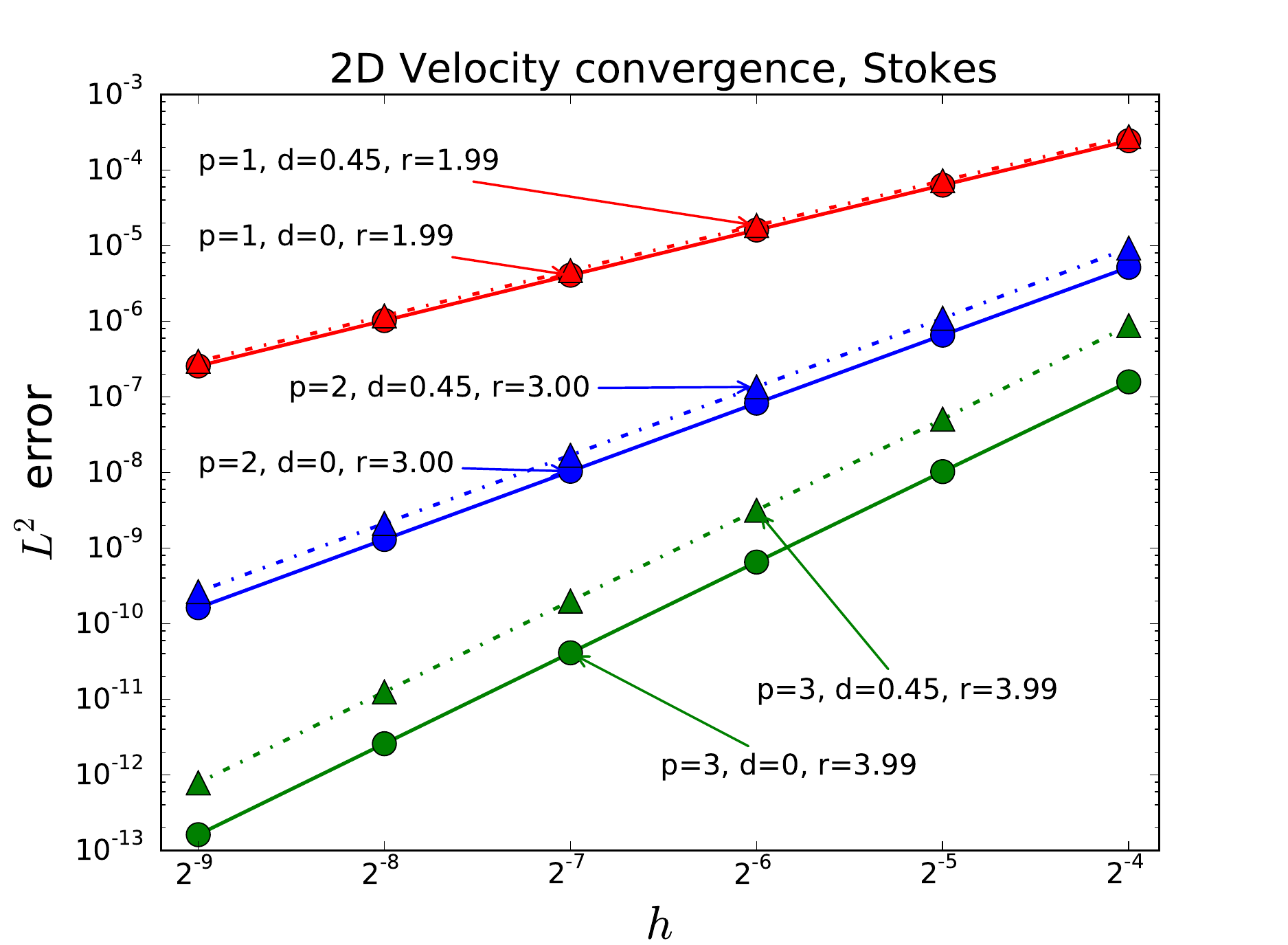}
}
\hspace{0.15cm}
\subfloat{
\includegraphics[trim={0.2cm 0 1.8cm 0.75cm},clip=true,width=7cm]{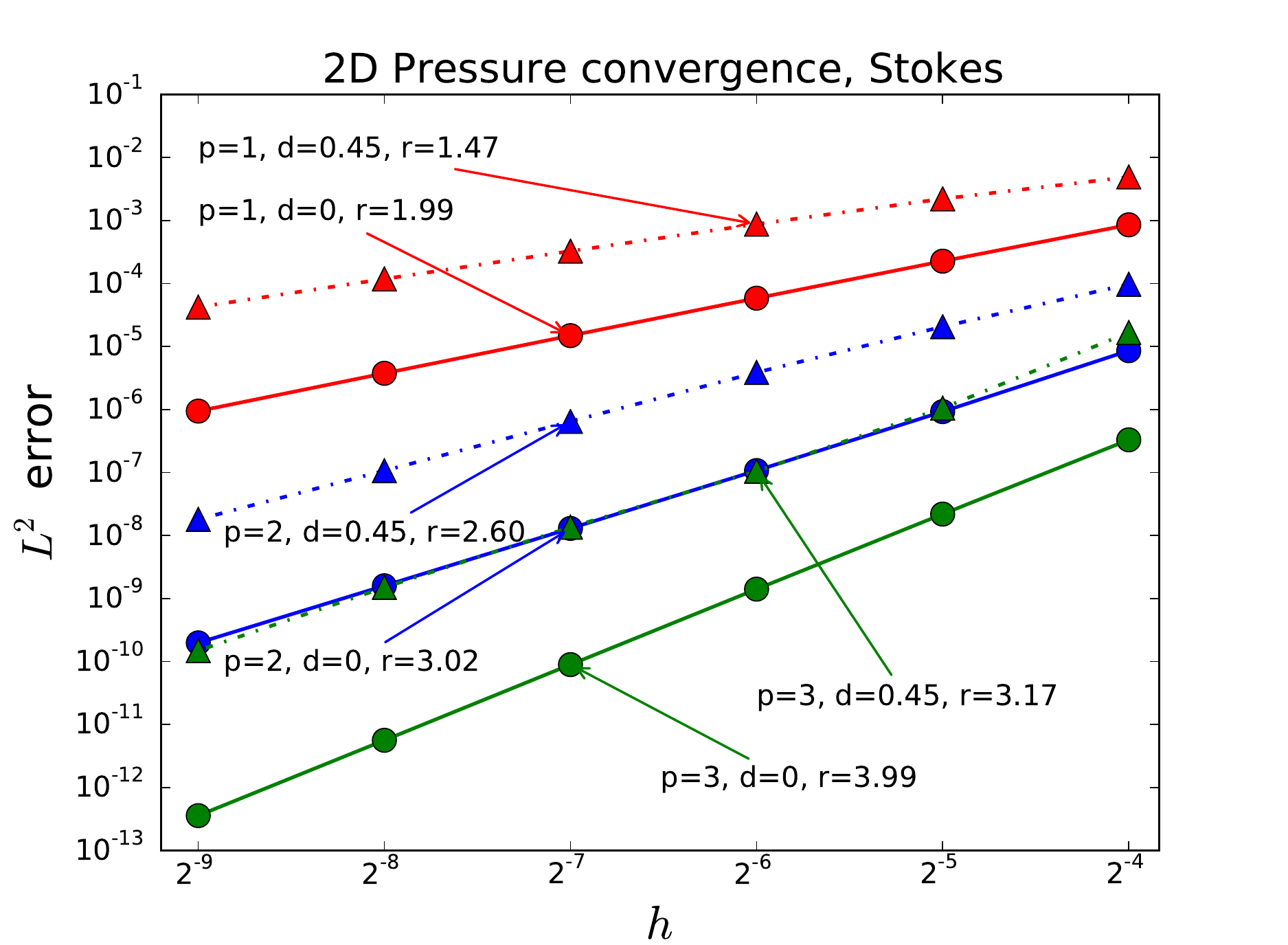}
}\\
\caption{Convergence test results for Stokes in the square problem.}
\label{fig:convsquareStokes}
\end{figure}

\begin{figure}[H]
\captionsetup{aboveskip=0.1cm}
\captionsetup{belowskip=-0.3cm}
\centering
\subfloat{
\includegraphics[trim={0.2cm 0 1.8cm 0.75cm},clip=true,width=7cm]{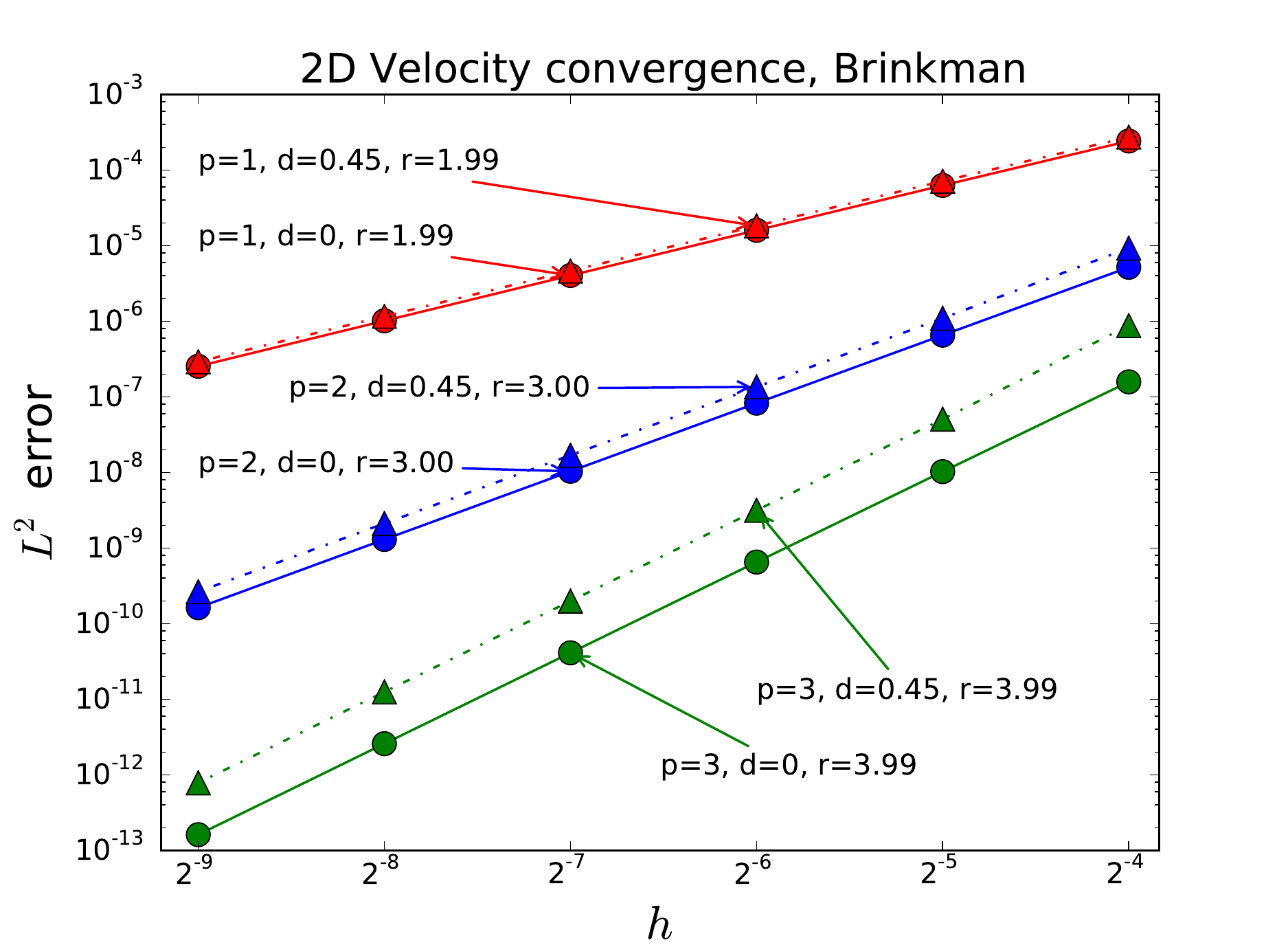}
}
\hspace{0.15cm}
\subfloat{
\includegraphics[trim={0.2cm 0 1.8cm 0.75cm},clip=true,width=7cm]{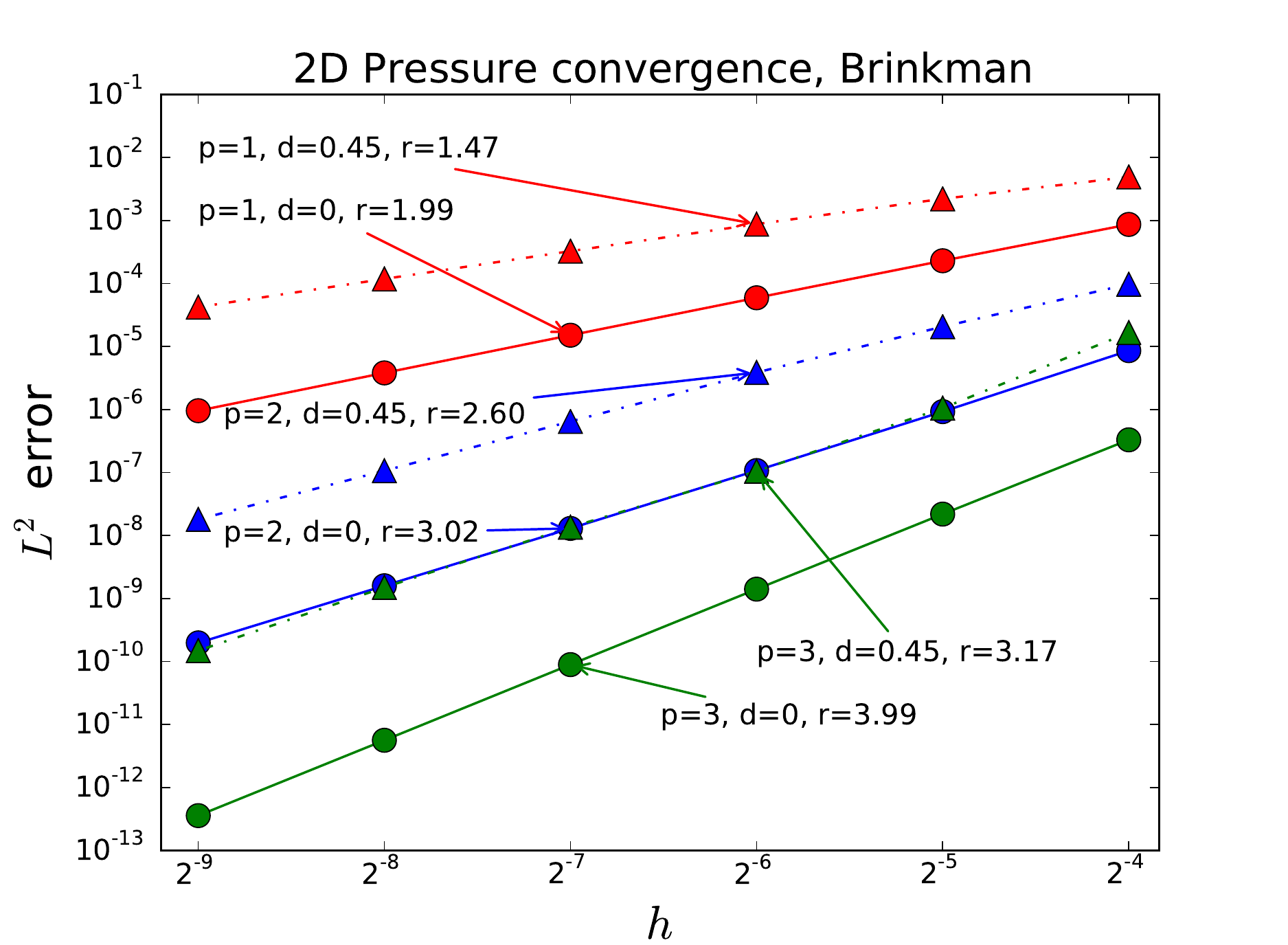}
}\\
\caption{Convergence test results for Brinkman in the square problem.}
\label{fig:convsquareBrinkman}
\end{figure}

\begin{figure}[H]
\captionsetup{aboveskip=0.1cm}
\captionsetup{belowskip=-0.3cm}
\centering
\subfloat{
\includegraphics[trim={0.2cm 0 1.8cm 0.75cm},clip=true,width=7cm]{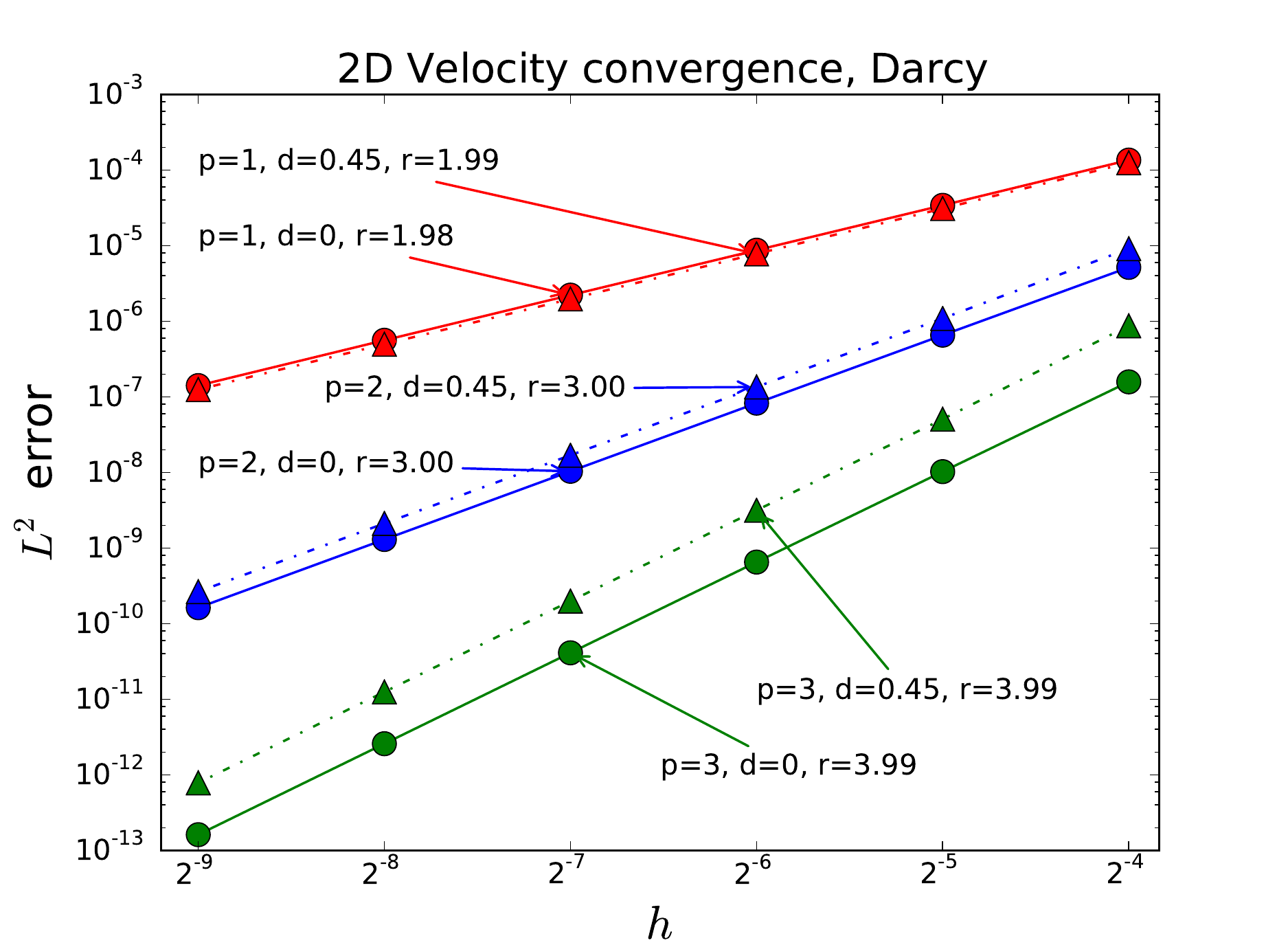}
}
\hspace{0.15cm}
\subfloat{
\includegraphics[trim={0.2cm 0 1.8cm 0.75cm},clip=true,width=7cm]{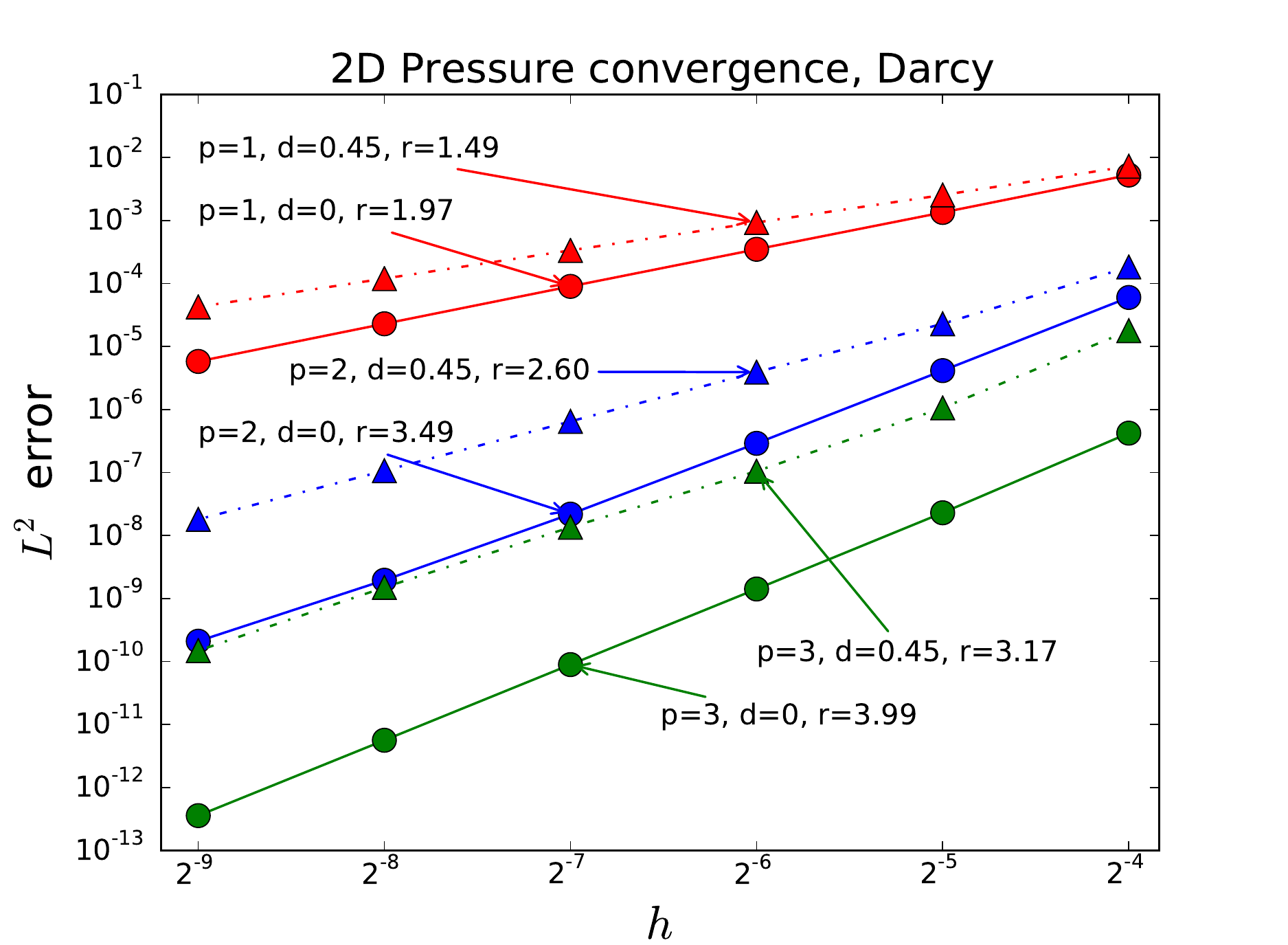}
}\\
\caption{Convergence test results for Darcy in the square problem.}
\label{fig:convsquareDarcy}
\end{figure}

\begin{figure}[H]
\captionsetup{aboveskip=0.1cm}
\captionsetup{belowskip=-0.3cm}
\centering
\subfloat{
\includegraphics[trim={0.2cm 0 1.8cm 0.75cm},clip=true,width=7cm]{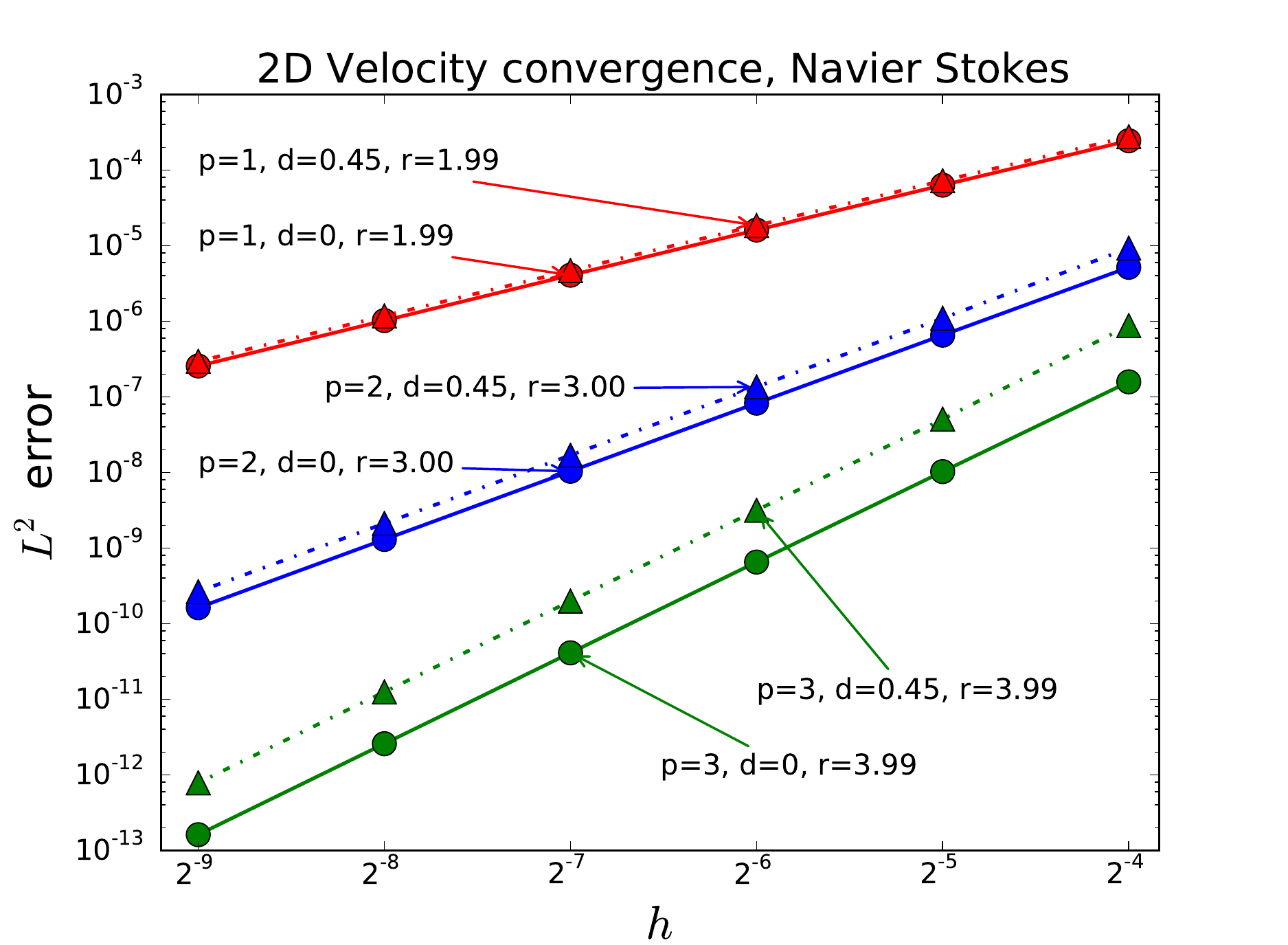}
}
\hspace{0.15cm}
\subfloat{
\includegraphics[trim={0.2cm 0 1.8cm 0.75cm},clip=true,width=7cm]{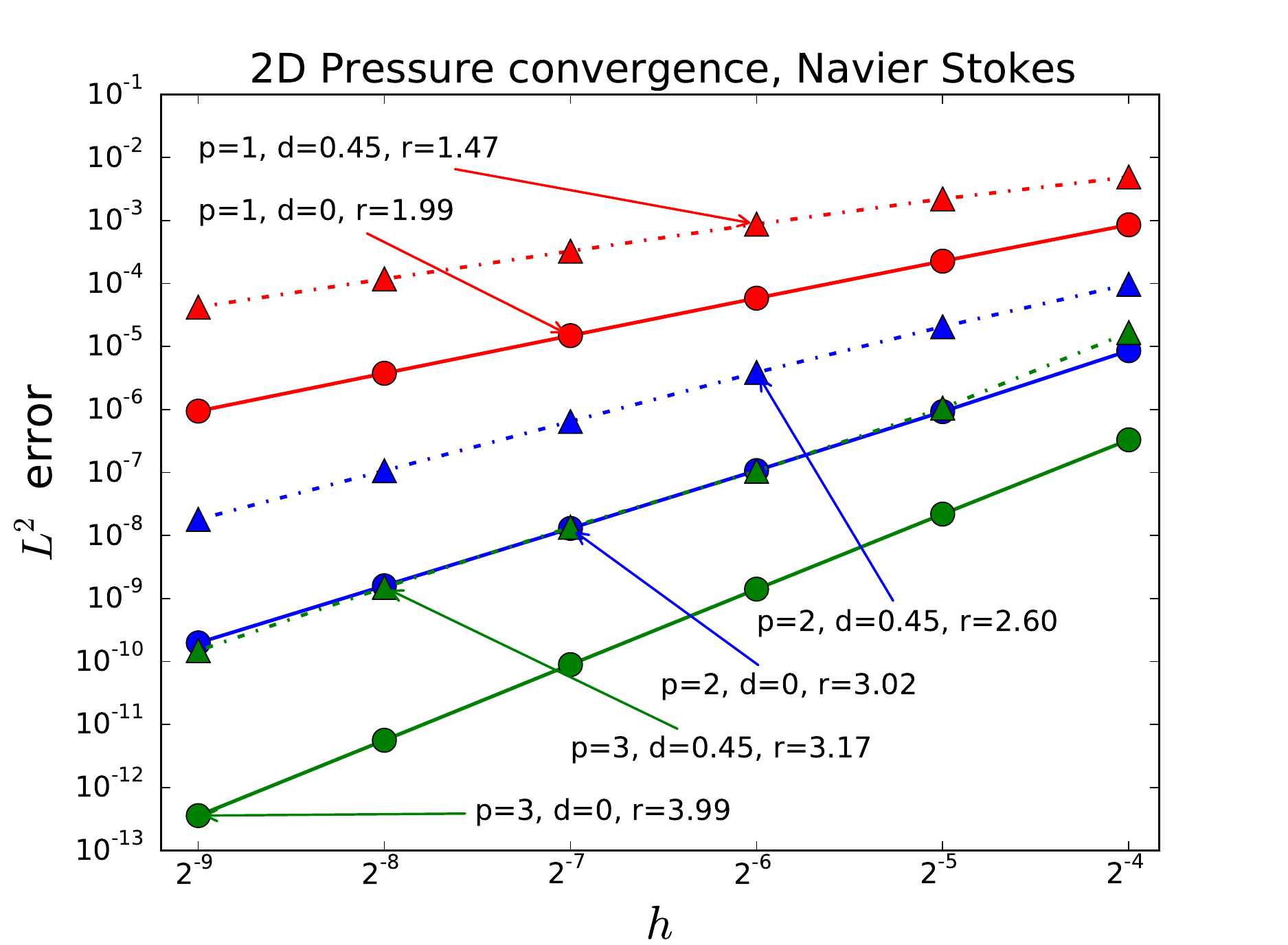}
}\\
\caption{Convergence test results for Navier-Stokes $Re\!=\!1$ in the square problem.}
\label{fig:convsquareNavierStokes}
\end{figure}

\begin{figure}[H]
\captionsetup{aboveskip=0.1cm}
\captionsetup{belowskip=-0.3cm}
\centering
\subfloat{
\includegraphics[trim={0.2cm 0 1.8cm 0.75cm},clip=true,width=7cm]{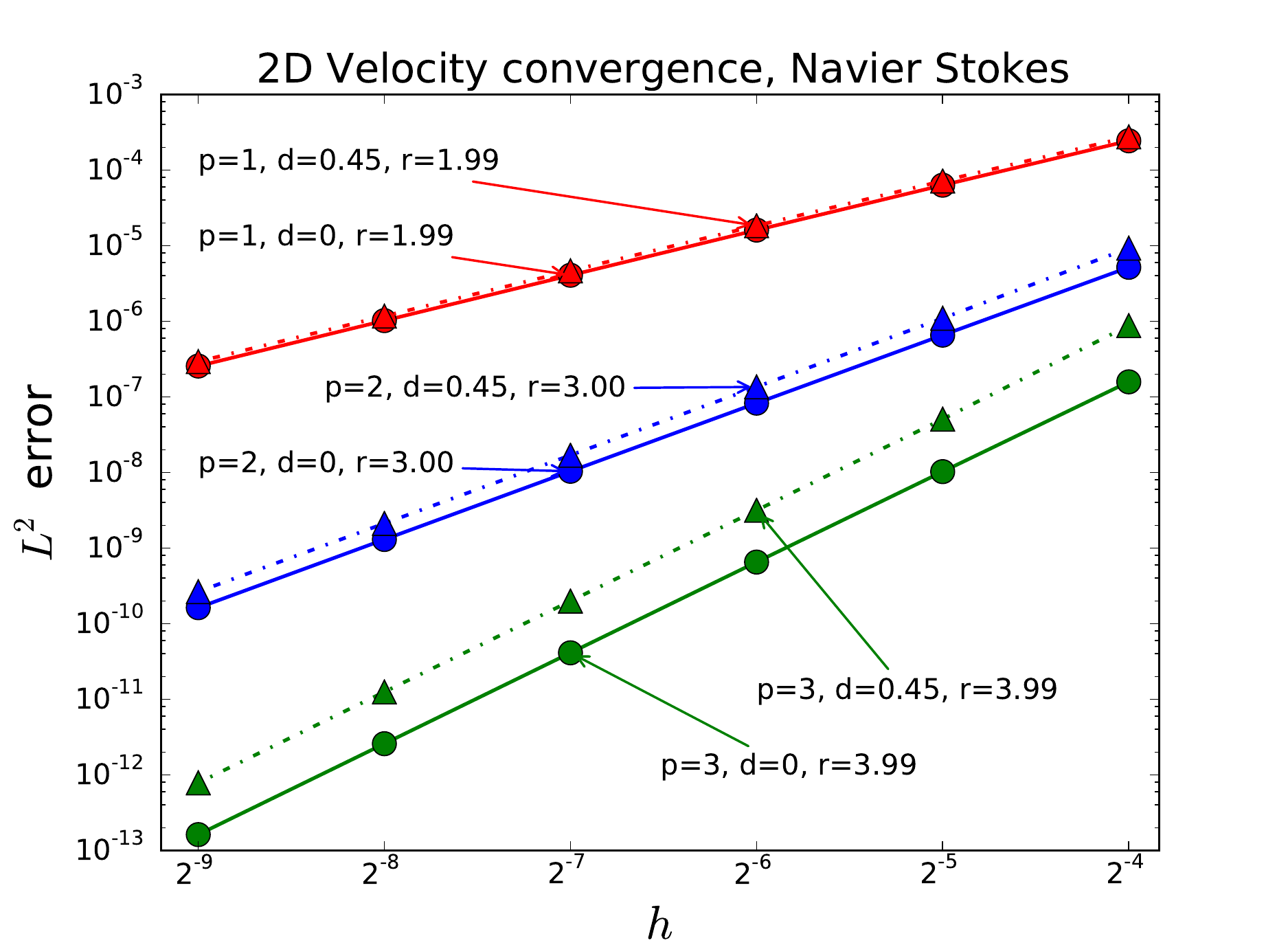}
}
\hspace{0.15cm}
\subfloat{
\includegraphics[trim={0.2cm 0 1.8cm 0.75cm},clip=true,width=7cm]{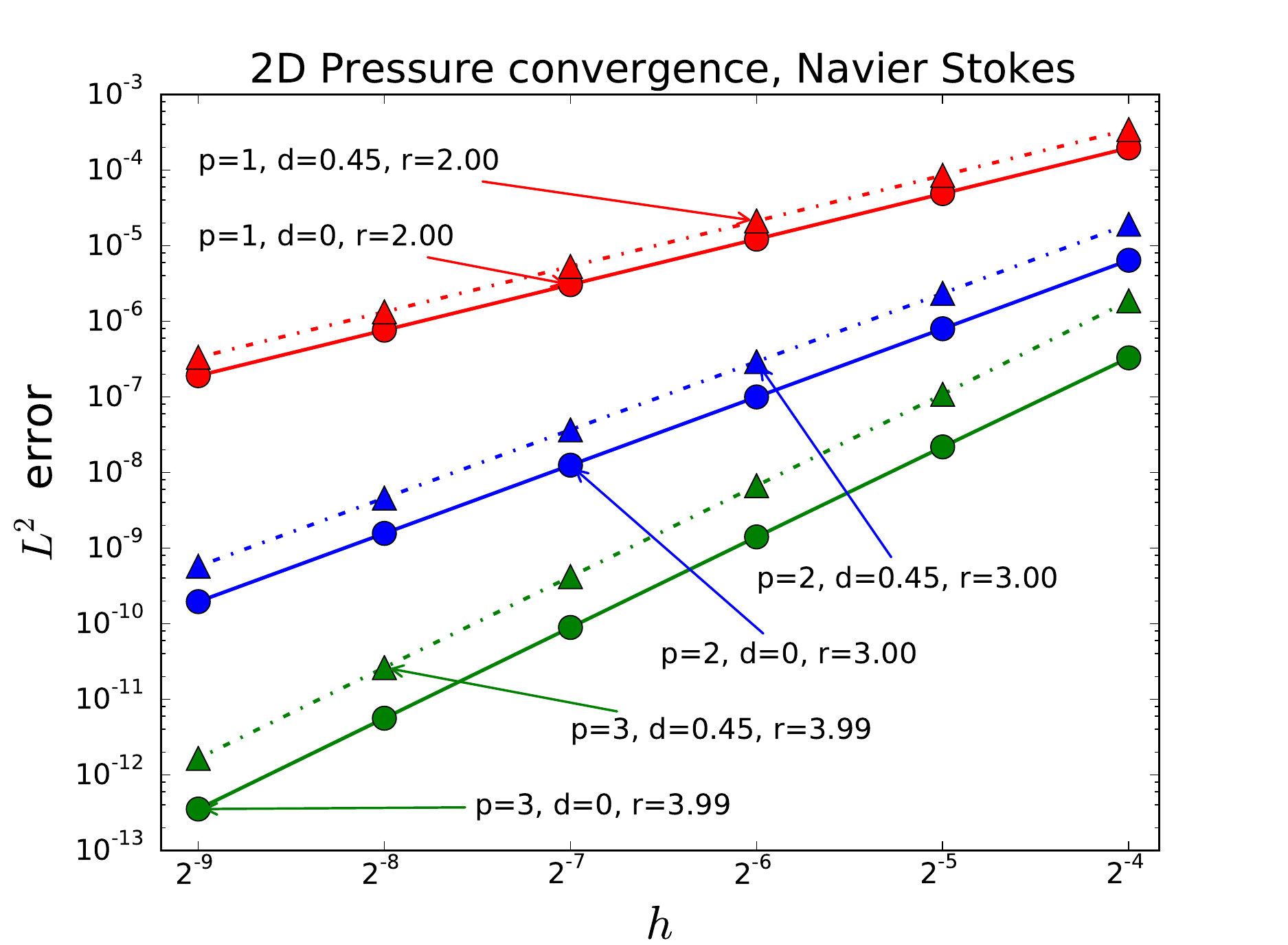}
}\\
\caption{Convergence test results for Navier-Stokes $Re\!=\!1000$ in the square problem.}
\label{fig:convsquareNavierStokes1000}
\end{figure}

Computing the $L^2$ norm of the error we verify the convergence rates of the method against the theoretical estimates. We solve for nested meshes from 16$\times$16 to 512$\times$512 elements, using the undistorted and distorted meshes seen in Figure \ref{fig:meshdist}, to prove convergence in the parametric and physical domains. The distorted mesh used for the convergence tests is created by moving the control points of a mesh with one element, polynomial order $p\!=\!2$ and continuity order $\varsigma\!=\!1$, a distance $d$ as shown in Figure \ref{fig:meshdist}(b), and then performing an $h$-refinement of the element. Results for three different polynomial orders with maximum continuity are shown in Figures \ref{fig:convsquareStokes} to \ref{fig:convsquareNavierStokes1000}, where the solid lines represent the results of the uniform meshes, and the dashed lines show the results of the distorted meshes. The asymptotic convergence rate $r$ is given for every mesh and discretization.

Figures \ref{fig:convsquareStokes} to \ref{fig:convsquareNavierStokes1000} show that the convergence rates $r$ for the error in the velocity is equal to $p\!+\!1$. These rates are not affected by the mesh distortion. Convergence rates for the error in the pressure when using uniform meshes are equal to $p\!+\!1$, while when using a distorted mesh the convergence rates deteriorate to $p$, being more notorious in the cases with high polynomial orders. The loss of convergence for the error in the pressure when using distorted meshes corroborates the a priori error estimates presented in \cite{EvansHughesStokes,EvansHughesNavierStokes}. These results show that the theory is sharp and that there is no superconvergence in the pressure.

\subsubsection{Remark.}
\begin{figure}[H]
\centering
\subfloat[Homogeneous reduction]{
\includegraphics[height=5.5cm]{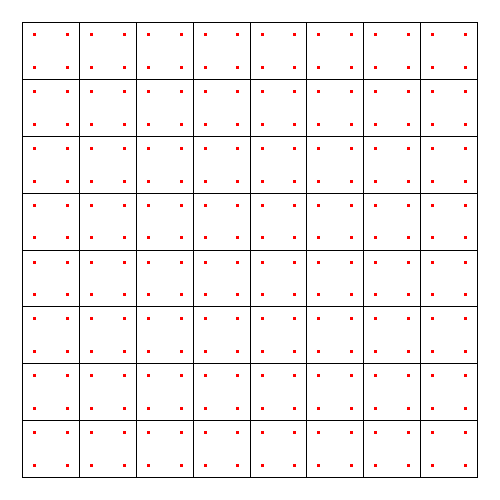}
}
\hspace{0.5cm}
\subfloat[Gradual reduction]{
\includegraphics[height=5.5cm]{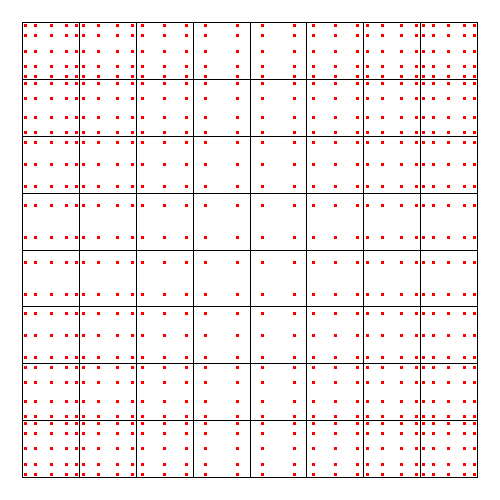}
}\\
\caption{Reduced quadrature schemes with a minimum of 2 points for a discretization using $p\!=\!3$, $\varsigma\!=\!2$ polynomials and a mesh of 8$\times$8 elements.}
\label{fig:reducequadrature}
\end{figure}

We evaluate the convergence rates using two reduced quadrature schemes, one in which we keep the exact quadrature of $p\!+\!2$ points for the elements at the boundaries, and gradually reduce the number of quadrature points by one to the contiguous elements, as they approach the center of the domain, until they reach a given minimum number of quadrature points per direction as shown in Figure \ref{fig:reducequadrature} (b). We also consider a homogeneous reduction of quadrature points. Both reduction schemes using $p\!+\!1$ quadrature points in every direction produce the same convergence rates as the exact quadrature, and no deterioration on the convergence constant. When both schemes reduce the number of quadrature points to $p$, velocity convergence remains equal to the exact quadrature, but pressure convergence rate and constant start decreasing. The matrix is not invertible in the case with the lowest order discretization ($p\!=\!1$, $\varsigma\!=\!0$) and the homogeneous reduction scheme using $p$ quadrature points. Any reduction beyond $p$ quadrature points deteriorates the convergence of both velocity and pressure.

\subsection{Two-dimensional lid-driven square cavity}
\begin{figure}
\centering
\begin{tikzpicture}[scale=1.05, transform shape]

\node[inner sep=0pt] (russell) at (0  ,0 ) {\includegraphics[height=6cm]{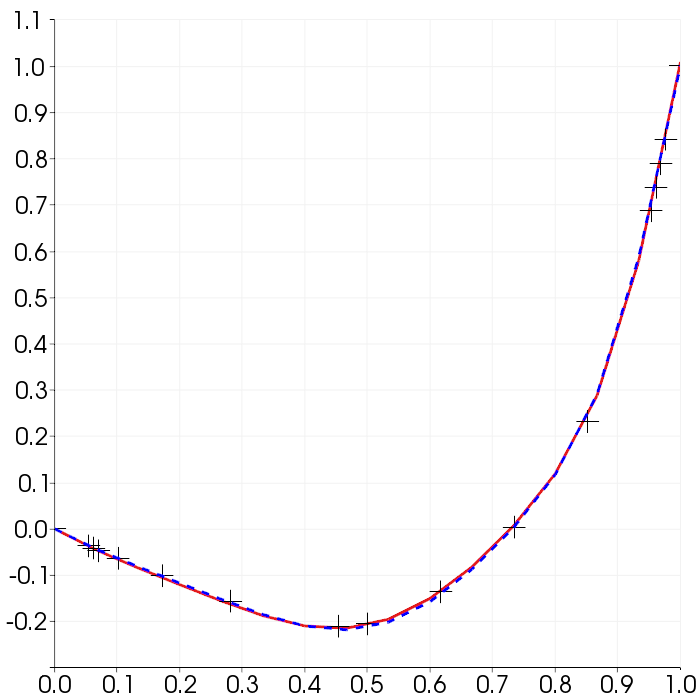}};

\draw (0.3,-3.2) node {$y$};
\draw (-3.2, 0.3) node[rotate=90] {$u$};

\draw (0,1  ) node (1u){\scriptsize $d\!=\!0.45$};
\draw (0,0.5) node (2u){\scriptsize $d\!=\!0.0 $};
\draw (0,0  ) node (3u){\scriptsize Ghia};

\draw[>=stealth,->,       red  ,line width=1pt](1u)--(2.5, 1  );
\draw[>=stealth,->,dashed,blue ,line width=1pt](2u)--(2.3, 0.2);
\draw[>=stealth,->,dotted,black,line width=1pt](3u)--(2  ,-0.6);

\node[inner sep=0pt] (russell) at (7  ,0 ) {\includegraphics[height=6cm]{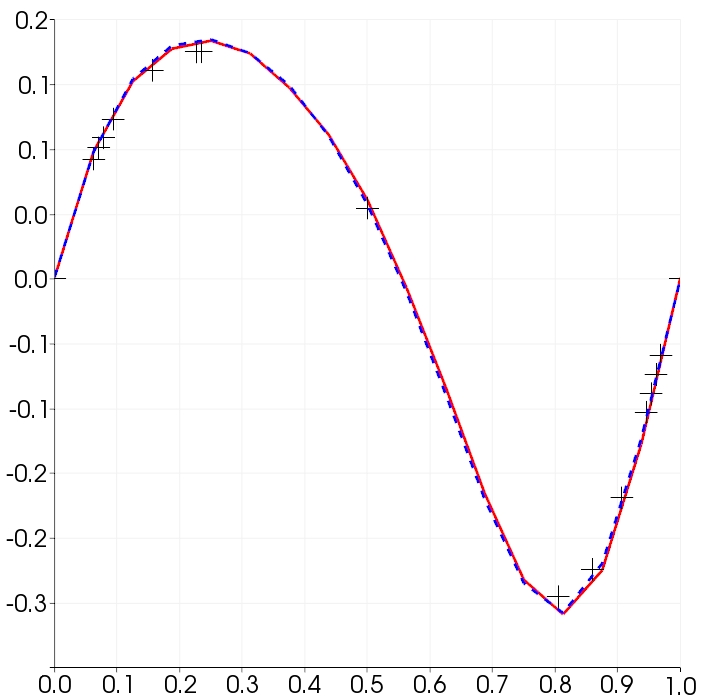}};

\draw (7.3,-3.2) node {$x$};
\draw (3.8, 0.3) node[rotate=90] {$v$};

\draw (6, 0  ) node (1v){\scriptsize $d\!=\!0.45$};
\draw (6,-0.5) node (2v){\scriptsize $d\!=\!0.0 $};
\draw (6,-1  ) node (3v){\scriptsize Ghia};

\draw[>=stealth,->,       red  ,line width=1pt](1v)--(7.7, 0  );
\draw[>=stealth,->,dashed,blue ,line width=1pt](2v)--(8.1,-1.1);
\draw[>=stealth,->,dotted,black,line width=1pt](3v)--(8.8,-2.1);

\end{tikzpicture}
\caption{Comparison of vertical and horizontal velocities along the horizontal and vertical centerlines, respectively, for uniform~($d\!=\!0.0$) and distorted~($d\!=\!0.45$) meshes, solving the Navier-Stokes problem in a unitary square with $Re\!=\!100$, using $h\!=\!$~\sfrac{1}{16} and $p\!=\!1$, $\varsigma\!=\!0$. Our numerical results with the two different meshes compare favorably to Ghia's benchmark \cite{Ghia}.}
\label{fig:Re100}
\end{figure}
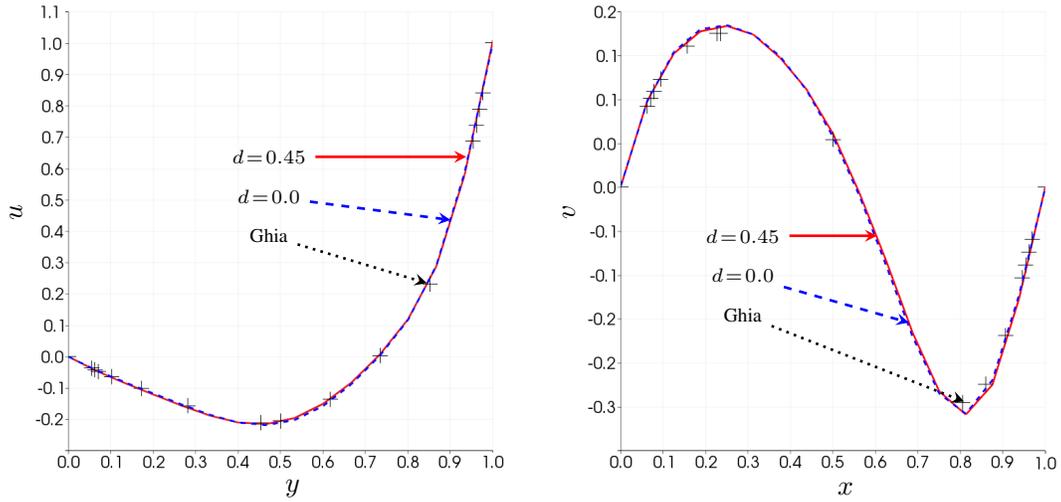

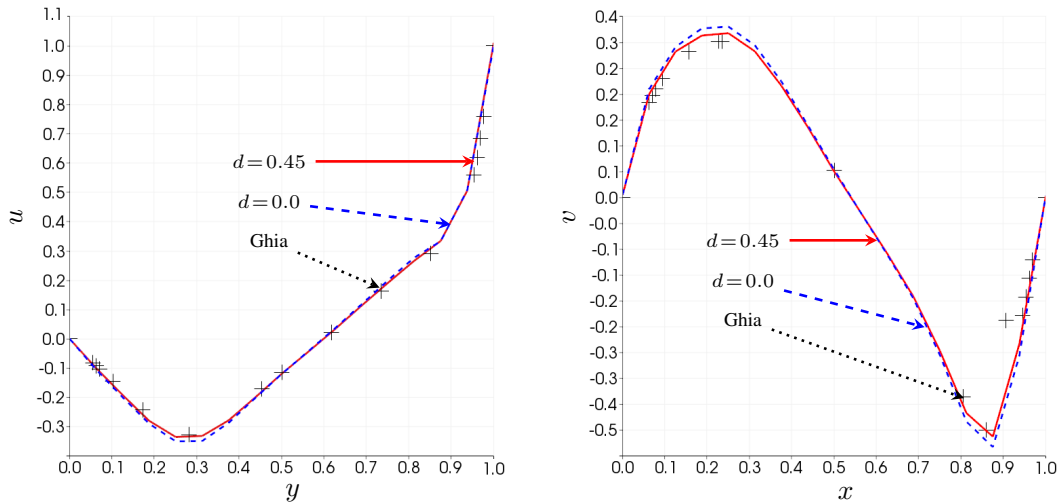
\begin{figure}
\centering
\begin{tikzpicture}[scale=1.05, transform shape]

\node[inner sep=0pt] (russell) at (0  ,0 ) {\includegraphics[height=6cm]{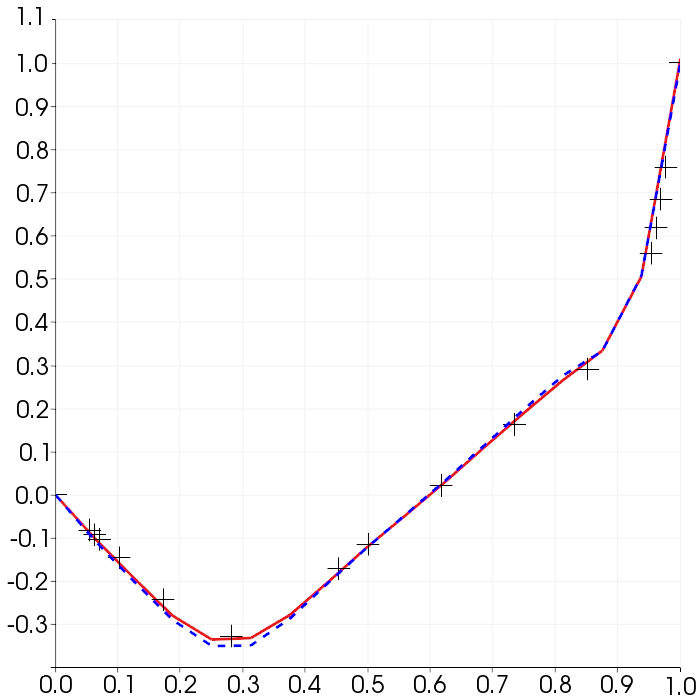}};

\draw (0.3,-3.2) node {$y$};
\draw (-3.2, 0.3) node[rotate=90] {$u$};

\draw (0,1  ) node (1u){\scriptsize $d\!=\!0.45$};
\draw (0,0.5) node (2u){\scriptsize $d\!=\!0.0 $};
\draw (0,0  ) node (3u){\scriptsize Ghia};

\draw[>=stealth,->,       red  ,line width=1pt](1u)--(2.6, 1  );
\draw[>=stealth,->,dashed,blue ,line width=1pt](2u)--(2.3, 0.2);
\draw[>=stealth,->,dotted,black,line width=1pt](3u)--(1.4,-0.6);

\node[inner sep=0pt] (russell) at (7  ,0 ) {\includegraphics[height=6cm]{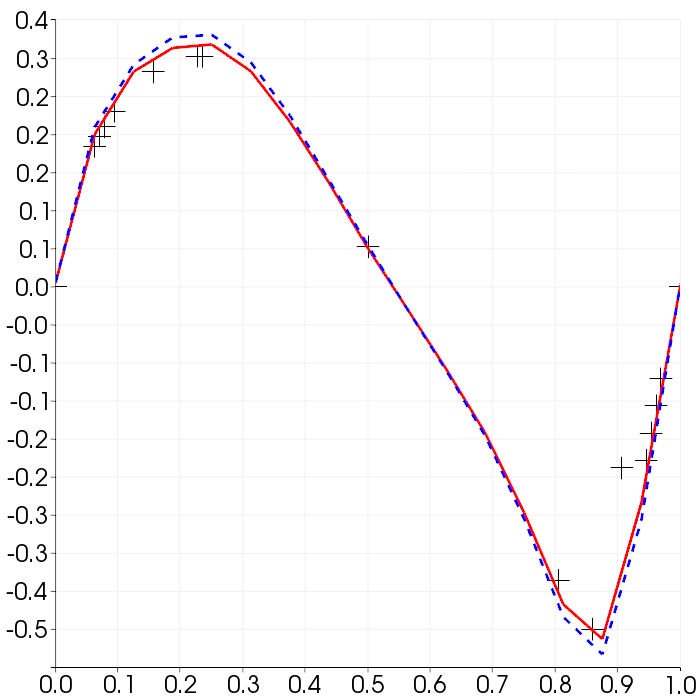}};

\draw (7.3,-3.2) node {$x$};
\draw (3.8, 0.3) node[rotate=90] {$v$};

\draw (6, 0  ) node (1v){\scriptsize $d\!=\!0.45$};
\draw (6,-0.5) node (2v){\scriptsize $d\!=\!0.0 $};
\draw (6,-1  ) node (3v){\scriptsize Ghia};

\draw[>=stealth,->,       red  ,line width=1pt](1v)--(7.7, 0  );
\draw[>=stealth,->,dashed,blue ,line width=1pt](2v)--(8.3,-1.1);
\draw[>=stealth,->,dotted,black,line width=1pt](3v)--(8.8,-2.0);

\end{tikzpicture}
\caption{Comparison of vertical and horizontal velocities along the horizontal and vertical centerlines, respectively, for uniform~($d\!=\!0.0$) and distorted~($d\!=\!0.45$) meshes, solving the Navier-Stokes problem in a unitary square with $Re\!=\!400$, using $h\!=\!$~\sfrac{1}{16} and $p\!=\!1$, $\varsigma\!=\!0$. Our numerical results with the two different meshes compare favorably to Ghia's benchmark \cite{Ghia}.}
\label{fig:Re400}
\end{figure}

We solve the two-dimensional lid-driven cavity test for the Stokes and Navier-Stokes equations, using the same set of nested meshes from 16 to 512 elements per side, as in the previous example to compare the solutions. The solutions found for the Stokes problem are compared to a spectral approach \cite{Botella} in Table \ref{tab:vort}, comparing the value of the vorticity at a specified point near the top right corner $(\bx\!=\!(1,0.95))$, for discretizations using $p\!=\!1,2,3$ and maximum continuity. The solution for the Navier-Stokes problem uses two different Reynolds numbers $Re\!=\!100$ and $Re\!=\!400$, and we compare with the results presented by Ghia in \cite{Ghia} and the spectral approach \cite{Botella}. Tables \ref{tab:maxvelRe100} and \ref{tab:maxvelRe400} compare the value and position of the minimum horizontal velocity along the vertical centerline $(x\!=\!0.5)$, and the value and position of the minimum and maximum vertical velocity along the horizontal centerline $(y\!=\!0.5)$. Figures \ref{fig:Re100} and \ref{fig:Re400} illustrate the effect of the mesh distortion for the case of Navier-Stokes, where the results found with the coarsest mesh ($h\!=\!$ \sfrac{1}{16}, $p\!=\!1$, $\varsigma\!=\!0$) without distortion, are compared to the results found using a distorted mesh ($d\!=\!0.45$) with the same discretization.

\begin{table}[H]
\footnotesize
\captionof{table}{Convergence of the vorticity @$x\!=\!(1,0.95)$ for the Stokes problem.} 
\label{tab:vort} 
\centering
\small
\begin{tabular}{C{1.5cm}C{0.5cm}C{1.5cm}C{1.5cm}C{1.6cm}C{1.6cm}C{1.6cm}C{1.6cm}C{0.1cm}}
\cline{1-8}
\multicolumn{1}{ c }{\multirow{2}{*}{Method} } &
\multicolumn{1}{ c }{\multirow{2}{*}{$h$}    } &
\multicolumn{2}{ c }{$p\!=\!1$, $\varsigma\!=\!0$}&\multicolumn{2}{c}{$p\!=\!2$, $\varsigma\!=\!1$}&\multicolumn{2}{c}{$p\!=\!3$, $\varsigma\!=\!2$}&\\[0.15cm]\cline{3-8}
\multicolumn{1}{ c }{}    						&& $d\!=\!0$ & $d\!=\!0.45$ & $d\!=\!0$ & $d\!=\!0.45$ & $d\!=\!0$ & $d\!=\!0.45$ &\\[0.15cm] \cline{1-8}
\multicolumn{1}{ c }{\multirow{6}{*}{IGA}}
		 						& \sfrac{1}{16}  & -0.528094 &  1.273373 & 12.947509 & 11.517694 & 32.790408 & 22.523328 &\\[0.15cm]
		 						& \sfrac{1}{32}  & 18.075800 &  9.838386 & 33.277310 & 23.387877 & 22.522894 & 29.289110 &\\[0.10cm]
		 						& \sfrac{1}{64}  & 19.186815 & 17.129618 & 35.017081 & 27.772500 & 30.291823 & 28.185708 &\\[0.10cm]
		 						& \sfrac{1}{128} & 23.479465 & 22.074412 & 25.848790 & 27.650488 & 29.324554 & 27.356594 &\\[0.10cm]
		 						& \sfrac{1}{256} & 25.425581 & 24.767150 & 27.342346 & 27.378879 & 27.642689 & 27.286833 &\\[0.10cm]
		 						& \sfrac{1}{512} & 26.371713 & 26.060740 & 27.294087 & 27.303925 & 27.278365 & 27.279689 &\\[0.10cm]\cline{1-8}
\multicolumn{1}{ c }{Spectral \cite{Botella}}	&& 27.27901  & -         & 27.27901  & -         & 27.27901  & -         &\\[0.10cm]\cline{1-8}
\end{tabular}
\end{table}

Table \ref{tab:vort} shows how the values found for the vorticity at a point that is near to the discontinuity of the velocity, located at the top right corner of the cavity, converge to the results found with a highly accurate spectral method of order 48 \cite{Botella}. Results found with a discretization of $p\!=\!1$ and $\varsigma\!=\!0$ converge to the benchmark solution at a slower rate than the other discretizations, a finer mesh than the ones tested here is needed to resolve the corner singularity in this case. Higher order discretizations converge to within the first two significant digits of the spectral solution with a mesh of $256\times256$ elements, and within four significant digits with $p\!=\!3$ and $\varsigma\!=\!2$ when using a mesh of $512\times512$ elements. Results found with the distorted mesh converge to those of the uniform mesh showing the robustness of the discretization used for the velocity field.

\noindent
\begin{table}[H]
\captionsetup{belowskip=-0.15cm}
\captionof{table}{Velocity extrema for the Navier-Stokes problem ($Re\!=\!100$) using uniform meshes ($d\!=\!0$).} 
\label{tab:maxvelRe100} 
\centering
\small 
\begin{tabular}{C{2cm}C{0.5cm}C{1.5cm}C{1.3cm}C{1.5cm}C{1.3cm}C{1.5cm}C{1.3cm}C{0.1cm}}
\cline{1-8}
Discretization&$h$									  &$u_{min}$	&$y_{min}$	&$v_{min}$		&$x_{min}$	&$v_{max}$		&$x_{max}$	&\\[0.15cm]\cline{1-8}
\multicolumn{1}{ c }{\multirow{2}{*}{$p\!=\!1$, $\varsigma\!=\!0$}}
							& \sfrac{1}{16}  		  &-0.2201506	&0.43750	&-0.2605222		&0.81249 	&0.1851086		&0.25000 	&\\[0.15cm]
							& \sfrac{1}{256} 		  &-0.2140707	&0.45703	&-0.2538092		&0.80859	&0.1795948		&0.23828	&\\[0.10cm]\cline{1-8}
\multicolumn{1}{ c }{\multirow{2}{*}{$p\!=\!2$, $\varsigma\!=\!1$}}
							& \sfrac{1}{16}  		  &-0.2142675	&0.45766	&-0.2537870		&0.81140	&0.1797504		&0.23706	&\\[0.10cm]
							& \sfrac{1}{256} 		  &-0.2140423	&0.45808	&-0.2538029		&0.81042	&0.1795728 		&0.23698	&\\[0.10cm]\cline{1-8}
\multicolumn{1}{ c }{\multirow{2}{*}{$p\!=\!3$, $\varsigma\!=\!2$}}
							& \sfrac{1}{16}  		  &-0.2140613	&0.45808	&-0.2539128		&0.81026	&0.1796009		&0.23679	&\\[0.10cm]
							& \sfrac{1}{256} 		  &-0.2140423	&0.45808	&-0.2538029		&0.81042	&0.1795728		&0.23698	&\\[0.10cm]\cline{1-8}
\multicolumn{2}{ c }{Spectral \cite{Botella}}		  &-0.2140424   &0.4581		&-0.2538030 	&0.8104 	&0.1795728 		&0.237  	&\\[0.10cm]
\multicolumn{2}{ c }{Finite differences \cite{Ghia}}  &-0.21090 	&0.4531		&-0.24533 		&0.8047 	&0.17527 		&0.2344 	&\\[0.10cm]\cline{1-8}
\end{tabular}
\end{table}

\noindent
\begin{table}[H]
\captionsetup{belowskip=-0.15cm}
\captionof{table}{Velocity extrema for the Navier-Stokes problem ($Re\!=\!400$) using uniform meshes ($d\!=\!0$).} 
\label{tab:maxvelRe400} 
\centering
\small
\begin{tabular}{C{2cm}C{0.5cm}C{1.5cm}C{1.3cm}C{1.5cm}C{1.3cm}C{1.5cm}C{1.3cm}C{0.1cm}}
\cline{1-8}
Discretization&$h$									  &$u_{min}$	&$y_{min}$	&$v_{min}$		&$x_{min}$	&$v_{max}$		&$x_{max}$	&\\[0.15cm]\cline{1-8}
\multicolumn{1}{ c }{\multirow{2}{*}{$p\!=\!1$, $\varsigma\!=\!0$}}
							& \sfrac{1}{16}  		  &-0.3523864 	&0.25000 	&-0.4920310 	&0.87499 	&0.3312674 		&0.24999	&\\[0.15cm]
							& \sfrac{1}{256} 		  &-0.3288927 	&0.28124 	&-0.4542830 	&0.86328 	&0.3039886 		&0.22656	&\\[0.10cm]\cline{1-8}
\multicolumn{1}{ c }{\multirow{2}{*}{$p\!=\!2$, $\varsigma\!=\!1$}}
							& \sfrac{1}{16}  		  &-0.3337101 	&0.28140 	&-0.4547631 	&0.85979 	&0.3078021 		&0.22429	&\\[0.10cm]
							& \sfrac{1}{256} 		  &-0.3287303 	&0.28002 	&-0.4540652 	&0.86220	&0.3038326 		&0.22530	&\\[0.10cm]\cline{1-8}
\multicolumn{1}{ c }{\multirow{2}{*}{$p\!=\!3$, $\varsigma\!=\!2$}}
							& \sfrac{1}{16}  		  &-0.3298355 	&0.28047 	&-0.4550065 	&0.86134 	&0.3047172 		&0.22599	&\\[0.10cm]
							& \sfrac{1}{256} 		  &-0.3287302 	&0.28002 	&-0.4540654		&0.86221 	&0.3038325 		&0.22530	&\\[0.10cm]\cline{1-8}
\multicolumn{2}{ c }{Finite differences \cite{Ghia}}  &-0.32726 	&0.2813 	&-0.44993  		&0.8594  	&0.30203  		&0.2266 	&\\[0.10cm]\cline{1-8}
\end{tabular}
\end{table}

Tables \ref{tab:maxvelRe100} and \ref{tab:maxvelRe400} compare the values of the velocity extrema when solving the Navier-Stokes system with a Reynolds number of one hundred, and four hundred, respectively, when using undistorted meshes, against the solution using a spectral method of order 96 \cite{Botella} for the case of $Re\!=\!100$ and a second order upwind finite differences method using $129\times129$ points \cite{Ghia} for both Reynolds numbers. We compare the values of the maximum horizontal velocity and its position along the vertical center line, and the values of the maximum and minimum vertical velocity along the horizontal center line, for the coarsest ($h\!=\!$ \sfrac{1}{16}) and the finest meshes ($h\!=\!$ \sfrac{1}{256}) used, and discretizations of $p\!=\!1,2,3$ and maximum continuity. For both Reynolds numbers considered, all the results are reasonably close to the benchmark values, with the exception of the coarsest mesh when using the $p\!=\!1$, $\varsigma\!=\!0$ discretization, which is the only that differs noticeably from the others. When using discretizations of $p\!>\!1$ the differences between the results from the coarsest and finest meshes become small, suggesting that a high order discretization with a coarse mesh may be enough to capture most of the features of the flow inside the domain.

Figures \ref{fig:Re100} and \ref{fig:Re400} illustrate the effect of the mesh distortion when using a $p\!=\!1, \varsigma\!=\!0$ discretization and the coarsest mesh, by comparing the results found when solving the Navier-Stokes problem for two Reynolds numbers ($Re\!=\!100$ and $Re\!=\!400$) with the undistorted and the distorted mesh against the results found with the second order upwind finite differences method using $129\times129$ points \cite{Ghia}. These comparisons show that the discretization is robust with respect to the mesh distortion and that even the coarsest discretization provides a fair approximation to the benchmark.

\setcounter{subfigure}{0}

\subsection{Cylindrical Couette flow}
\begin{figure}
\centering
\subfloat[Velocity magnitude.]{
\includegraphics[height=5cm]{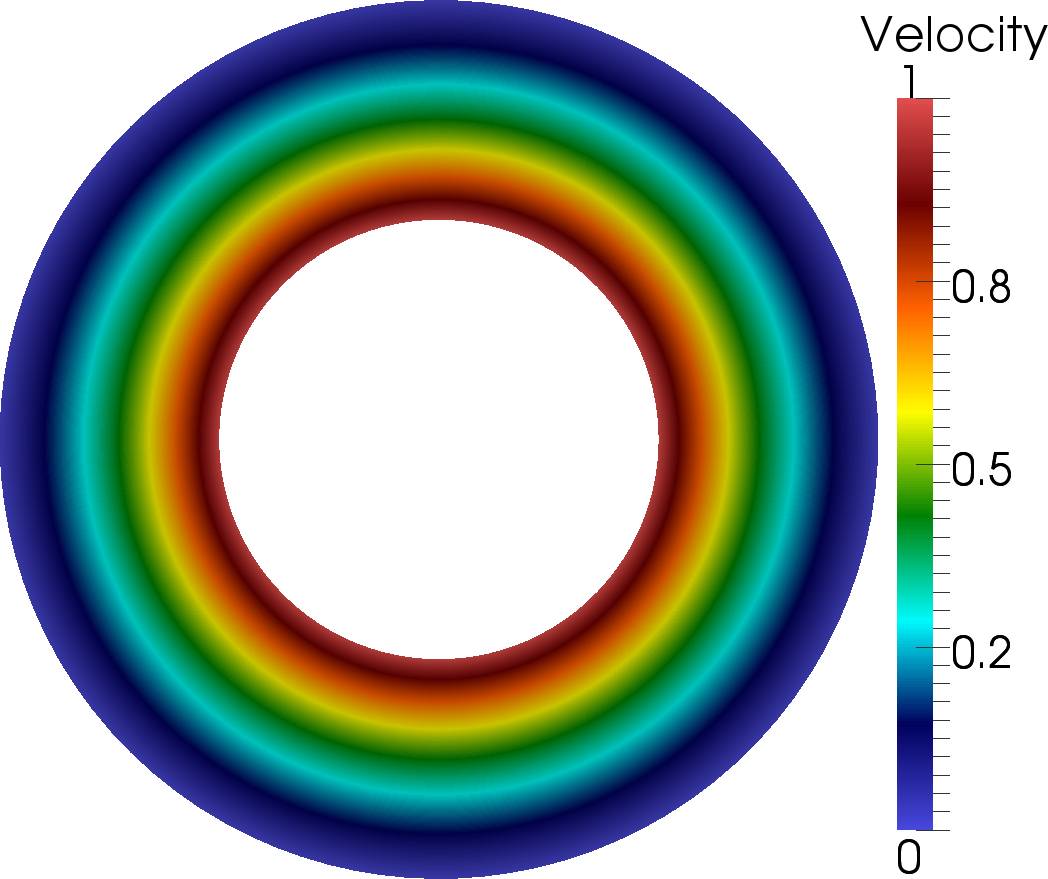}
}
\hspace{0.7cm}
\subfloat[Pressure.]{
\includegraphics[height=5cm]{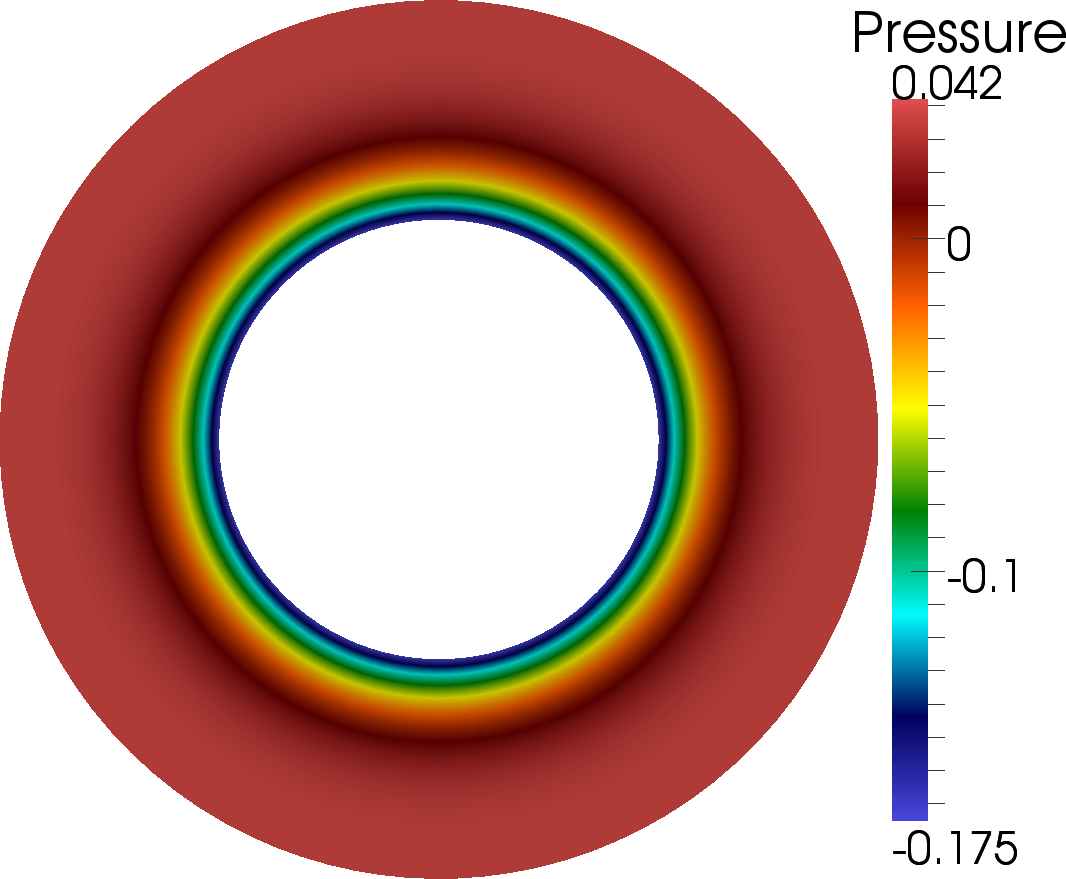}
}
\caption{Analytical solution for the Couette flow for the Navier-Stokes system.}
\label{fig:resultsannulusanaly}
\end{figure}

We present results for a Couette flow in an annulus to test convergence in a physical domain different than a square. The flow is driven by a boundary condition of a unitary tangential velocity on the inner face of the annulus. We test the solutions for the Stokes and Navier-Stokes problems for two different Reynold numbers $Re\!=\!1$ and $Re\!=\!100$. The analytical solution for the velocity when considering $Da\!=\!0$, as shown in Figure \ref{fig:resultsannulusanaly}(a) is given by the following expression:
\begin{align*}
\overline{\bu}=&\begin{bmatrix}
 (Ar+\frac{B}{r})\sin(\theta)\\[0.15cm]
 (Ar+\frac{B}{r})\cos(\theta)
\end{bmatrix}
\end{align*}
where $r$ and $\theta$ correspond to the polar coordinates, and
\begin{align*}
A=-\frac{U\,\delta^2}{r_{in}(1-\delta^2)},\quad
B=\frac{U\,r_{in}}{(1-\delta^2)},\quad
\delta=\frac{r_{in}}{r_{out}}.
\end{align*}

The analytical solution for the pressure for the case of Stokes is equal to zero in all the domain, and for the Navier-Stokes case as shown in Figure \ref{fig:resultsannulusanaly}(b), is given by the following expression:
\begin{align*}
\frac{\partial \overline{p}}{\partial r}=\frac{\left(Ar+\frac{B}{r}\right)^2}{r}
\end{align*}

\begin{figure}[H]
\captionsetup{aboveskip=0.1cm}
\captionsetup{belowskip=-0.3cm}
\centering
\subfloat{
\includegraphics[trim=3.4cm 2.9cm 2.9cm 2.9cm, clip=true,height=5cm]{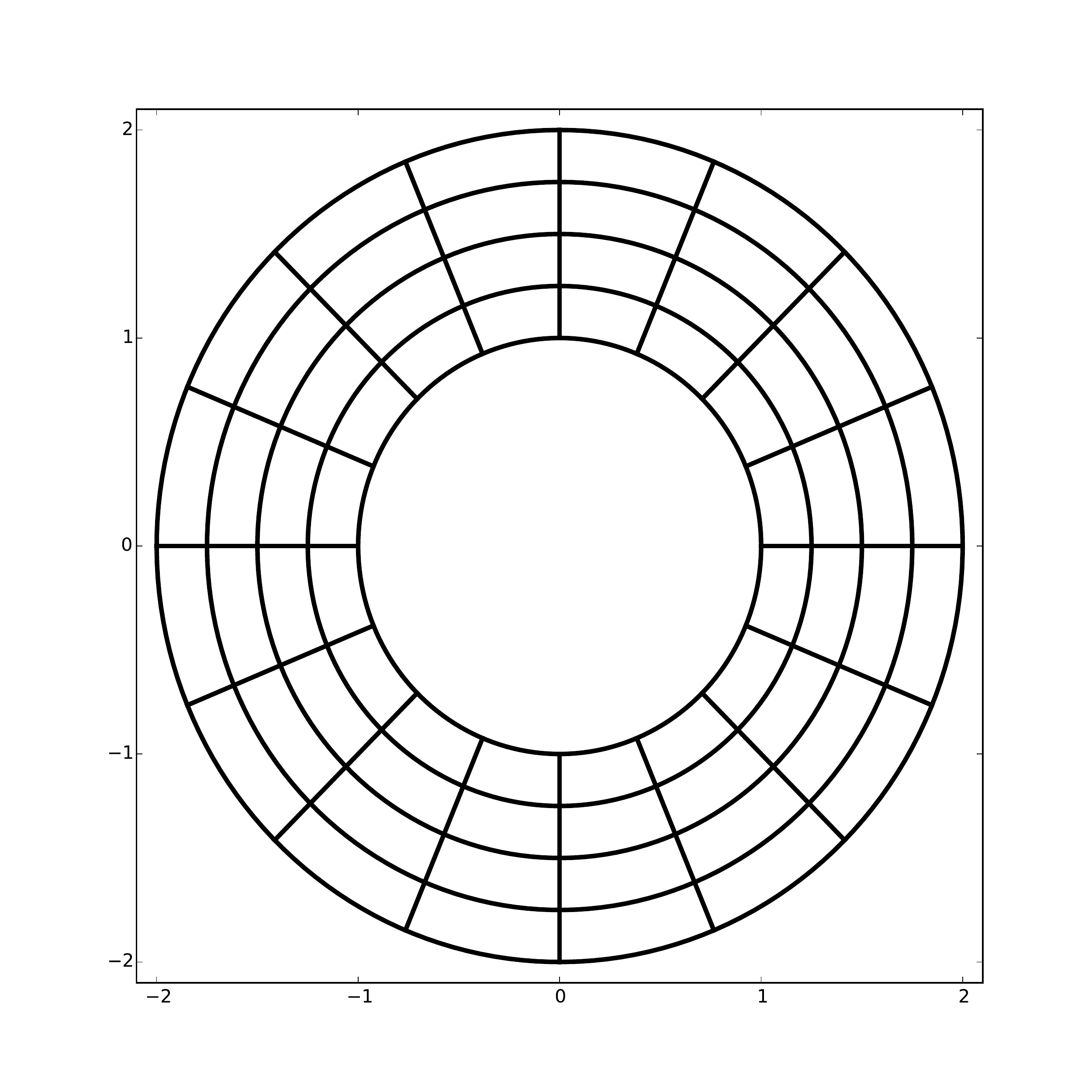}
}
\hspace{1cm}
\subfloat{
\includegraphics[trim=3.4cm 2.9cm 2.9cm 2.9cm, clip=true,height=5cm]{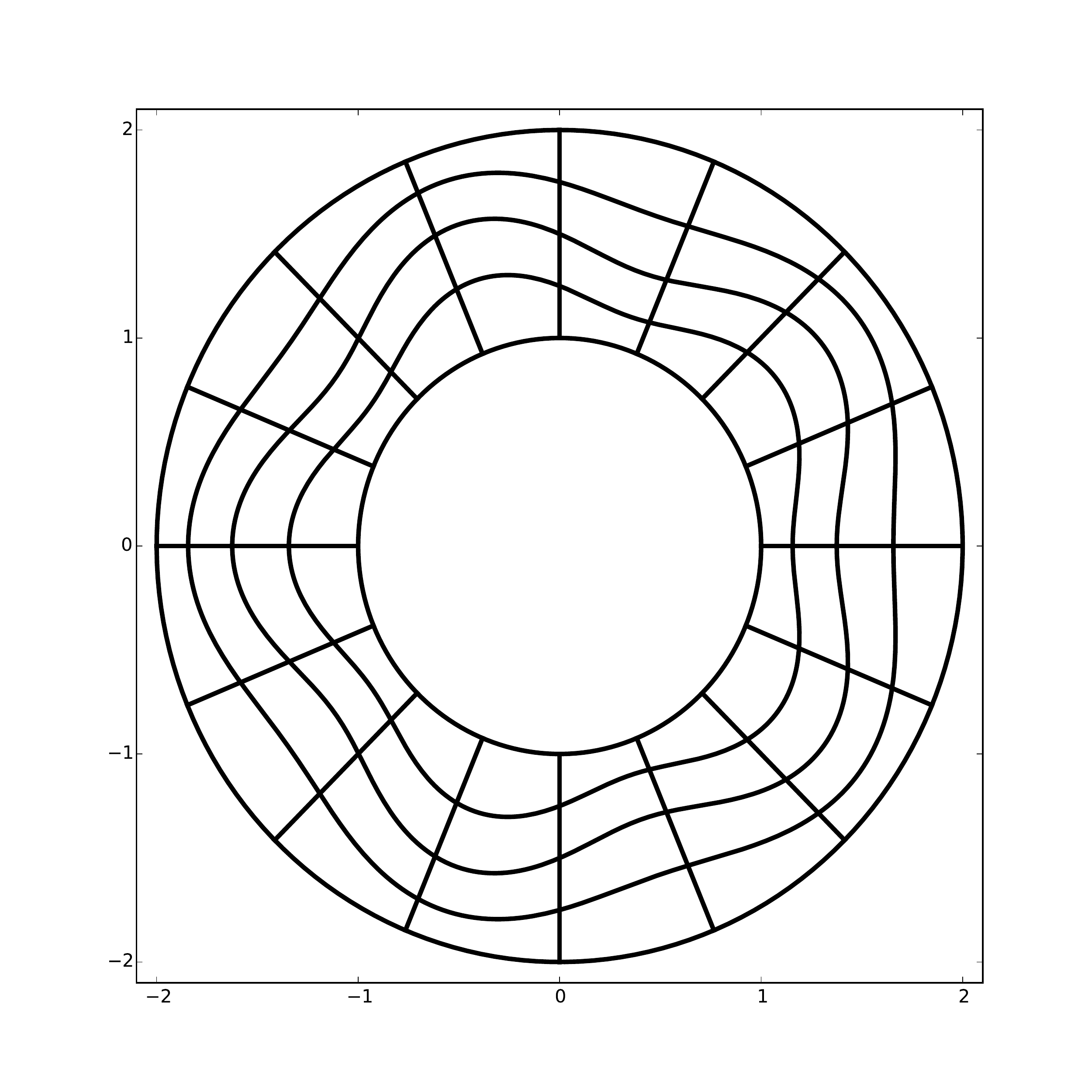}
}
\caption{Meshes of 4$\times$16 elements used to discretize the domain using the analytical mapping. Mesh with $d\!=\!0$ on the left, and $a\!=\!5$, $d\!=\!0.5$ on the left.}
\label{fig:meshannulus}
\end{figure}

\begin{figure}[H]
\captionsetup{aboveskip=0.1cm}
\captionsetup{belowskip=-0.3cm}
\centering
\subfloat{
\includegraphics[trim={0.2cm 0 1.8cm 0.75cm},clip=true,width=7cm]{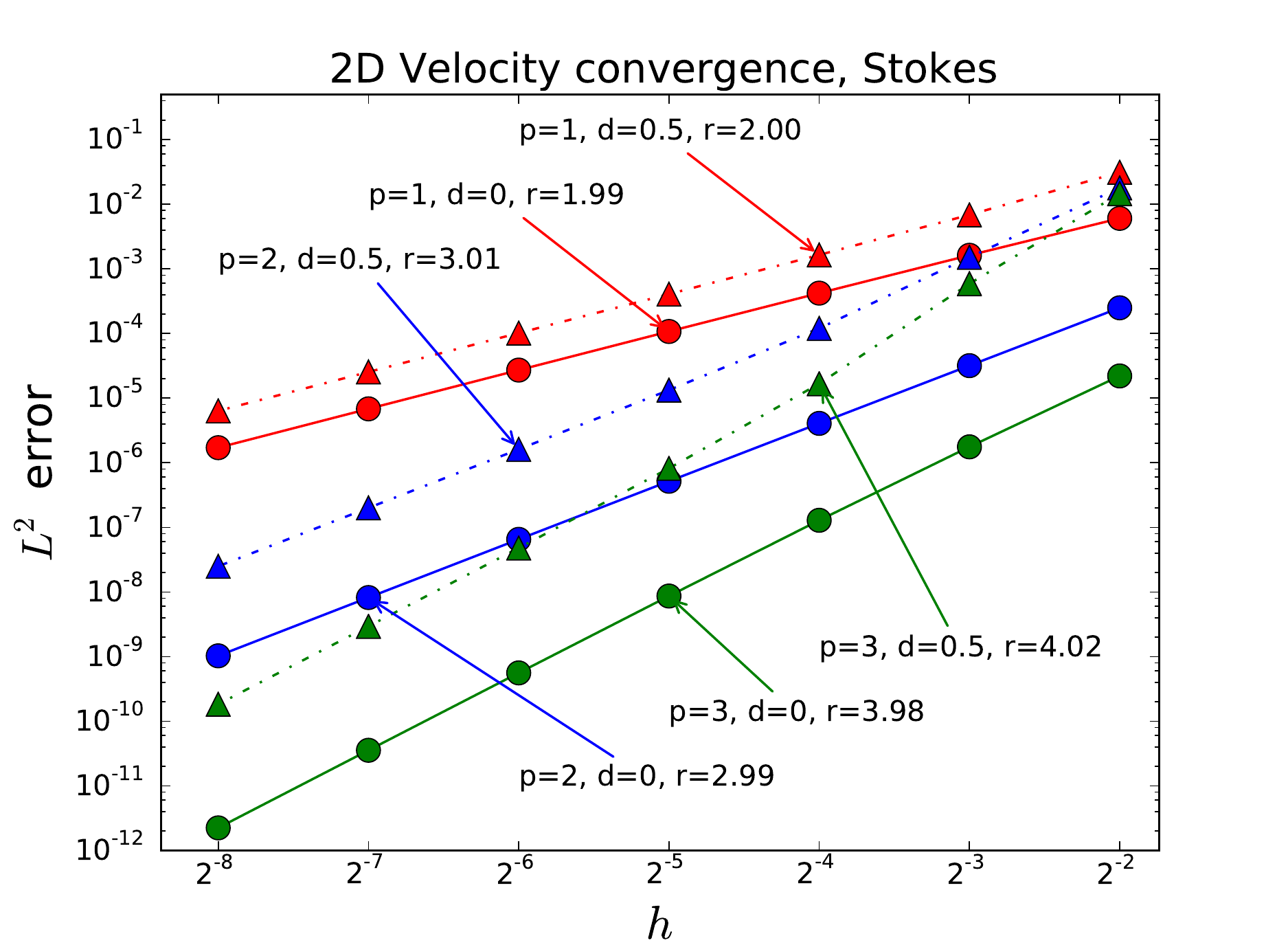}
}\\
\caption{Convergence test results for Stokes in the Couette flow problem using an analytical mapping.}
\label{fig:anannulusStokes}
\end{figure}

\begin{figure}[H]
\captionsetup{aboveskip=0.1cm}
\captionsetup{belowskip=-0.3cm}
\centering
\subfloat{
\includegraphics[trim={0.2cm 0 1.8cm 0.75cm},clip=true,width=7cm]{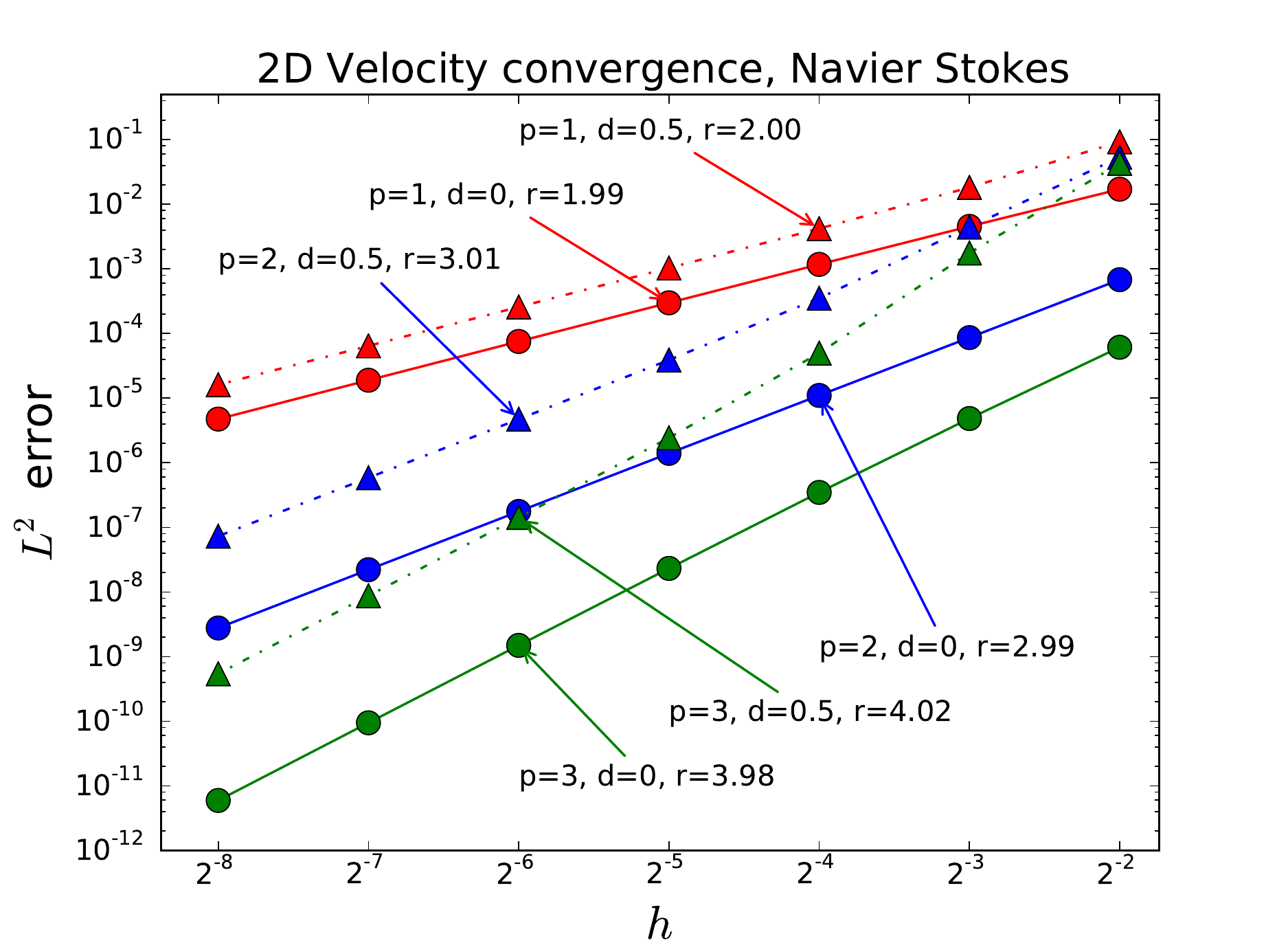}
}
\hspace{0.15cm}
\subfloat{
\includegraphics[trim={0.2cm 0 1.8cm 0.75cm},clip=true,width=7cm]{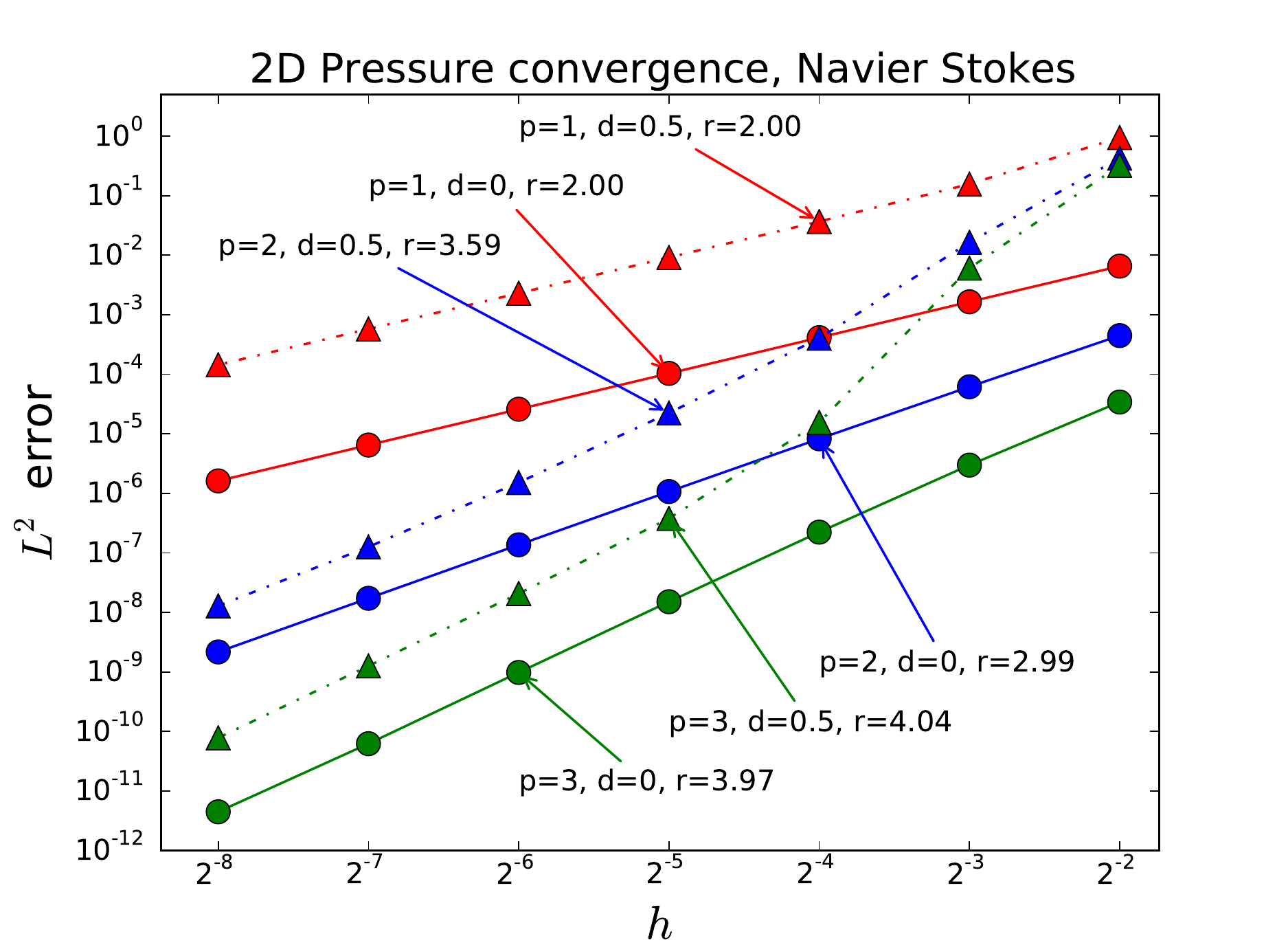}
}\\
\caption{Convergence test results for Navier-Stokes $Re\!=\!1$ in the Couette flow problem using an analytical mapping.}
\label{fig:anannulusNSRe1}
\end{figure}

\begin{figure}[H]
\captionsetup{aboveskip=0.1cm}
\captionsetup{belowskip=-0.3cm}
\centering
\subfloat{
\includegraphics[trim={0.2cm 0 1.8cm 0.75cm},clip=true,width=7cm]{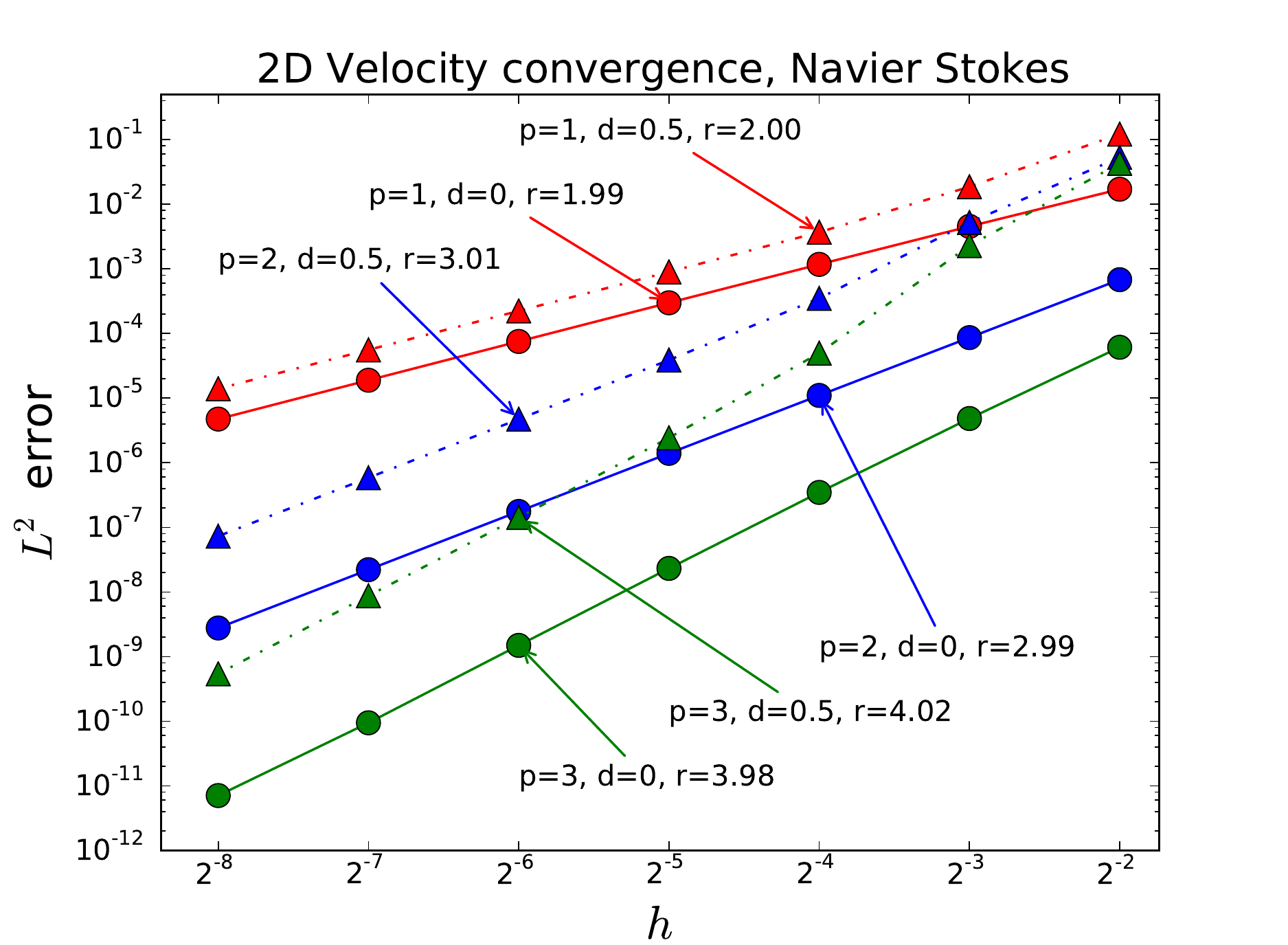}
}
\hspace{0.15cm}
\subfloat{
\includegraphics[trim={0.2cm 0 1.8cm 0.75cm},clip=true,width=7cm]{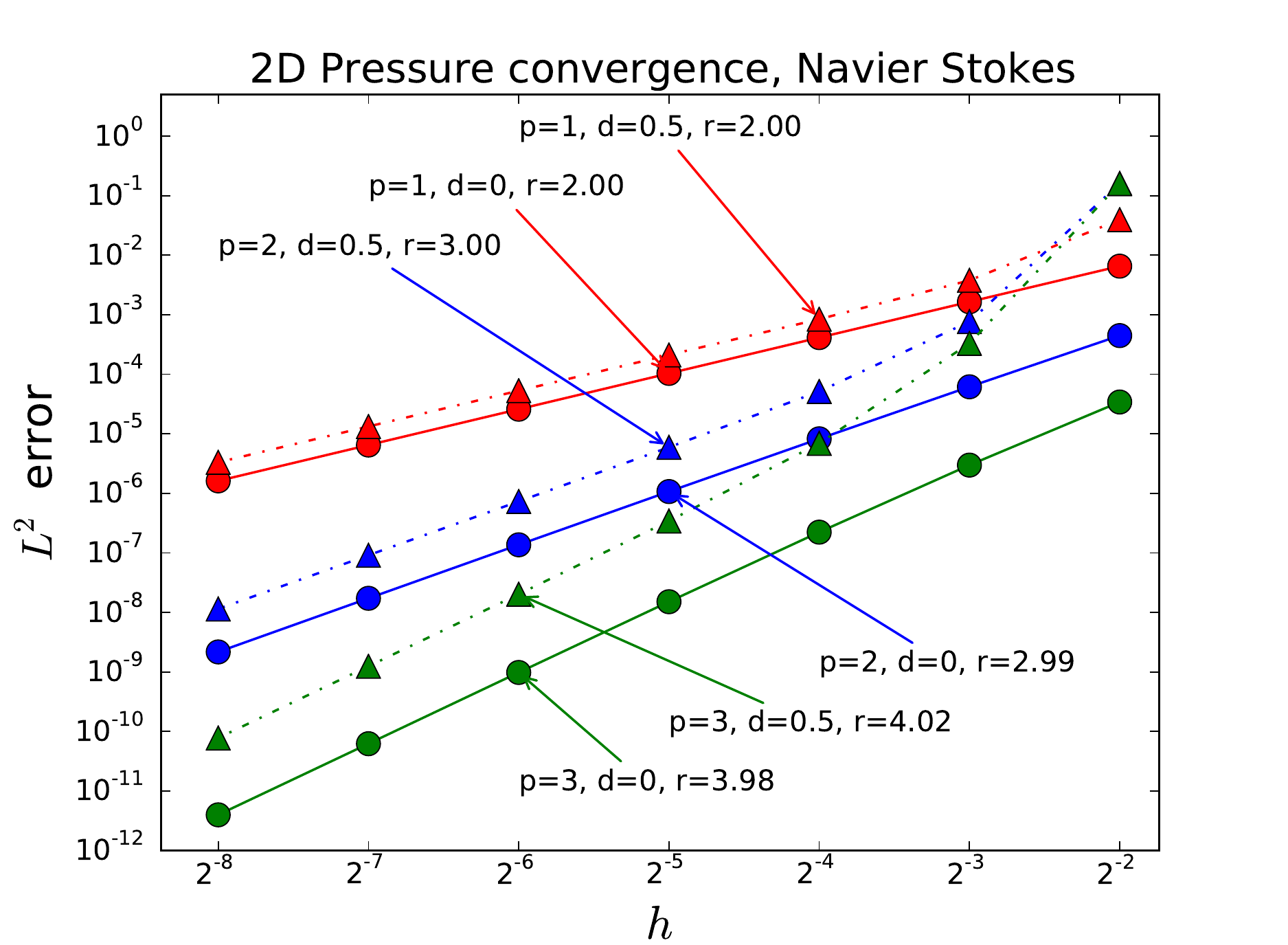}
}\\
\caption{Convergence test results for Navier-Stokes $Re\!=\!100$ in the Couette flow problem using an analytical mapping.}
\label{fig:anannulusNSRe100}
\end{figure}

The domain is defined by the inner radius $r_{in}\!=\!1$ and the outer radius $r_{out}\!=\!2$. The simulations use the analytical mapping described in equation \eqref{eq:map}, where $d$ indicates the distortion from the polar mapping, generating a mesh as shown in Figure \ref{fig:meshannulus}. The results found using the analytical mapping are shown in Figure \ref{fig:anannulusStokes} for the Stokes problem, and in Figures \ref{fig:anannulusNSRe1} and \ref{fig:anannulusNSRe100} for the Navier-Stokes one.

\begin{equation}\label{eq:map}
\mathbf{F}(\xi_1,\xi_2)=\begin{bmatrix}
(d\>\cos(2a\pi\xi_2)(\xi_1^2\!-\!\xi_1)\!+\!\xi_1\!+\!1)\cos(2\pi\xi_2)\\[0.15cm] 
(d\>\cos(2a\pi\xi_2)(\xi_1^2\!-\!\xi_1)\!+\!\xi_1\!+\!1)\sin(2\pi\xi_2)
\end{bmatrix}, \forall (\xi_1,\xi_2) \in \widehat{\Omega}, a\in\mathbb{Z}, d\in[-1,1]
\end{equation}

Figures \ref{fig:anannulusStokes} to \ref{fig:anannulusNSRe100} show that in this case, both velocity and pressure converge at a rate $r\!=\!p+1$, for uniform and distorted meshes. Convergence for pressure in the case of the Stokes problem is omitted since both the uniform and the distorted meshes solve the homogeneous zero condition exactly.

\subsection{Solution in a unitary cube}
\begin{figure}
\centering
\subfloat[Velocity magnitude.]{
\includegraphics[height=5cm]{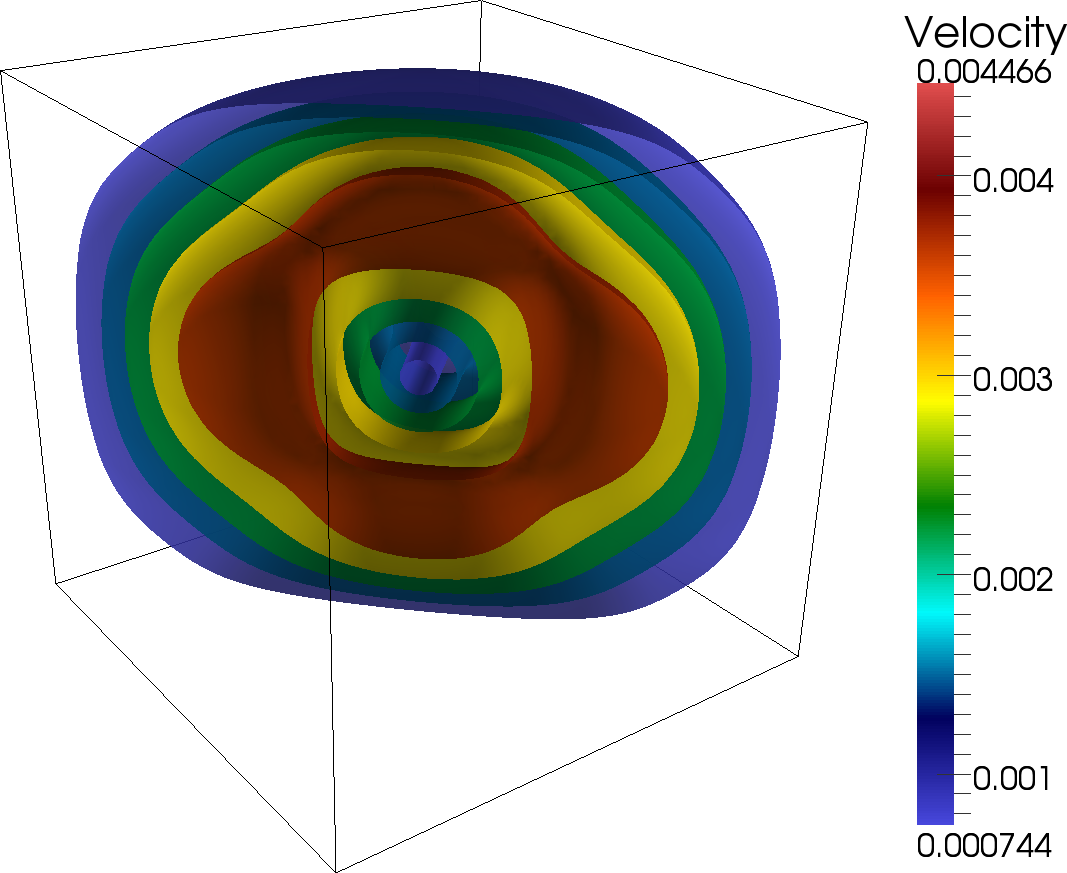}
}
\hspace{0.5cm}
\subfloat[Pressure.]{
\includegraphics[height=5cm]{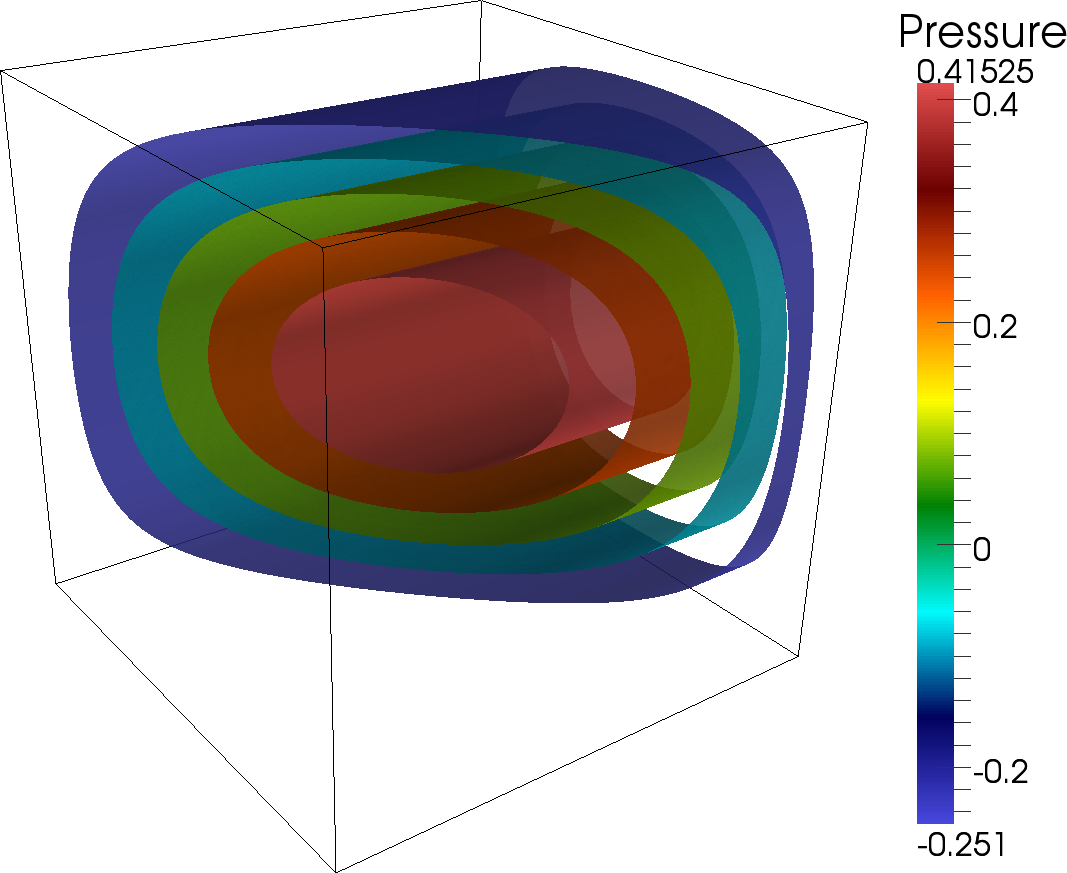}
}\\
\caption{Isocontours of velocity magnitude and pressure in a unitary cube.}
\label{fig:cubeanalytic}
\end{figure}

Here we test our 3D implementation of the Darcy, Stokes, Brinkman, and Navier-Stokes problems against a three-dimensional manufactured solution. The forcing is applied with homogeneous boundary conditions, and the analytic solution as shown in Figure \ref{fig:cubeanalytic} is given by

\begin{align*}
\overline{\bu}=&\nabla \times \overline{\phi}\\
\overline{p}=&\sin(\pi x)\sin(\pi y)-\frac{4}{\pi^2}
\end{align*}
where
\begin{align*}
\overline{\phi}=\begin{bmatrix}
x(x-1)y^2(y-1)^2z^2(z-1)^2\\
0\\
x^2(x-1)^2y^2(y-1)^2z(z-1).
\end{bmatrix}
\end{align*}

\begin{figure}
\captionsetup{aboveskip=0.1cm}
\captionsetup{belowskip=-0.3cm}
\centering
\subfloat[Undistorted mesh.]{
\includegraphics[height=5cm]{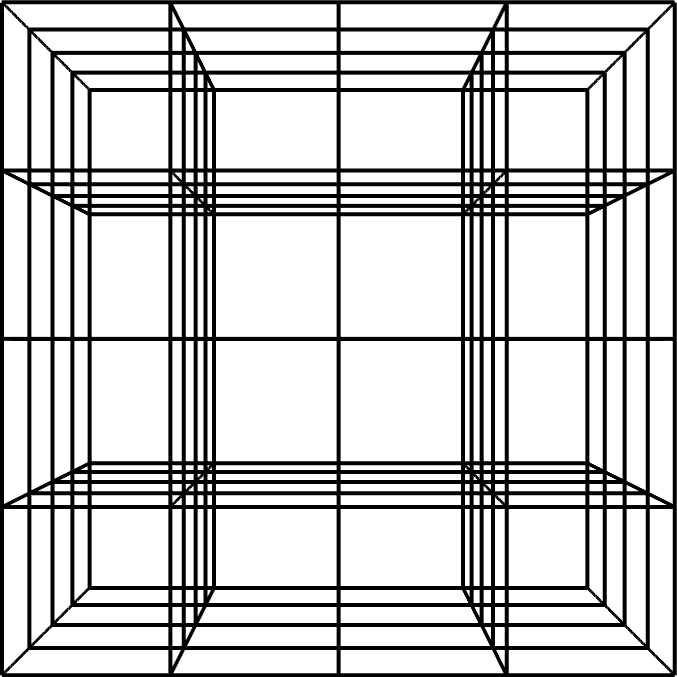}
}
\hspace{2cm}
\subfloat[Distorted mesh.]{
\includegraphics[height=5cm]{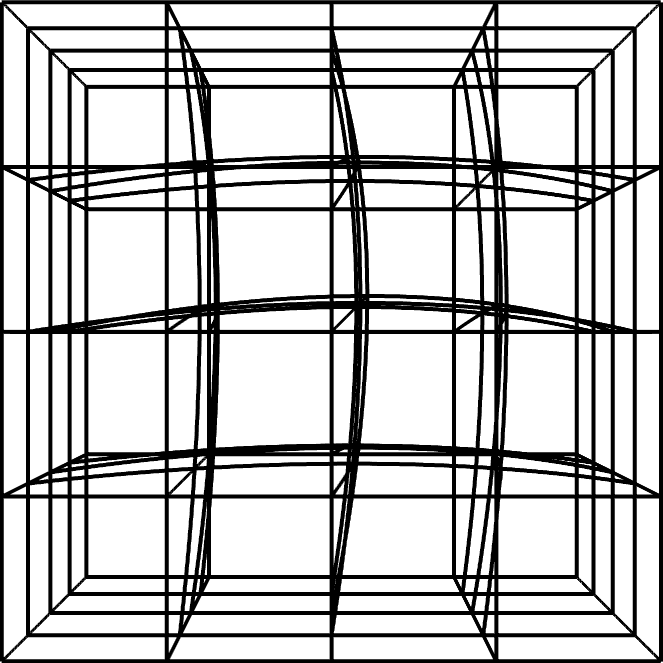}
}
\caption{Meshes of 4$\times$4$\times$4 elements used to discretize the physical domain.}
\label{fig:meshcube}
\end{figure}

To test our implementation in the parametric and physical domains a set of nested uniform and distorted meshes, from $4\times4\times4$ to $64\times64\times64$ elements, were considered for this case. The distorted mesh was built by moving the central control point of a single second order element, a distance $d$ in the positive direction of every axis and then performing an $h$-refinement. An example of the meshes used is shown in Figure \ref{fig:meshcube}.

\begin{figure}
\captionsetup{aboveskip=0.1cm}
\captionsetup{belowskip=-0.3cm}
\centering
\subfloat{
\includegraphics[trim={0.2cm 0 1.8cm 0.75cm},clip=true,width=7cm]{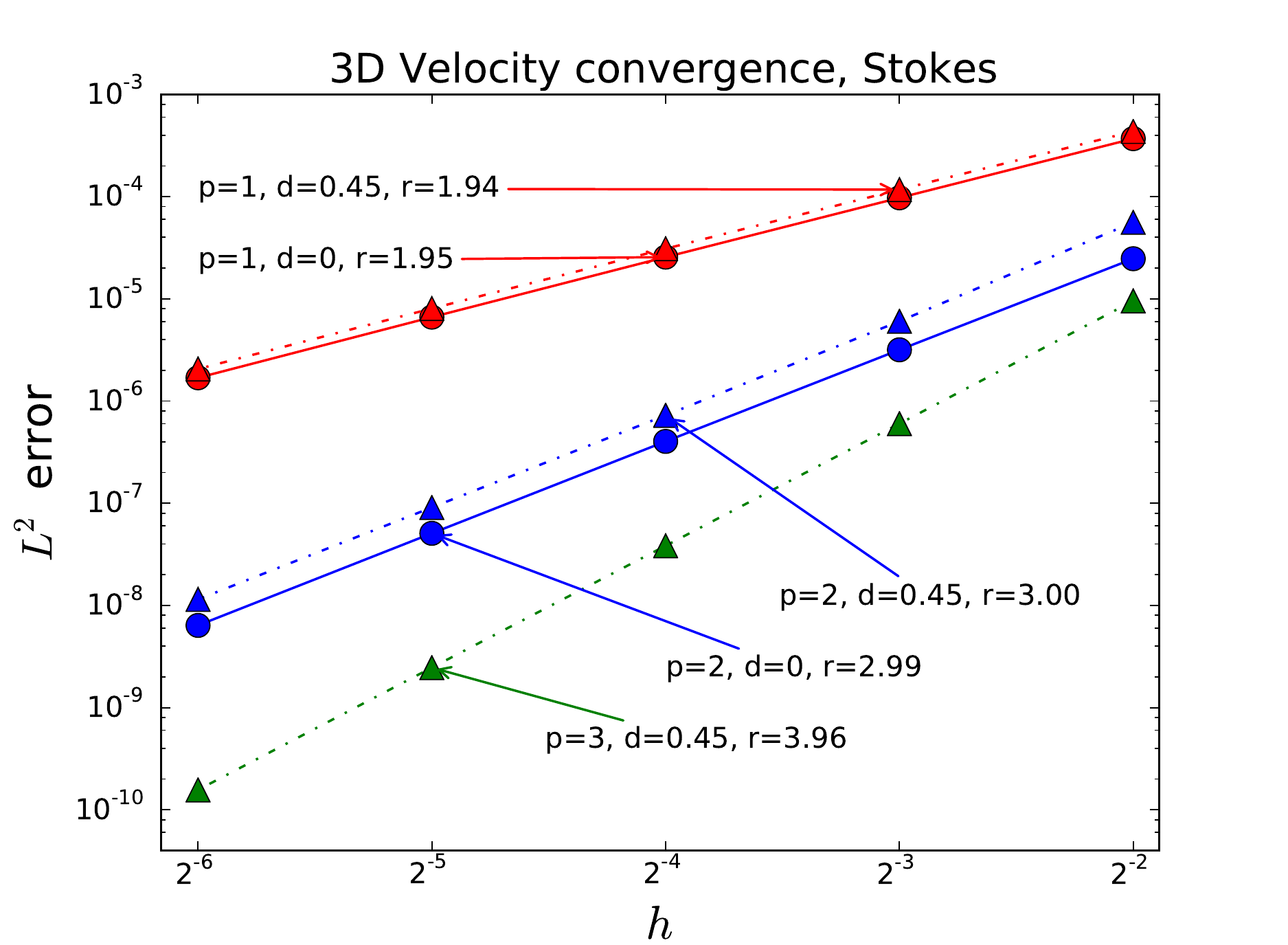}
}
\hspace{0.15cm}
\subfloat{
\includegraphics[trim={0.2cm 0 1.8cm 0.75cm},clip=true,width=7cm]{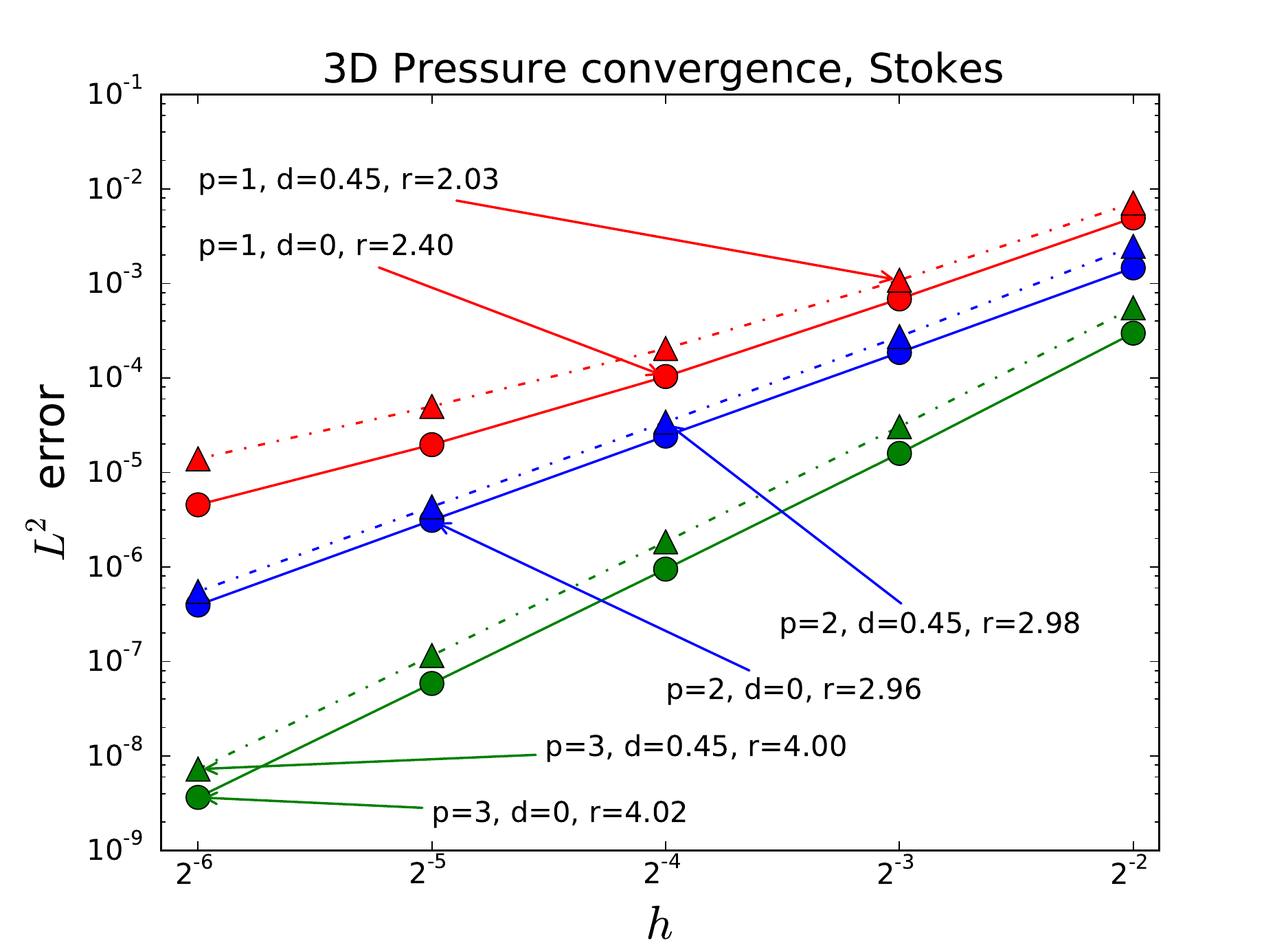}
}\\
\caption{Convergence test results for Stokes in a unitary cube.}
\label{fig:cubeconvStokes}
\end{figure}

\begin{figure}
\captionsetup{aboveskip=0.1cm}
\captionsetup{belowskip=-0.3cm}
\centering
\subfloat{
\includegraphics[trim={0.2cm 0 1.8cm 0.75cm},clip=true,width=7cm]{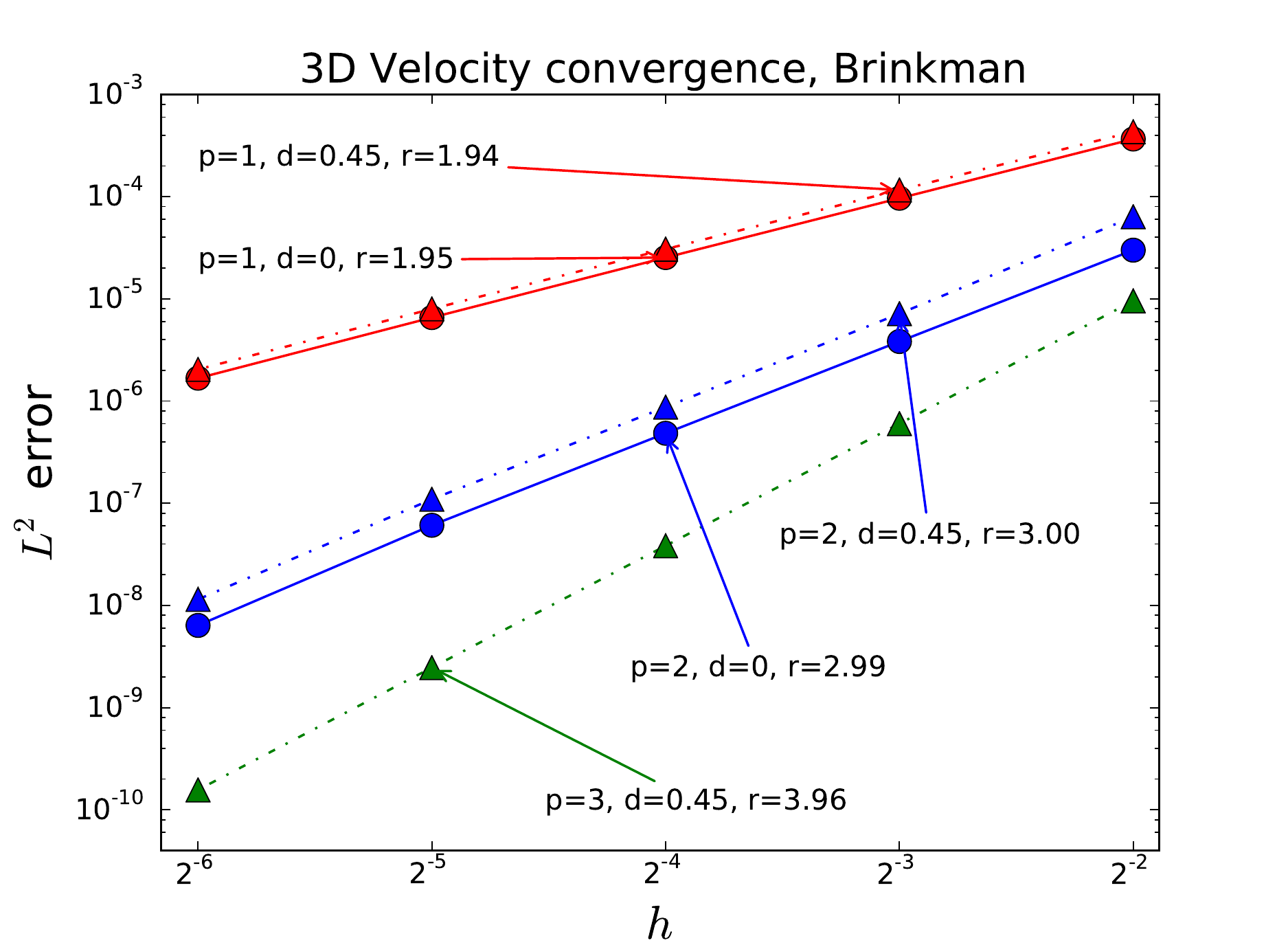}
}
\hspace{0.15cm}
\subfloat{
\includegraphics[trim={0.2cm 0 1.8cm 0.75cm},clip=true,width=7cm]{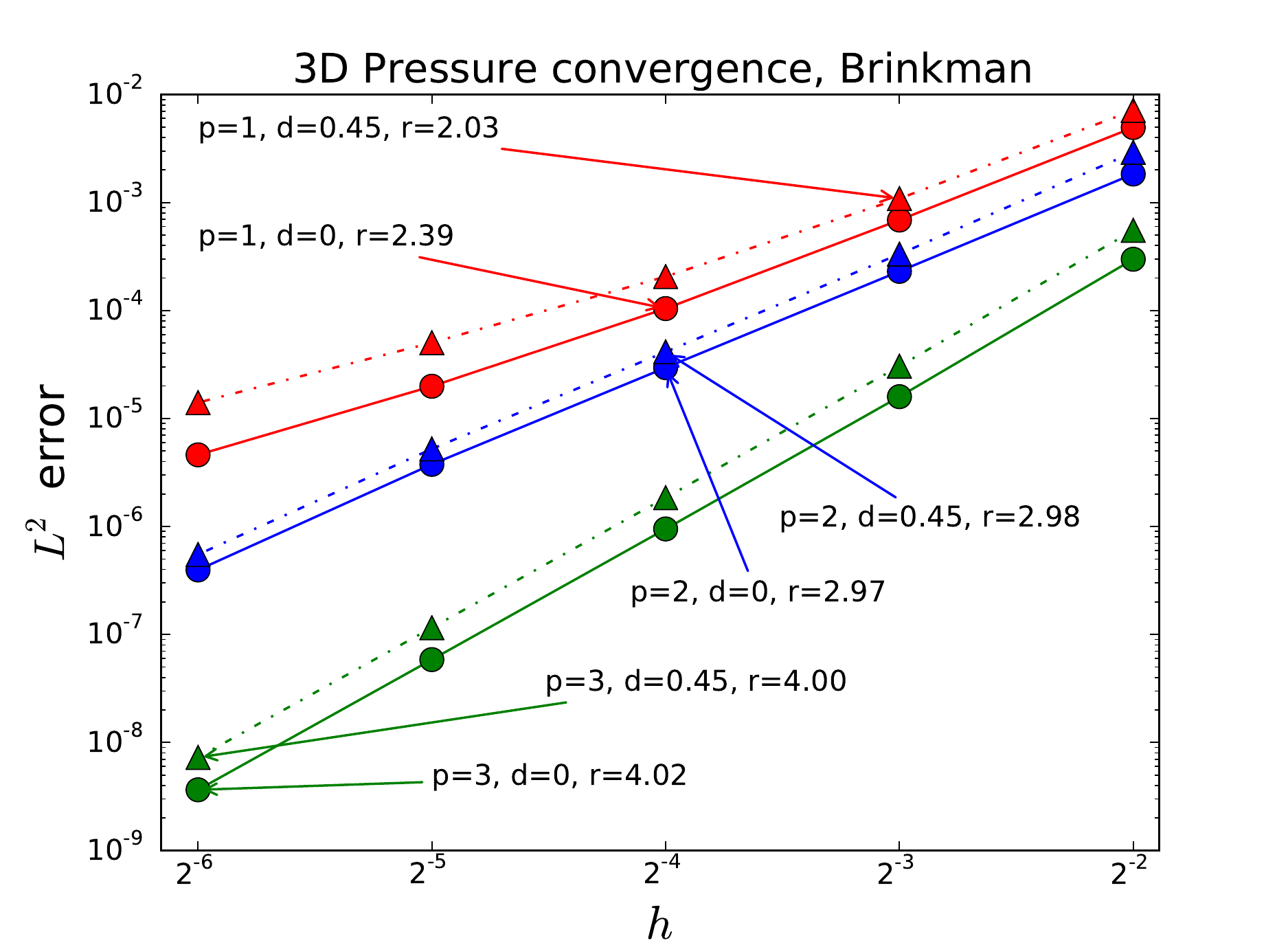}
}\\
\caption{Convergence test results for Brinkman in a unitary cube.}
\label{fig:cubeconvBrinkman}
\end{figure}

\begin{figure}
\captionsetup{aboveskip=0.1cm}
\captionsetup{belowskip=-0.3cm}
\centering
\subfloat{
\includegraphics[trim={0.2cm 0 1.8cm 0.75cm},clip=true,width=7cm]{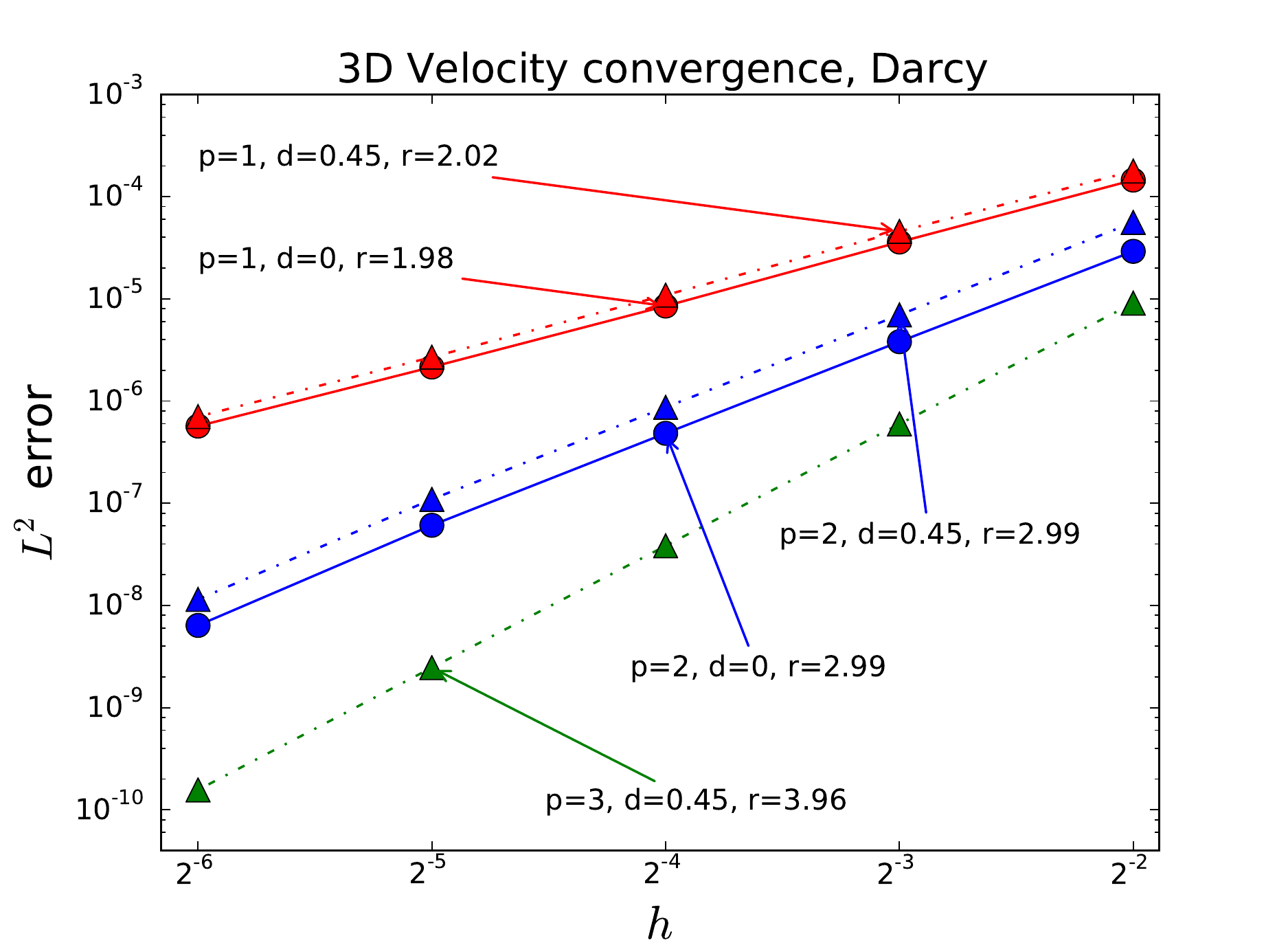}
}
\hspace{0.15cm}
\subfloat{
\includegraphics[trim={0.2cm 0 1.8cm 0.75cm},clip=true,width=7cm]{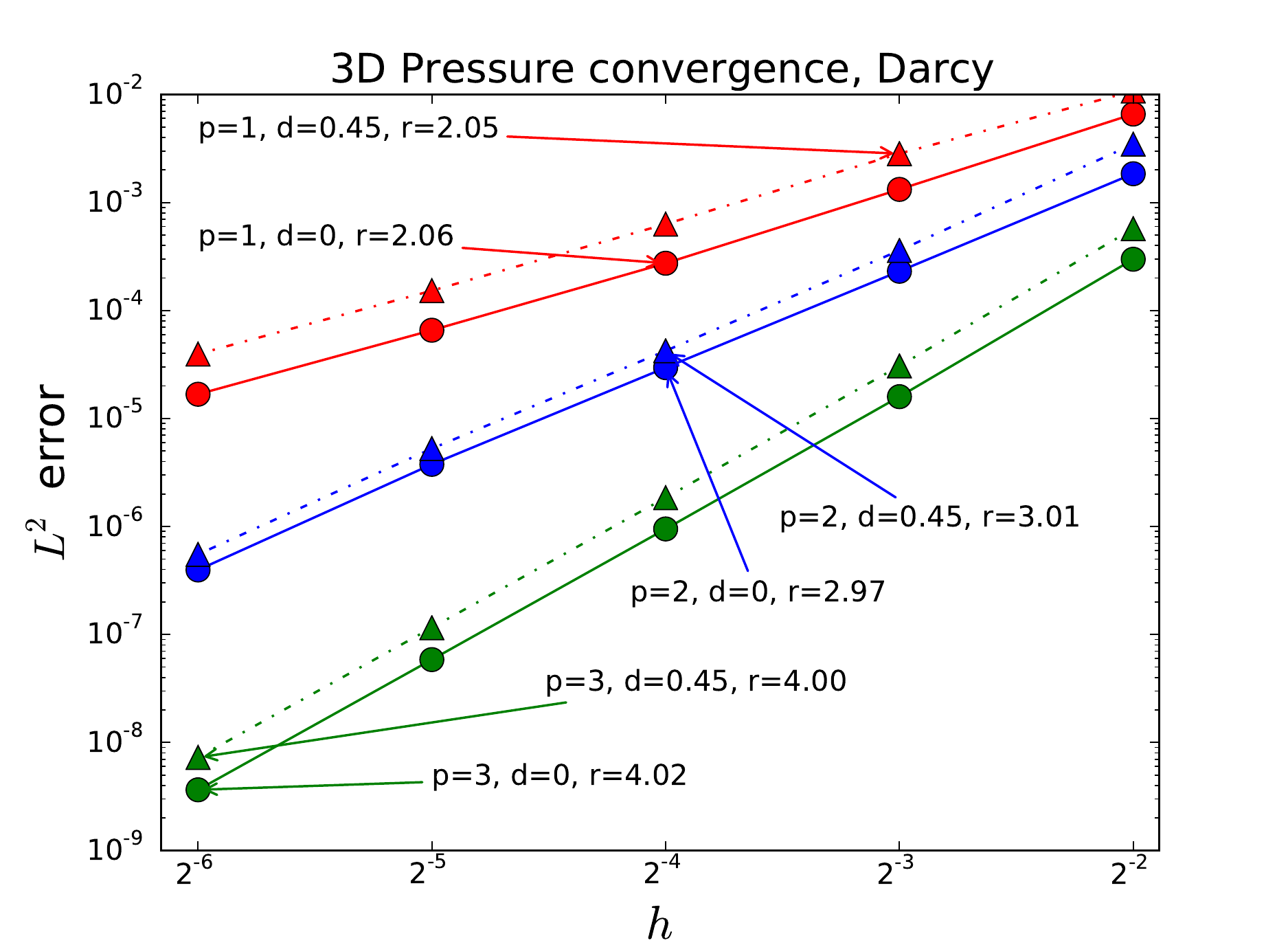}
}\\
\caption{Convergence test results for Darcy in a unitary cube.}
\label{fig:cubeconvDarcy}
\end{figure}

\begin{figure}
\captionsetup{aboveskip=0.1cm}
\captionsetup{belowskip=-0.3cm}
\centering
\subfloat{
\includegraphics[trim={0.2cm 0 1.8cm 0.75cm},clip=true,width=7cm]{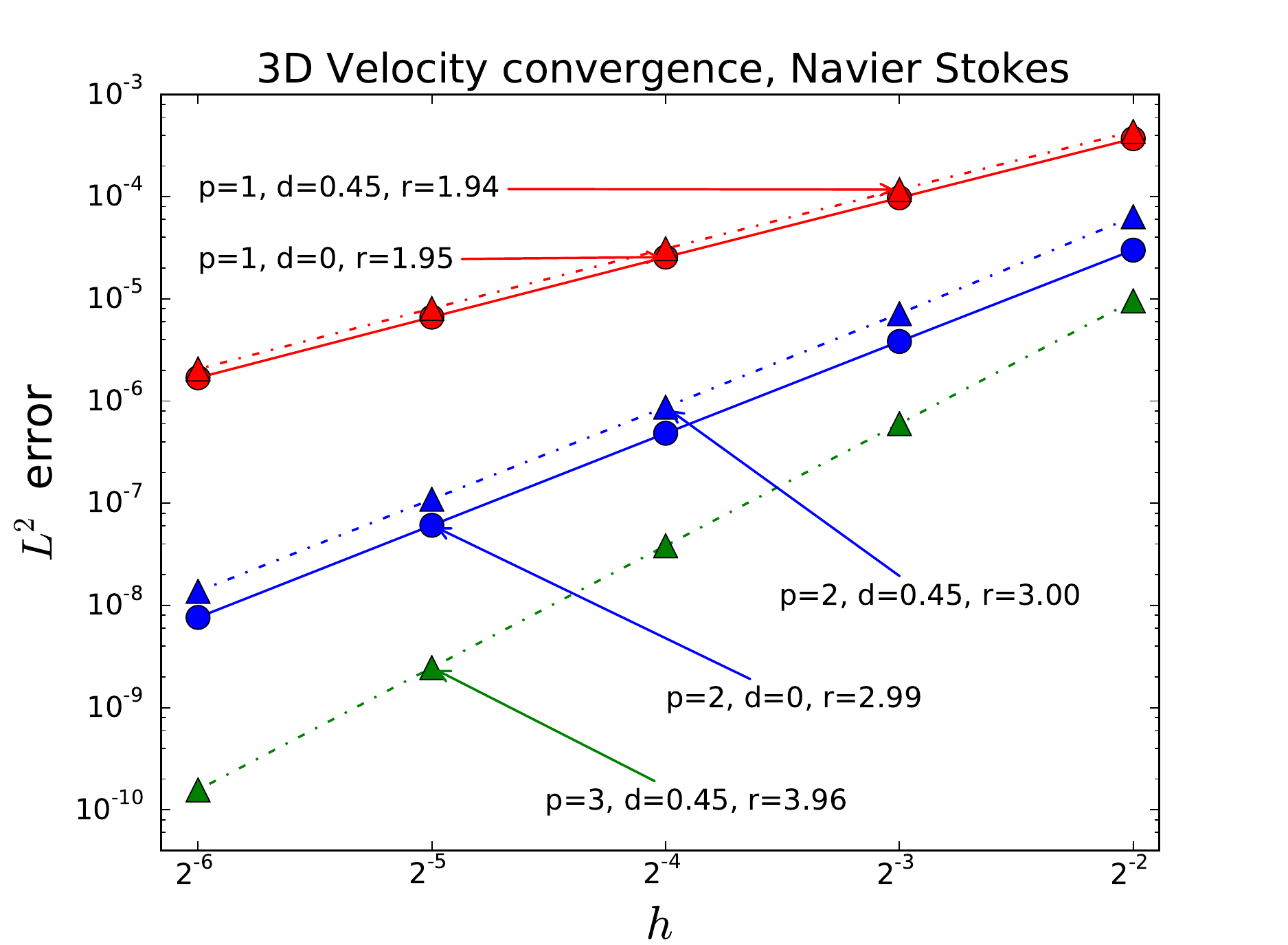}
}
\hspace{0.15cm}
\subfloat{
\includegraphics[trim={0.2cm 0 1.8cm 0.75cm},clip=true,width=7cm]{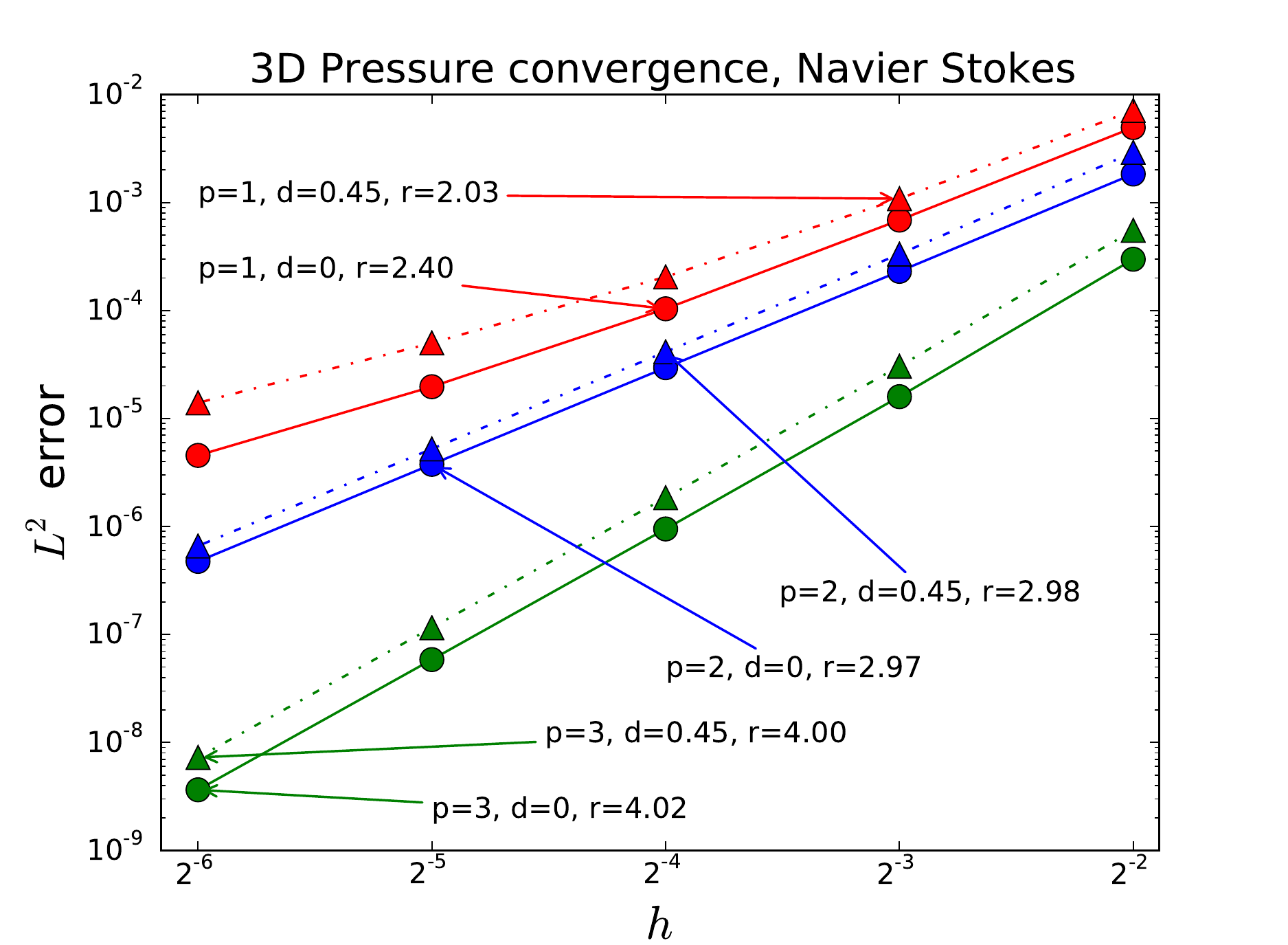}
}\\
\caption{Convergence test results for Navier-Stokes $Re\!=\!1$ in a unitary cube.}
\label{fig:cubeconvNStokes}
\end{figure}

Results shown in Figures \ref{fig:cubeconvStokes} to \ref{fig:cubeconvNStokes}, present convergence rates for velocity and pressure equal to $p\!+\!1$ for both uniform and distorted meshes, the exception being pressure with a uniform mesh and a discretization of $p\!=\!1$ when the convergence rate falls between $p\!+\!1$ and $p\!+\!2$. Uniform meshes using a discretization with $p\!=\!3$ solve exactly the fourth order polynomial given by the analytic solution of the velocity, while the same order discretization using the distorted mesh has a convergence order of $p\!+\!1$, due to the B-spline based mapping.

\subsection{Three-dimensional lid-driven cavity}
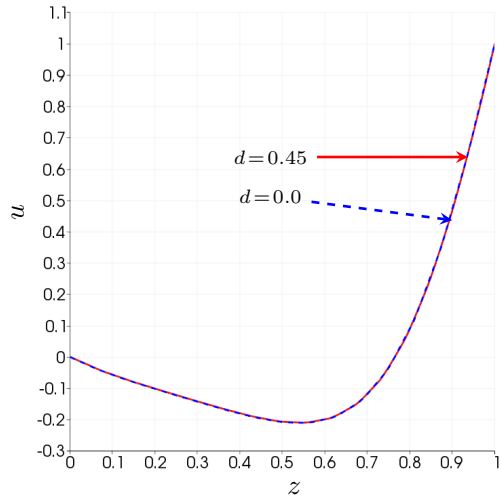
\begin{figure}
\centering
\begin{tikzpicture}[scale=1.05, transform shape]

\node[inner sep=0pt] (russell) at (0  ,0 ) {\includegraphics[height=6cm]{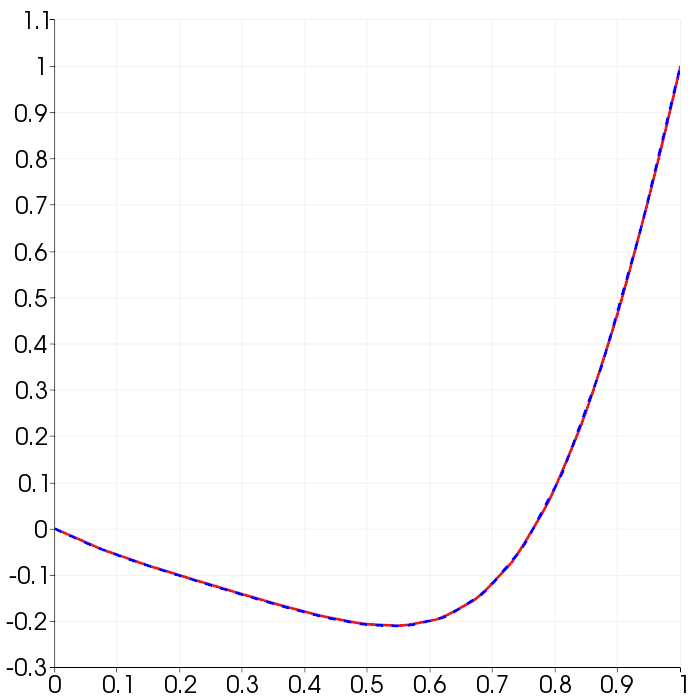}};

\draw (0.3,-3.2) node {$z$};
\draw (-3.2, 0.3) node[rotate=90] {$u$};

\draw (0,1  ) node (1u){\scriptsize $d\!=\!0.45$};
\draw (0,0.5) node (2u){\scriptsize $d\!=\!0.0 $};

\draw[>=stealth,->,       red  ,line width=1pt](1u)--(2.5 , 1  );
\draw[>=stealth,->,dashed,blue ,line width=1pt](2u)--(2.3 , 0.2);

\end{tikzpicture}
\caption{Comparison of horizontal velocity $u$ along the vertical centerline for uniform~($d\!=\!0.0$) and distorted~($d\!=\!0.45$) meshes, solving the Stokes problem in a cube, using $h\!=\!$~\sfrac{1}{16} and $p\!=\!1$, $\varsigma\!=\!0$.}
\label{fig:Stokescube}
\end{figure}

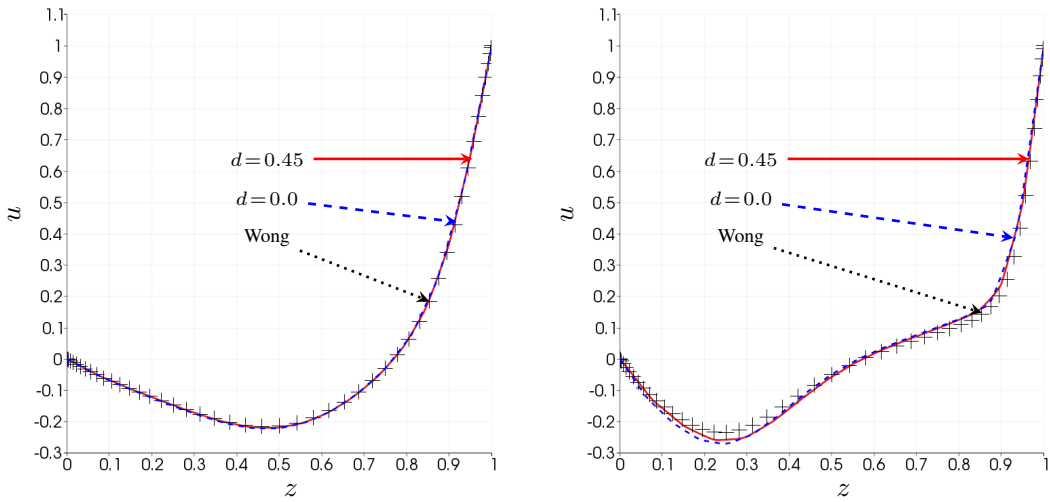
\begin{figure}
\centering
\begin{tikzpicture}[scale=1.05, transform shape]

\node[inner sep=0pt] (russell) at (0  ,0 ) {\includegraphics[height=6cm]{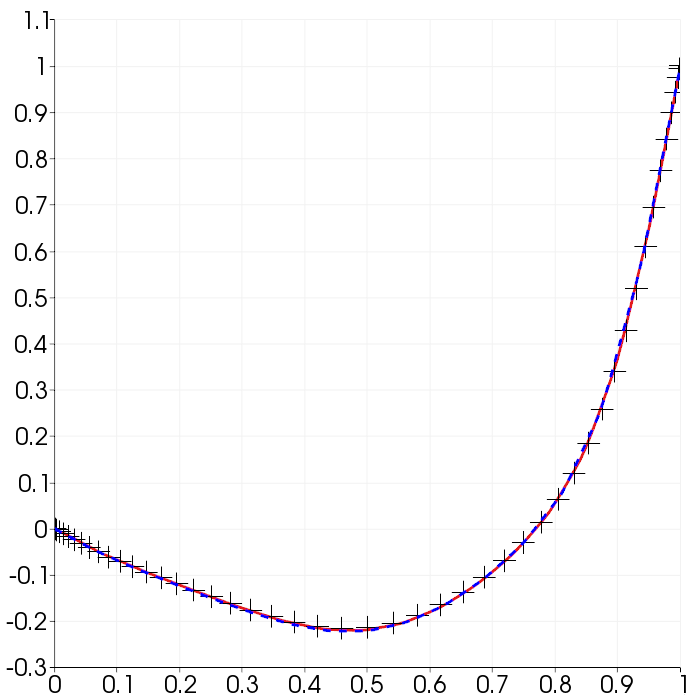}};

\draw (0.3,-3.2) node {$z$};
\draw (-3.2, 0.3) node[rotate=90] {$u$};

\draw (0,1  ) node (1u){\scriptsize $d\!=\!0.45$};
\draw (0,0.5) node (2u){\scriptsize $d\!=\!0.0 $};
\draw (0,0  ) node (3u){\scriptsize Wong};

\draw[>=stealth,->,       red  ,line width=1pt](1u)--(2.6, 1  );
\draw[>=stealth,->,dashed,blue ,line width=1pt](2u)--(2.4 , 0.2);
\draw[>=stealth,->,dotted,black,line width=1pt](3u)--(2.05,-0.8);

\node[inner sep=0pt] (russell) at (7  ,0 ) {\includegraphics[height=6cm]{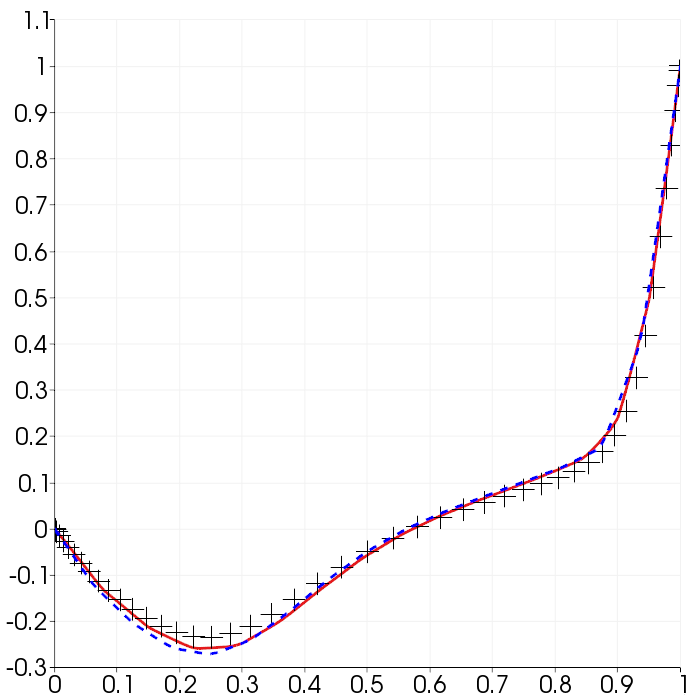}};

\draw (7.3,-3.2) node {$z$};
\draw (3.8, 0.3) node[rotate=90] {$u$};

\draw (6,1  ) node (1u){\scriptsize $d\!=\!0.45$};
\draw (6,0.5) node (2u){\scriptsize $d\!=\!0.0 $};
\draw (6,0  ) node (3u){\scriptsize Wong};

\draw[>=stealth,->,       red  ,line width=1pt](1u)--(9.65, 1  );
\draw[>=stealth,->,dashed,blue ,line width=1pt](2u)--(9.45, 0.0);
\draw[>=stealth,->,dotted,black,line width=1pt](3u)--(9.05,-0.95);

\end{tikzpicture}
\caption{Comparison of horizontal velocity along the vertical centerline for uniform~($d\!=\!0.0$) and distorted~($d\!=\!0.45$) meshes, solving the Navier-Stokes problem in a cube with $Re\!=\!100$ and $Re\!=\!400$, using $h\!=\!$~\sfrac{1}{16} and $p\!=\!1$, $\varsigma\!=\!0$. Our numerical results with the two different meshes compare favorably to Wong's benchmark \cite{Wong}.}
\label{fig:Re100cube}
\end{figure}

We solve the three-dimensional lid-driven cavity test for the Stokes and Navier-Stokes systems, using a set of nested meshes, with and without distortion, from $4\times4\times4$ to $32\times32\times32$ elements. We compare our solution for the Stokes problem with a differential quadrature method using a mesh of $25\times25\times25$ \cite{Lo1}. We solve the Navier-Stokes problem for two different Reynolds numbers $Re\!=\!100$ and $Re\!=\!400$, and compared with the results presented by Lo in \cite{Lo2} using a finite difference method, solving the case for $Re\!=\!100$ and $Re\!=\!400$ with a $51\times51\times51$ and $101\times101\times101$ mesh, respectively. We compare our simulation results with those presented by Wong in \cite{Wong}, using the finite element method to solve the velocity-vorticity formulation with $48\times48\times48$ elements. Table \ref{tab:maxvelcube} compares the values and positions of the minimum horizontal velocity $u$ along the vertical centerline $(x\!=\!0.5,\,y\!=\!0.5)$. Figures \ref{fig:Stokescube} and \ref{fig:Re100cube} show the results of the horizontal velocity along the vertical centerline, found with a discretization using $h\!=\!$ \sfrac{1}{16}, $p\!=\!1$, $\varsigma\!=\!0$, comparing the results of a uniform mesh against the ones from a distorted mesh ($d\!=\!0.45$).

\noindent
\begin{table}[H]
\captionsetup{belowskip=-0.15cm}
\captionof{table}{Convergence of the velocity extrema for the Stokes and Navier-Stokes problem in a cube using uniform meshes ($d\!=\!0$).} 
\label{tab:maxvelcube} 
\centering
\small
\begin{tabular}{C{2cm}C{1cm}C{1.3cm}C{1.3cm}C{1.3cm}C{1.3cm}C{1.3cm}C{1.3cm}C{0.1cm}}
\cline{1-8}
\multicolumn{1}{ C{2cm}   }{\multirow{2}{*}{Discretization} } &
\multicolumn{1}{ C{1cm}   }{\multirow{2}{*}{$h$} } &
\multicolumn{2}{ C{2.6cm} }{\vspace{0.2cm}Stokes      }&
\multicolumn{2}{ C{2.6cm} }{\vspace{0.2cm}$Re\!=\!100$}&
\multicolumn{2}{ C{2.6cm} }{\vspace{0.2cm}$Re\!=\!400$}\\[0.15cm] \cline{3-8}
 								          			   &&$u_{min}$		&$z_{min}$	&$u_{min}$		&$z_{min}$	&$u_{min}$		&$z_{min}$	&\\[0.15cm]\cline{1-8}
\multicolumn{1}{ c }{\multirow{4}{*}{$p\!=\!1$, $\varsigma\!=\!0$}}
							& \sfrac{1}{4}   			&-0.28512		&0.50000	&-0.34646		&0.50000	&-0.32065		&0.50000	&\\[0.15cm]
							& \sfrac{1}{8}   			&-0.21946		&0.50000	&-0.24737		&0.50000	&-0.34257		&0.25000	&\\[0.10cm]
							& \sfrac{1}{16}  			&-0.21070		&0.56249	&-0.22176		&0.50000	&-0.27280		&0.25000	&\\[0.10cm]
							& \sfrac{1}{32}  			&-0.20868		&0.53125	&-0.21764		&0.46875	&-0.24539		&0.25000	&\\[0.10cm]\cline{1-8}
\multicolumn{1}{ c }{\multirow{4}{*}{$p\!=\!2$, $\varsigma\!=\!1$}}
							& \sfrac{1}{4}   			&-0.21778		&0.53553	&-0.22364		&0.43358	&-0.30577		&0.30559	&\\[0.15cm]
							& \sfrac{1}{8}   			&-0.20796		&0.53237	&-0.21774		&0.46565	&-0.25739		&0.24299	&\\[0.10cm]
							& \sfrac{1}{16}  			&-0.20771		&0.53468	&-0.21575		&0.46830	&-0.24086		&0.24008	&\\[0.10cm]
							& \sfrac{1}{32}  			&-0.20776		&0.53605	&-0.21560		&0.46923	&-0.23702		&0.23885	&\\[0.10cm]\cline{1-8}
\multicolumn{2}{ c }{Differential quadrature \cite{Lo1}}&-0.231			&-		 	&-0.215 		&-			&-0.236 		&-			&\\[0.10cm] 
\multicolumn{2}{ c }{Finite differences \cite{Lo2}}		&- 				&-		 	&-0.2163		&0.46		&-0.2334		&0.26		&\\[0.10cm] 
\multicolumn{2}{ c }{Finite elements \cite{Wong}}       &-				&-		 	&-0.2154		&0.4592 	&-0.2349 		&0.2509  	&\\[0.10cm]\cline{1-8}
\end{tabular}
\end{table}

Table \ref{tab:maxvelcube} compares the minimum horizontal velocity and its position over the vertical centerline, against solutions using differential quadrature \cite{Lo1}, finite differences \cite{Lo2} and finite element method \cite{Wong}. For these Reynolds numbers the simulation results are reasonably close to the benchmark values when using meshes with $h\!=\!$ \sfrac{1}{16} and $h\!=\!$ \sfrac{1}{32}, especially for the case of $p\!=\!2$, $\varsigma\!=\!1$. When using discretizations of $p\!=\!2$ the differences between the results of using meshes with $h\!=\!$ \sfrac{1}{16} and $h\!=\!$ \sfrac{1}{32}, are small, suggesting that a mesh $h\!=\!$ \sfrac{1}{32} is enough to represent the flow inside the domain.

Figures \ref{fig:Stokescube} and \ref{fig:Re100cube} compare the results found with uniform and distorted meshes, with the ones reported by Wong \cite{Wong}. In the three cases, the results obtained using a distorted mesh are very similar to those of the uniform meshes, showing no major effect of the mesh distortion over the velocity convergence. Results obtained for the Navier-Stokes equations compare well to those of Wong, being indistinguishable from the benchmark in the case of $Re\!=\!100$.

\section{Acknoledgments} \label{sect:Acknoledgments}

This publication was made possible in part by a National Priorities Research Program grant \mbox{7-1482-1-278} from the Qatar National Research Fund (a member of The Qatar Foundation), by the European Union's Horizon 2020 Research and Innovation Program of the Marie Sk{\l}odowska-Curie grant agreement \mbox{No. 644602}and the Center for Numerical Porous Media at King Abdullah University of Science and Technology (KAUST). L. Dalcin was partially supported by Agencia Nacional de Promoci\'on Cient\'{i}fica y Tecnol\'ogica grants \mbox{PICT~2014--2660} and \mbox{PICT-E~2014--0191}. The J. Tinsley Oden Faculty Fellowship Research Program at the Institute for Computational Engineering and Sciences (ICES) of the University of Texas at Austin has partially supported the visits of VMC to ICES.

\section{Conclusions} \label{sect:Conclusions}

We introduce a framework, called PetIGA-MF, for multi-field high-performance isogeometric analysis. PetIGA-MF provides structure-preserving vector field discretizations to solve multi-physics problems. This framework allows us to use different approximation spaces for each component of the discrete fields. This flexibility simplifies the implementation and discretization of complex multi-field problems while guaranteeing stability. We extend PetIGA and adapt PETSc to manage the parallelism, and offer access to a significant variety of solvers and preconditioners.

We test our simulation framework with numerous benchmarks and evaluate the effect of distorting the mesh on the convergence rates, finding optimal convergence rates for the velocity and pressure in all cases. When using uniform meshes the convergence rates for the pressure are equal to those of the velocity, although its discretization uses spaces with one order lower polynomials. Under mesh distortion, the convergence rates for the pressure decreases, by almost one order, which is closer to the predicted limit of the a priori error estimates. These results lead us to conclude that the error estimates for pressure presented in \cite{EvansHughesStokes} are not conservative, but only distorted meshes and non-trivial geometries may present strictly optimal convergence rates. An in-depth analysis of the effects of mesh distortion on the error estimates is required to understand better the circumstances under which this loss of superconvergence may occur.

Our discrete space choices for velocity and pressure, namely the divergence-conforming spaces, guaranty the conservation of mass at every point in the domain discretely and an accurate solution of the flow. Weak imposition of boundary conditions yields accurate results when focusing on flow near boundaries while avoiding instabilities due to over restricting the velocity space.



\section{References}

\bibliographystyle{paper}	
\bibliography{paper.bib}	

\begin{thebibliography}{10}
\expandafter\ifx\csname url\endcsname\relax
  \def\url#1{\texttt{#1}}\fi
\expandafter\ifx\csname urlprefix\endcsname\relax\def\urlprefix{URL }\fi
\expandafter\ifx\csname href\endcsname\relax
  \def\href#1#2{#2} \def\path#1{#1}\fi

\bibitem{ExtCalc}
D.N. Arnold, R.S. Falk, R.~Winther, Finite element exterior calculus: from
  {H}odge theory to numerical stability, Bulletin American Mathematical Society
  47 (2010) 281--354.
\newblock \href {http://dx.doi.org/10.1090/S0273-0979-10-01278-4}
  {\path{doi:10.1090/S0273-0979-10-01278-4}}.

\bibitem{monk2003finite}
P.~Monk, Finite Element Methods for Maxwell's Equations, Numerical Mathematics
  and Scientific Computation, Oxford University Press, 2003.

\bibitem{Dem2000}
L.~Demkowicz, P.~Monk, L.~Vardapetyan, W.~Rachowicz, {D}e {R}ham diagram for hp
  finite element spaces, Computers \& Mathematics with Applications 39~(7--8)
  (2000) 29 -- 38.
\newblock \href
  {http://dx.doi.org/http://dx.doi.org/10.1016/S0898-1221(00)00062-6}
  {\path{doi:http://dx.doi.org/10.1016/S0898-1221(00)00062-6}}.

\bibitem{HBCBOOK}
J.~A. Cottrell, T.~J.~R Hughes, Y.~Bazilevs, Isogeometric Analysis: Toward
  Integration of {CAD} and {FEA}, Wiley, 2009.
\newblock \href {http://dx.doi.org/10.1002/9780470749081}
  {\path{doi:10.1002/9780470749081}}.

\bibitem{Buffa2011}
A.~Buffa, J.~Rivas, G.~Sangalli, R.~V{\'a}zquez, Isogeometric discrete
  differential forms in three dimensions, SIAM Journal Numerical Analysis
  49~(2) (2011) 818--844.
\newblock \href {http://dx.doi.org/10.1137/100786708}
  {\path{doi:10.1137/100786708}}.

\bibitem{Buffa2010}
A.~Buffa, G.~Sangalli, R.~V{\'a}zquez, Isogeometric analysis in
  electromagnetics: {B-splines} approximation, Computer Methods in Applied
  Mechanics and Engineering 199~(17--20) (2010) 1143--1152.
\newblock \href {http://dx.doi.org/10.1016/j.cma.2009.12.002}
  {\path{doi:10.1016/j.cma.2009.12.002}}.

\bibitem{Buffa2010stokes}
A.~Buffa, C.~de~Falco, G.~Sangalli, Isogeometric analysis: stable elements for
  the {2D} {Stokes} equation, International Journal for Numerical Methods in
  Fluids 65~(11-12) (2011) 1407--1422.
\newblock \href {http://dx.doi.org/10.1002/fld.2337}
  {\path{doi:10.1002/fld.2337}}.

\bibitem{EvansHughesStokes}
J.~A. Evans, T.~J.~R Hughes, Isogeometric divergence-conforming {B-splines} for
  the {Darcy-Stokes-Brinkman} equations, Mathematical Models and Methods in
  Applied Sciences 23~(04) (2013) 671--741.
\newblock \href {http://dx.doi.org/10.1142/S0218202512500583}
  {\path{doi:10.1142/S0218202512500583}}.

\bibitem{EvansHughesNavierStokes}
J.~A. Evans, T.~J.~R Hughes, Isogeometric divergence-conforming {B-splines} for
  the steady {Navier-Stokes} equations, Mathematical Models and Methods in
  Applied Sciences 23~(08) (2013) 1421--1478.
\newblock \href {http://dx.doi.org/10.1142/S0218202513500139}
  {\path{doi:10.1142/S0218202513500139}}.

\bibitem{EvansNSUnsteady}
J.~A. Evans, T.~J.~R Hughes, Isogeometric divergence-conforming {B-splines} for
  the unsteady {Navier-Stokes} equations, Journal of Computational Physics
  241~(0) (2013) 141--167.
\newblock \href {http://dx.doi.org/10.1016/j.jcp.2013.01.006}
  {\path{doi:10.1016/j.jcp.2013.01.006}}.

\bibitem{petiga}
L.~{Dalcin}, N.~{Collier}, P.~{Vignal}, A.~M.~A. {Cortes}, V.~M. {Calo},
  {PetIGA: A Framework for High-Performance Isogeometric Analysis}, ArXiv
  e-prints\href {http://arxiv.org/abs/1305.4452} {\path{arXiv:1305.4452}}.

\bibitem{Rudraraju}
S.~Rudraraju, A.~Van der Ven, K.~Garikipati, Three-dimensional isogeometric
  solutions to general boundary value problems of toupin’s gradient
  elasticity theory at finite strains, Computer Methods in Applied Mechanics
  and Engineering 278 (2014) 705 -- 728.
\newblock \href {http://dx.doi.org/http://dx.doi.org/10.1016/j.cma.2014.06.015}
  {\path{doi:http://dx.doi.org/10.1016/j.cma.2014.06.015}}.

\bibitem{Vignal}
P.~{Vignal}, L.~{Dalcin}, D.~L. {Brown}, N.~{Collier}, V.~M. {Calo}, {An
  energy-stable convex splitting for the phase-field crystal equation},
  Computers \& Structures 158 (2015) 355 -- 368.
\newblock \href
  {http://dx.doi.org/http://dx.doi.org/10.1016/j.compstruc.2015.05.029}
  {\path{doi:http://dx.doi.org/10.1016/j.compstruc.2015.05.029}}.

\bibitem{Wozniak}
M.~Woźniak, K.~Kuźnik, M.~Paszyński, V.~M. Calo, D.~Pardo, Computational
  cost estimates for parallel shared memory isogeometric multi-frontal solvers,
  Computers \& Mathematics with Applications 67~(10) (2014) 1864 -- 1883.
\newblock \href
  {http://dx.doi.org/http://dx.doi.org/10.1016/j.camwa.2014.03.017}
  {\path{doi:http://dx.doi.org/10.1016/j.camwa.2014.03.017}}.

\bibitem{Yokota}
R.~{Yokota}, J.~{Pestana}, H.~{Ibeid}, D.~{Keyes}, {Fast Multipole
  Preconditioners for Sparse Matrices Arising from Elliptic Equations}, ArXiv
  e-prints\href {http://arxiv.org/abs/1308.3339} {\path{arXiv:1308.3339}}.

\bibitem{Cortes2015}
A.~M.~A. Cortes, A.~L.~G.~A. Coutinho, L.~Dalcin, V.~M. Calo, {Performance
  evaluation of block-diagonal preconditioners for the divergence-conforming
  B-spline discretization of the Stokes system}, Journal of Computational
  Science 11 (2015) 123 -- 136.
\newblock \href
  {http://dx.doi.org/http://dx.doi.org/10.1016/j.jocs.2015.01.005}
  {\path{doi:http://dx.doi.org/10.1016/j.jocs.2015.01.005}}.

\bibitem{nsch15a}
P.~Vignal, A.~F. Sarmiento, A.~M.~A. C{\^o}rtes, L.~Dalcin, V.~M. Calo,
  Coupling {N}avier--{S}tokes and {C}ahn--{H}illiard equations in a
  two-dimensional annular flow configuration, Procedia Computer Science 51
  (2015) 934--943.

\bibitem{nsch15b}
L.~F.~R. {Espath}, A.~F. {Sarmiento}, P~{Vignal}, B.~O.~N. {Varga}, A.~M.~A.
  {Cortes}, L~{Dalcin}, V.~M. {Calo}, {Energy Exchange Analysis in Droplet
  Dynamics via the Navier--Stokes--Cahn--Hilliard Model}, ArXiv e-prints\href
  {http://arxiv.org/abs/1512.02249} {\path{arXiv:1512.02249}}.

\bibitem{LesPiegl}
Les~A. Piegl, W.~Tiller, The {NURBS} book, 2nd Edition, Springer, 1996.

\bibitem{Gonz}
Oscar Gonzalez, Andrew~M. Stuart, A First Course in Continuum Mechanics,
  Cambridge University Press, 2008.

\bibitem{WeakImpos1}
Y.~Bazilevs, C.~Michler, V.~M. Calo, T.~J.~R. Hughes, Weak {Dirichlet} boundary
  conditions for wall-bounded turbulent flows, Computer Methods in Applied
  Mechanics and Engineering 196~(49-52) (2007) 4853--4862.
\newblock \href {http://dx.doi.org/10.1016/j.cma.2007.06.026}
  {\path{doi:10.1016/j.cma.2007.06.026}}.

\bibitem{WeakImpos2}
Y.~Bazilevs, C.~Michler, V.~M. Calo, T.~J.~R Hughes, Isogeometric variational
  multiscale modeling of wall-bounded turbulent flows with weakly enforced
  boundary conditions on unstretched meshes, Computer Methods in Applied
  Mechanics and Engineering 199~(13-16) (2010) 780--790.
\newblock \href {http://dx.doi.org/10.1016/j.cma.2008.11.020}
  {\path{doi:10.1016/j.cma.2008.11.020}}.

\bibitem{petsc}
Satish Balay, Shrirang Abhyankar, Mark~F. Adams, Jed Brown, Peter Brune, Kris
  Buschelman, Lisandro Dalcin, Victor Eijkhout, William~D. Gropp, Dinesh
  Kaushik, Matthew~G. Knepley, Lois~Curfman McInnes, Karl Rupp, Barry~F. Smith,
  Stefano Zampini, Hong Zhang.
\newblock \href{http://www.mcs.anl.gov/petsc}{{PETS}c {W}eb page} [online]
  (2015).
\newline\urlprefix\url{http://www.mcs.anl.gov/petsc}

\bibitem{Collocation}
F.~Auricchio, L.~Beirão Da~Veiga, T.~J.~R. Hughes, A.~Reali, G.~Sangalli,
  Isogeometric collocation methods, Mathematical Models and Methods in Applied
  Sciences 20~(11) (2010) 2075--2107.
\newblock \href {http://dx.doi.org/10.1142/S0218202510004878}
  {\path{doi:10.1142/S0218202510004878}}.

\bibitem{Botella}
O.~Botella, R.~Peyret, Benchmark spectral results on the lid-driven cavity
  flow, Computers \& Fluids 27~(4) (1998) 421 -- 433.
\newblock \href
  {http://dx.doi.org/http://dx.doi.org/10.1016/S0045-7930(98)00002-4}
  {\path{doi:http://dx.doi.org/10.1016/S0045-7930(98)00002-4}}.

\bibitem{Ghia}
U~Ghia, K.N Ghia, C.T Shin, {High-Re solutions for incompressible flow using
  the Navier-Stokes equations and a multigrid method }, Journal of
  Computational Physics 48~(3) (1982) 387 -- 411.
\newblock \href
  {http://dx.doi.org/http://dx.doi.org/10.1016/0021-9991(82)90058-4}
  {\path{doi:http://dx.doi.org/10.1016/0021-9991(82)90058-4}}.

\bibitem{Lo1}
D.~C. Lo, D.~L. Young, K.~Murugesan, {An accurate numerical solution algorithm
  for 3D velocity--vorticity Navier--Stokes equations by the DQ method},
  Communications in Numerical Methods in Engineering 22~(3) (2006) 235--250.
\newblock \href {http://dx.doi.org/10.1002/cnm.817}
  {\path{doi:10.1002/cnm.817}}.

\bibitem{Lo2}
D.~C. Lo, K.~Murugesan, D.~L. Young, {Numerical solution of three-dimensional
  velocity--vorticity Navier--Stokes equations by finite difference method},
  International Journal for Numerical Methods in Fluids 47~(12) (2005)
  1469--1487.
\newblock \href {http://dx.doi.org/10.1002/fld.822}
  {\path{doi:10.1002/fld.822}}.

\bibitem{Wong}
K.~L. Wong, A.~J. Baker, {A 3D incompressible Navier--Stokes
  velocity--vorticity weak form finite element algorithm}, International
  Journal for Numerical Methods in Fluids 38~(2) (2002) 99--123.
\newblock \href {http://dx.doi.org/10.1002/fld.204}
  {\path{doi:10.1002/fld.204}}.

\end{thebibliography}


\end{document}